\documentclass[twoside,11pt]{article}

\usepackage[english]{babel}
\usepackage{hhline}
\usepackage{bbm}
\usepackage{amsmath}
\usepackage{mathtools}
\usepackage{caption}
\usepackage[framemethod=TikZ]{mdframed} %
\usepackage[linesnumbered,ruled]{algorithm2e}
\usepackage{url}
\usepackage{enumerate}
\graphicspath {{figures/}}

\usepackage{jmlr2e}

\newtheorem{defi}{Definition}
\newtheorem{ass}{Assumption}
\newtheorem{lem}{Lemma}
\newtheorem{ex}{Example}
\newtheorem{prop}{Proposition}
\newtheorem{thm}{Theorem}
\newtheorem{cor}{Corollary}

\renewcommand\epsilon{\varepsilon}
\renewcommand\subset{\subseteq}

\newcommand\N{\mathbb{N}}

\newcommand\R{\mathbb{R}}

\renewcommand\P{\mathbb{P}}
\newcommand\E{\mathbb{E}}

\newcommand{\B}[1]{\mathbf{#1}}

\DeclareMathOperator*{\argmax}{arg\,max}
\newcommand{\given}{\,|\,}

\DeclareMathOperator{\pa}{PA}
\newcommand\PA[2]{\pa^{#2}(#1)}

\newcommand\indep{\protect\mathpalette{\protect\independenT}{\perp}}
\def\independenT#1#2{\mathrel{\rlap{$#1#2$}\mkern2mu{#1#2}}}

\DeclarePairedDelimiterX{\card}[1]{\lvert}{\rvert}{#1}
\DeclarePairedDelimiterX{\norm}[1]{\lVert}{\rVert}{#1}
 
\usepackage{tikz}
\usetikzlibrary{arrows}
\usetikzlibrary{positioning}
\newdimen\arrowsize
\pgfarrowsdeclare{arcsq}{arcsq}
{
  \arrowsize=0.2pt
  \advance\arrowsize by .5\pgflinewidth
  \pgfarrowsleftextend{-4\arrowsize-.5\pgflinewidth}
  \pgfarrowsrightextend{.5\pgflinewidth}
}
{
  \arrowsize=1.5pt
  \advance\arrowsize by .5\pgflinewidth
  \pgfsetdash{}{0pt} %
  \pgfsetroundjoin   %
  \pgfsetroundcap    %
  \pgfpathmoveto{\pgfpoint{0\arrowsize}{0\arrowsize}}
  \pgfpatharc{-90}{-140}{4\arrowsize}
  \pgfusepathqstroke
  \pgfpathmoveto{\pgfpointorigin}
  \pgfpatharc{90}{140}{4\arrowsize}
  \pgfusepathqstroke
}

\tikzset{every picture/.style={line width=0.6pt}}
\tikzset{every picture/.style={outer sep=.4mm}}
\tikzstyle{graphnode} = 
   [circle,draw=black,minimum size=25pt,text centered,inner sep=0.2mm] 
\tikzstyle{observed}   =[graphnode,fill=white,text=black]
\tikzstyle{unobserved}   =[graphnode,fill=white,text=black,style=dashed]
\tikzstyle{graphnodesmall} = 
   [circle,draw=black,minimum size=20pt,text centered,inner sep=0.2mm] 
\tikzstyle{observedsmall}   =[graphnodesmall,fill=white,text=black]
\tikzstyle{unobservedsmall}   =[graphnodesmall,fill=white,text=black,style=dashed]
\tikzstyle{graphnodestiny} = 
   [circle,draw=black,minimum size=15pt,text centered,inner sep=0.2mm] 
\tikzstyle{observedtiny}   =[graphnodestiny,fill=white,text=black]
\tikzstyle{unobservedtiny}   =[graphnodestiny,fill=white,text=black,style=dashed]

\definecolor{light-gray}{gray}{0.97}
\definecolor{commentgray}{gray}{0.5}
 
\DeclareRobustCommand\sampleline[1]{%
  \tikz\draw[#1] (0,0) (0,\the\dimexpr\fontdimen22\textfont2\relax)
  -- (1em,\the\dimexpr\fontdimen22\textfont2\relax);%
}

\begin{document}

\title{Switching Regression Models and Causal Inference in the Presence of Discrete Latent Variables}
       
\author{\name Rune Christiansen\thanks{\textit{Address for correspondence}: Rune Christiansen, Department of Mathematical Sciences, Universitetsparken 5, 2100 Copenhagen, Denmark.} \email krunechristiansen@math.ku.dk \\
	   \name Jonas Peters \email jonas.peters@math.ku.dk \\
       \addr Department of Mathematical Sciences \\
       University of Copenhagen\\
       Copenhagen, Denmark}

\editor{}

\maketitle

\begin{abstract}%
Given a response $Y$ and a vector $X = (X^1, \dots, X^d)$ of $d$ predictors,
we investigate the problem of inferring direct causes of $Y$ among the vector $X$.
Models for $Y$ that use 
all of its
causal covariates as predictors enjoy the property of being invariant across different environments or interventional settings. 
Given data from such environments, this property has been exploited for causal discovery.
Here,
we extend this inference principle to situations in which some 
(discrete-valued)
direct causes of $ Y $ are unobserved.
Such cases naturally give rise to switching regression models.
We provide sufficient conditions for the existence, consistency and asymptotic normality of the MLE in linear switching regression models with Gaussian noise, and construct a test for the equality of such models.
These results
allow us to prove that the proposed causal discovery method obtains asymptotic false discovery control under mild conditions. 
We provide an algorithm, make available code, and test our method on simulated data. It is robust against model violations and outperforms state-of-the-art approaches. 
We further apply our method to a real data set, where we show that it does not only output causal predictors, but also
a process-based clustering of data
points, which could be of additional interest to practitioners.
\end{abstract}

\begin{keywords}
causal discovery, invariance, switching regression models, hidden Markov models, latent variables
\end{keywords}

\section{Introduction} \label{sec:intro}
\subsection{Causality}
In many real world applications, we are often interested in causal rather than purely statistical relations. 
In the last decades, 
seminal work by 
\citet{Imbens2015}, \citet{Spirtes2000}, and \citet{Pearl2009}
has provided a solid mathematical basis for 
formalizing 
causal questions.
They often start from a given causal model in the form of a structural causal model (SCM) or potential outcomes. 
In practice, we often do not know the underlying causal model, and the field of causal discovery aims at inferring causal models from data. 
There are several lines of work that are based on different assumptions. 
Among them are constraint-based methods \citep{Spirtes2000, Pearl2009, Maathuis2009}, score-based methods \citep{Chickering2002, 
Silander2006,
Koivisto2006,
Cussens2011},
methods based on restricted SCMs \citep{Shimizu2006, 
Mooijetal16,
Peters2017book}, and methods based on the independence of causal mechanisms \citep{Janzing2012,Steudel2010a}.
The problem of hidden variables 
has been addressed in several works \citep[e.g.,][]{
Spirtes1995,
Silva2006,
Silva2009,
Sgouritsa2013,
Claassen2013,
Ogarrio2016, 
evans2016causal, 
richardson2017nested, 
tsagris2018constraint}. 
These methods usually consider slightly different setups than our work does; e.g., they concentrate on full causal discovery (rather than estimating causal parents), and consider different model classes. 

In this work, 
instead of aiming to learn all of the data generating structure, 
we consider the subproblem of inferring the set of causal parents of a target variable 
$Y$ among a set of variables $X = (X^1, \dots, X^d)$. 
We furthermore assume that some of the causal 
predictors are unobserved. 
While in general, this is a notoriously hard problem to solve, we will 
constrain the influence of the hidden variables by assuming that they take only few different values. 
Such a model is applicable whenever the system may be in one of several unobserved states and was motivated by an example from Earth system science, see Section~\ref{sec:SIF}.
We further assume that the data are not purely observational but come from different environments.

For the case when all causal parents are observed, \citet{peters2016causal} recently proposed the method \textit{invariant causal prediction} (ICP). 
Under the assumption that the causal mechanism generating $Y$ from its causal predictors remains the same in all environments (``invariant prediction''), it is possible to obtain the following guarantee: with large probability, the inferred set is a subset of the true set of causal predictors. A concise description of the method is provided in Section~\ref{sec:principle}.

If some of the causal predictors are unobserved, the above guarantee will, in general, not hold anymore. 
Under the additional assumption of faithfulness, one can still prove that ICP infers a subset of the causal ancestors of the target $Y$.
 In many cases, however, the method of ICP infers the empty set, which is not an incorrect, but certainly an uninformative answer.
This paper extends the idea of invariant models to situations, in which relevant parts of the system are unobserved. 
In particular, we suggest a relaxation of the invariance assumption and introduce the formal framework of $h$-invariance (``hidden invariance''). If the influence of the hidden variable is not too complex, e.g., because it takes only a few discrete values, this property is restrictive enough to be exploited for causal discovery. 
The assumption of $h$-invariance gives rise to switching regression models, 
where each value of the hidden variable corresponds to a 
different regression coefficient
(we provide more details in Section~\ref{sec:srm}).
For building an invariance-based procedure, we require a test for the equality of switching regression models. 
In this paper, we provide such a test and show that it satisfies asymptotic level guarantees. This result allows us to prove that our causal discovery procedure is asymptotically correct under mild assumptions. 
In case of sequential data, we allow for the possibilities that the hidden variables follow  an i.i.d.\ structure or a hidden Markov model \citep[e.g.,][]{Zucchini2016}. 
We suggest efficient algorithms, provide code and test our method on simulated and real data.

\subsection{Switching Regression Models} \label{sec:srm}
Switching regression models
are often used to model statistical dependencies that are subject to unobserved ``regime switches'', and 
can be viewed as ordinary regression models that include interactions with a discrete hidden variable. 
Roughly speaking, each data point $(X_i,Y_i)$ is assumed to follow one of several different regression models;
a formal definition is given in Definition~\ref{def:mixreg}. 
Switching regression models
have been used in various disciplines, e.g., to model stock returns \citep{sander2018market}, energy prices \citep{langrock2017markov} or the propagation rate of plant infections \citep{turner2000estimating}. Statistical inference in switching regression models is a challenging problem for several reasons: switching regression models are non-identifiable (permuting mixture components does not change the modeled conditional distribution),
and their likelihood function is unbounded 
(one may consider one of the regression models containing a single point with noise variance shrinking toward zero)
and non-convex. 
In this paper, we circumvent the problem of an unbounded likelihood function by
imposing parameter constraints on the error variances of the mixture components 
\citep[e.g.,][]{hathaway1985constrained, goldfeld1973estimation}.
We then
construct a test for the equality of switching regression models by evaluating
the joint overlap of the Fisher confidence regions (based on the maximum likelihood estimator) of the respective parameter vectors of the different models. 
We establish an asymptotic level guarantee for this test by providing sufficient conditions for (i) the existence, (ii) the consistency and (iii) the asymptotic normality of the maximum likelihood estimator. 
To the best of our knowledge, each of these three results is novel and may be of interest in itself.
We further discuss two ways of numerically optimizing the likelihood function.

Without parameter constraints, the likelihood function is unbounded and
global maximum likelihood estimation is an ill-posed problem \citep[e.g.,][]{DeVeaux1989}. 
Some analysis
has therefore been done on 
using 
local maxima of the likelihood function instead. \citet{kiefer1978discrete} show that there exists a sequence of roots of the likelihood equations that yield a consistent estimator, but provide no information on which root, in case there is more than one, is consistent. 
Another popular approach is to impose parameter constraints on the error variances of the mixture components. 
In the case of ordinary, univariate Gaussian mixture models, \citet{hathaway1985constrained} formulate such a constrained optimization problem and prove the existence of a global optimum.
In this paper, we present a similar result for switching regression models. 
The proof of \citet{hathaway1985constrained} uses the fact that the maximum likelihood estimates of all mean parameters are bounded by the smallest and the largest observation. 
This reasoning cannot be applied to the regression coefficients in switching regression models and therefore requires a modified argument.
We also provide sufficient conditions for the consistency and the asymptotic normality (both up to label permutations) of the proposed constrained maximum likelihood estimator. Our proofs 
are based on the proofs provided by
\citet{bickel1998asymptotic} and \citet{jensen1999asymptotic}, who show similar results for the maximum likelihood estimator in hidden Markov models with finite state space. Together, (ii) and (iii) prove the asymptotic coverage of Fisher confidence regions and ensure the asymptotic level guarantee of our proposed test.

Readers mainly interested in inference in switching regression models, may want to skip directly to Section~\ref{sec:inferenceSR}.
Additionally, Sections~\ref{sec:tests}~and~\ref{sec:ci}
contain
our proposed test for the equality of switching regression models that is available as the function \texttt{test.equality.sr} in our code package.

\subsection{The Principle of Invariant Causal Prediction} \label{sec:principle}
This section follows the presentation provided by \citet{pfister2017invariant}.
Suppose that we observe several instances 
$(Y_1, X_1), \ldots, (Y_n, X_n)$ 
of a response or target variable $Y \in \R$ and covariates $X \in \R^{1 \times d}$. 
We assume that the instances 
stem from different environments $e \subseteq \{1, \dots, n \}$, and use $\mathcal{E}$ to denote the collection of these, i.e., $\dot\bigcup_{e \in \mathcal{E}} e = \{1, \ldots, n\}$.
These environments can, for example, correspond to different physical or geographical settings in which the system is embedded, or controlled experimental designs in which some of the variables have been intervened on. 
The crucial assumption is then that there exists a subset $ S^* \subseteq \{1, \dots, d\}$ of variables from $ X $ that yield a predictive model for $ Y $ that is invariant across all environments.

More formally, one assumes the existence of a set $S^* \subset \{1, \dots, d\}$, such that
for all $x$ and all $1 \leq s, t \leq n$, we have 
\begin{equation} \label{eq:Sstar}
Y_s \given (X^{S^*}_s = x) \stackrel{d}{=} Y_t \given (X^{S^*}_t = x),
\end{equation}
where $X_t^{S^*}$ denotes the 
covariates in $S^*$ at instance $t$. For simplicity, the reader may think about~\eqref{eq:Sstar} in terms of conditional densities.
Also, the reader might benefit from 
thinking about the set $S^*$ in the context of causality, which is why we will below refer to the set $S^*$ as the set of (observable) direct causes of the target variable.
If, for example, data come from a structural causal model
(which we formally define in Appendix~\ref{app:SCM}), 
and different interventional settings,
a sufficient condition for \eqref{eq:Sstar} to hold is 
that the structural assignment for $Y$ remains the same across all observations, i.e., 
there are no interventions occurring directly on $Y$.
In Section~\ref{sec:causality}, 
we will discuss the relationship to causality in more detail.
Formally, however, this paper does not rely on the definition of the term ``direct causes''. 

Since each instance is only observed once, it is usually hard to test whether 
Equation~\eqref{eq:Sstar} holds. 
We therefore make use of the environments.
Given a set $S \subset \{1, \dots, d\}$, we implicitly assume that for every $e \in \mathcal{E}$, the conditional distribution $P_{Y_t \vert X_t^S}$\footnote{We use $P_{Y_t \vert X^S_t}$ as shorthand notation for the family $\left( P_{Y_t \vert (X^S_t = x)} \right)_{x}$ of conditional distributions.} is the same for all $t \in e$, say $P^e_{Y \vert X^S}$, and check whether 
for all $e,f \in \mathcal{E}$, we have that
\begin{equation} \label{eq:Sstar2}
P^e_{Y \vert X^S} = P^f_{Y \vert X^S}.
\end{equation}
In the population case, 
Equation~\eqref{eq:Sstar2} can be used to recover (parts of) $S^*$ from the conditional distributions $P^e_{Y \vert X^S}$: 
for each subset $S \subseteq \{1, \dots, d\}$ 
of predictors 
we check the validity of 
\eqref{eq:Sstar2}
and output
the set 
\begin{equation} \label{eq:Stilde}
\tilde S := \bigcap_{S \text{ satisfies } \eqref{eq:Sstar2}} S
\end{equation}
of variables that are necessary to obtain predictive stability. 
Under assumption \eqref{eq:Sstar}, $\tilde{S}$ only contains variables from $S^*$. For purely observational data, 
i.e.,
$(Y_t, X_t) \stackrel{d}{=} (Y_s, X_s)$
for all $s,t$,
Equation~\eqref{eq:Sstar2} is trivially satisfied for any set $S \subseteq \{1, \dots, d\}$ and thus $\tilde{S} = \emptyset$. It is the different heterogeneity patterns of the data in different environments that allow for causal discovery. 
If only a single i.i.d.\ data set is available, the method's result would not be incorrect, but it would not be informative either.
Based on a sample from $(Y_t, X_t)_{t \in e}$ for each environment, \citet{peters2016causal} propose an estimator $\hat S$ of $\tilde{S}$ that comes with a statistical guarantee: with controllable (large) probability, the estimated set $\hat S$ is contained in $S^*$. In other words, whenever the method outputs a set of predictors, they are indeed causal with high certainty.

In this paper, we consider cases in which the full set of direct causes of $Y$ is not observed. We then aim to infer the set of \textit{observable} causal variables $S^* \subseteq \{1, \dots, d\}$. Since the invariance assumption~\eqref{eq:Sstar} cannot be expected to hold in this case, the principle of invariant prediction is inapplicable. We therefore introduce the concept of $h$-invariance, a relaxed version of assumption~\eqref{eq:Sstar}. If the the latent variables are constrained to take only few values,
the $h$-invariance property can, similarly to \eqref{eq:Stilde}, be used for the inference of~$S^*$. 

\subsection{Organization of the Paper} \label{sec:organization}
The remainder 
of the paper is organized as follows. 
Section~\ref{sec:icph} explains in which sense 
the principle of invariant causal prediction 
breaks down in the presence of hidden variables and proposes an adaptation of the inference principle. It also contains hypothesis tests that are suitable for the setting with hidden variables. 
In Section~\ref{sec:inferenceSR}, we establish asymptotic guarantees for these tests. This section contains all of our theoretical results on the inference in switching regression models, and can be read independently of the problem of causal inference. 
In Section~\ref{sec:covgarant}, we 
combine the results of the preceding sections
into our overall causal discovery method (ICPH), provide an algorithm and prove the asymptotic false discovery control of ICPH.
The experiments on simulated data in Section~\ref{sec:experiments} support these theoretical findings. They further show that even for sample sizes that are too small for the asymptotic results to be effective, the overall method 
generally keeps
the type I error control. The method is robust against a wide range of model misspecifications and outperforms other approaches. We apply our method to a real world data set on photosynthetic activity and vegetation type.
Proofs of our theoretical results are contained in Appendix~\ref{app:proofs}.
All our code is available as an \verb|R| package at \url{https://github.com/runesen/icph}, and can be installed by \verb|devtools::install_github("runesen/icph/code")|, for example. 
Scripts reproducing all simulations can be found at the same url.

\section{Invariant Causal Prediction in the Presence of Latent Variables} \label{sec:icph}
Consider a collection $(\B{Y}, \B{X}, \B{H}) = (Y_t, X_t, H_t)_{t \in \{1, \dots, n\}}$ of 
triples of a target variable $Y_t \in \R$, observed covariates $X_t \in \R^{1 \times d}$ and some latent variables $H_t \in \R^{1 \times k}$.
For simplicity, we refer to the index $t$ as time, but we also allow for an i.i.d.\ setting; see Section~\ref{sec:time_dep} for details. 
When referring to properties of the data that hold true for all $ t $, we sometimes omit the index altogether.

In analogy to Section~\ref{sec:principle}, we start by assuming the existence of an invariant predictive model for $ Y $, but do not require all relevant variables to be observed. That is, we assume the existence of a set $S^* \subseteq \{1, \dots, d \}$ and a subvector $H^*$ of $H$ such that the conditional distribution of $Y_t \given (X^{S^*}_t, H^*_t)$ is the same for all time points $ t $. Based on the observed data $(\B{Y}, \B{X})$, we then aim to infer the set $S^*$.

Section~\ref{sec:violation} shows why the original version of invariant causal prediction is inapplicable. In Sections~\ref{sec:hiddeninvprop} and~\ref{sec:infh_inv} we introduce the formal concept of $h$-invariance and present an adapted version of the inference principle discussed in Section~\ref{sec:principle}. 
In Sections~\ref{sec:tests}~and~\ref{sec:ci} we then present tests for $h$-invariance of sets $S \subseteq \{1, \dots, d\}$, which are needed for the construction of an empirical estimator $\hat S $ of $S^*$. A causal interpretation of the $h$-invariance property is given in Section~\ref{sec:causality}.

\subsection{Latent Variables and Violation of Invariance} \label{sec:violation}
The inference principle
described in Section~\ref{sec:principle} relies on the invariance assumption~\eqref{eq:Sstar}. The following example shows that if some of the invariant predictors of $Y$ are unobserved, we cannot expect this assumption to hold.
The principle of ordinary invariant causal prediction is therefore inapplicable. 
\begin{ex}[Violation of invariance assumption due to latent variables] \label{ex:icp_fail}
We consider a linear model for the data $(Y_t, X^1_t, X^2_t, H_t^*)_{t \in \{1, \dots, n\}} \in \R^{n \times 4}$. Assume there exist i.i.d.\ zero-mean noise variables $\varepsilon_1, \dots, \varepsilon_n$ such that for all $t$, $(X^1_t, H^*_t, \varepsilon_t)$ are jointly independent and
\begin{equation*}
Y_t = X^1_t + H^*_t + \epsilon_t.
\end{equation*}
Assume furthermore that the distribution of the latent variable $H_t^*$ changes over time, 
say $\E[H_{r}^*] \not = \E[H_{s}^*]$ for some $r,s$. Then, with $ S^* := \{1\} $, the conditional distribution $P_{Y_t \vert (X_t^{S^*}, H_t^*)}$ is time-homogeneous, but 
\begin{equation*}
\E[Y_{r} \vert X_{r}^{S^*} = x] = x+ \E[H^*_{r}] \not = x + \E[H^*_{s}] = \E[Y_{s} \vert X_{s}^{S^*} = x],
\end{equation*}
which shows 
that
$P_{Y_t \vert X_t^{S^*}}$ 
is not time-homogeneous,
i.e., $S^*$ 
does not satisfy \eqref{eq:Sstar}. 
\end{ex}
The above example shows that in the presence of hidden variables,
assumption \eqref{eq:Sstar}
may be too strong.
The distribution in the above example, however, allows for a different invariance. For all $t,s$ and all $x,h$ we have that\footnote{In the remainder of this work, we implicitly assume that for every $t$, $(Y_t, X_t, H_t)$ is abs.\ continuous w.r.t.\ a product measure. This ensures the existence of densities $f_t(y,x,h)$ for $(Y_t, X_t, H_t)$. The marginal density $f_t(x, h)$ can be chosen strictly positive on the support of $(X_t, H_t)$ and thereby defines a set of conditional distributions $\{Y_t \given (X_t = x, H_t=h)\}_{(x,h) \in \text{supp}((X_t, H_t))}$ via the conditional densities $f_t(y \given x,h) = f_t(y,x,h) / f_t(x,h)$. Strictly speaking, we therefore assume that the conditional distributions \textit{can be chosen} s.t.\ \eqref{eq:SstarH} holds for all $(x,h) \in \text{supp}((X_t^{S^*}, H^*_t)) \cap \text{supp}((X_s^{S^*}, H^*_s))$.} 
\begin{equation} \label{eq:SstarH}
Y_t \given (X_t^{S^*} = x, H^*_t = h) \stackrel{d}{=} Y_s \given (X_s^{S^*} = x, H^*_s = h).
\end{equation}
Ideally, we would like to directly exploit this property for the inference of $S^*$. Given a candidate set $S \subseteq \{1, \dots, d\}$, we need to check if there exist $H_1^*, \dots, H_n^*$ such that \eqref{eq:SstarH} holds true for $S^* = S$. 
Similarly to \eqref{eq:Stilde}, the idea is then to output the intersection of all sets for which this is the case. Without further restrictions on the influence of the latent variables, however, the result will always be the empty set.
\begin{prop}[Necessity of constraining the influence of $H^*$] \label{prop:restrictions}
Let %
$S \subset \{1, \dots, d\}$ be an arbitrary subset of the predictors $X_t$. Then, there exist variables $H_1, \dots, H_n$ such that \eqref{eq:SstarH} is satisfied for $S^* = S$ and $(H_t^*)_{t \in \{1, \dots, n\}} = (H_t)_{t \in \{1, \dots, n\}}$. 
\end{prop}
The proof is immediate by 
choosing latent variables with non-overlapping support (e.g., such that 
for all $t$,
$P(H_t = t)=1$).
Proposition~\ref{prop:restrictions} shows that without constraining the influence of $H^*$, \eqref{eq:SstarH} cannot be used to identify $S^*$.
Identifiability improves, however, 
for univariate, discrete latent variables $H^* \in \{1, \dots, \ell \}$ with relatively few states $\ell \geq 2$. Equation~\eqref{eq:SstarH} then translates into the following assumption on the observed conditional distributions $P_{Y_t \given X_t^{S^*}}$: for all $t,x$ it holds that
\begin{equation} \label{eq:Pxj}
P_{Y_t \vert (X_t^{S^*} = x)} = \sum_{j=1}^\ell \lambda^j_{xt} P^j_x ,
\end{equation}
for some $\lambda^1_{xt}, \dots, \lambda^\ell_{xt} \in (0,1)$ with $\sum_{j=1}^\ell \lambda^j_{xt}=1$ and distributions $P^1_x, \dots, P^\ell_x$ that do not depend on $t$. 
This fact can be
seen by expressing the conditional density of $P_{Y_t \vert (X_t^{S^*} = x)}$ as $f_t(y \given x) = \int f_t(y \given x,h) f_t(h \given x) dh$. By \eqref{eq:SstarH}, $f_t(y \given x,h)$ does not depend on $t$. Property \eqref{eq:Pxj} then follows by taking $\lambda^j_{xt} = P(H^*_t=j \given X^{S^*}_t = x)$ and letting $P^j_x$ denote the distribution of $Y_1 \given (X_1^{S^*} = x, H^*_1 = j)$.

The conditional distributions of $Y_t \given (X^{S^*}_t = x)$ are thus assumed to follow mixtures of $\ell$ distributions, each of which remains invariant across time. The mixing proportions $\lambda_{xt}$ may vary over time. In the following subsection, we translate property~\eqref{eq:Pxj} into the framework of mixtures of linear regressions with Gaussian noise. The invariance assumption on $P^1_x, \dots P_x^\ell$ then corresponds to time-homogeneity of the regression parameters of all mixture components.

\subsection{Hidden Invariance Property} \label{sec:hiddeninvprop}
As motivated by Proposition~\ref{prop:restrictions}, we will from now on assume that $ H^*$ only takes a small number of different values. We now formalize the dependence of $Y$ on $(X^{S^*}, H^*)$ by a parametric function class. 
We purposely refrain from modeling the dependence between observations of different time points, and come back to that topic in Section~\ref{sec:time_dep}. Since the inference principle described in Section~\ref{sec:principle} requires us to evaluate \eqref{eq:Pxj} 
for different candidate sets $S$, we state the following definition in terms of a general $p$-dimensional vector $X$ (which will later play the role of the subvectors $X^S$, see Definition~\ref{def:h_inv}).
\begin{defi}[Switching regression] \label{def:mixreg}
Let $X$ be a $p$-dimensional random vector, $\ell \in \N$ and $\lambda \in (0,1)^\ell$ with $\sum_{j=1}^\ell \lambda_j = 1$. Let furthermore $\Theta$ be a matrix of dimension $(p+2) \times \ell$ with columns $\Theta_{\cdot j} = (\mu_j, \beta_j, \sigma_j^2) \in \R \times \R^p  \times \R_{>0}$, for $j \in \{1, \dots, \ell\}$. The joint distribution $P$ of $(Y, X) \in \R^{(1+p)}$ is said to follow a \emph{switching regression} of degree $\ell$ with parameters $(\Theta, \lambda)$, if there exist $H \sim \emph{Multinomial}(1, \lambda)$ and $\epsilon_j \sim \mathcal{N}(0, \sigma_j^2), j \in \{1, \dots, \ell\}$, with $(\epsilon_1, \dots, \epsilon_\ell) \indep X$, such that
\begin{equation*}
Y =  \sum_{j=1}^\ell (\mu_j + X \beta_j +  \epsilon_j) \mathbbm{1}_{\{H=j\}},
\end{equation*}
where $\mathbbm{1}_{\{H=j\}}$ denotes the indicator function for the event $H=j$. 
\end{defi}
A few remarks are in place. First, we will as of now let $\ell \geq 2$ be fixed. The reader is encouraged to think of $\ell = 2$, which is also the case to be covered in most examples and experiments. (Non-binary latent variables are considered in Appendix~\ref{app:non_binary}.) Second, it will be convenient to parametrize the matrix $\Theta$ by a map $\theta \mapsto \B{\Theta}(\theta)$, $\theta \in \mathcal{T}$, where $\mathcal{T}$ is a subset of a Euclidean space. This allows for a joint treatment of different types of parameter contraints such as requiring all intercepts or all variances to be equal. We will use $\mathcal{SR}_{\B{\Theta}}(\theta, \lambda \given X)$ (``Switching Regression'') to denote the distribution $P$ over $(Y,X)$ satisfying Definition~\ref{def:mixreg} with parameters $(\B{\Theta}(\theta), \lambda)$, although we will often omit the implicit dependence on $\B{\Theta}$ and simply write $\mathcal{SR}(\theta, \lambda \given X)$. For now, the reader may think of $(\B{\Theta}, \mathcal{T})$ as the unconstrained parametrization, where $\mathcal{T} = (\R \times \R^p  \times \R_{>0})^\ell$ and where $\B{\Theta}$ consists of the coordinate projections $\B{\Theta}_{ij}(\theta) = \theta_{(j-1) (p+2) + i}$. Finally, we will for the rest of this paper disregard the intercept terms $\mu_j$  as they can be added without loss of generality by adding a constant predictor to $X$.

The following definition and assumption translate~\eqref{eq:Pxj} into the model class $\mathcal{SR}$. 
\begin{defi}[$h$-invariance] \label{def:h_inv}
A set $S \subset \{1, \dots, d\}$ is called \emph{$h$-invariant} w.r.t.\ $(\B{Y}, \B{X}) = (Y_t, X_t)_{t \in \{1, \dots, n\}}$ if there exist $\theta$ and $\lambda_1, \dots, \lambda_n$ such that, for all $t$, $P_{(Y_t, X_t^S)} = \mathcal{SR} (\theta, \lambda_t \given  X_t^S)$.
\end{defi}
Definition~\ref{def:h_inv} describes an invariance in the regression parameters $\theta$ and makes no restriction on the mixing proportions $\lambda_1, \dots, \lambda_n$. This allows the influence of the latent variable to change over time. 
From now on, we assume the existence of an $h$-invariant set $S^*$.
\begin{ass} \label{ass:h_inv}
There exists a set $S^* \subset \{1, \dots, d \}$ which is $h$-invariant w.r.t.~$(\mathbf{Y}, \mathbf{X})$.
\end{ass}
This assumption is at the very core of the proposed methodology, with the unknown $h$-invariant set $S^*$ as inferential target. 
In Section~\ref{sec:causality} 
we show that
if the data $(\B{Y}, \B{X}, \B{H})$ are generated by different interventions in an SCM
(see Appendix~\ref{app:SCM}),
in which the variable $H^* \in \{1, \dots, \ell\}$ acts on $Y$, Assumption~\ref{ass:h_inv} is 
satisfied by the set $S^* = \PA{Y}{0}$ of observable parents of $Y$.
Here, interventions are 
allowed to act on the latent variables,~and~thus indirectly on the target $Y$. For illustrations of the $h$-invariance property, see Figures~\ref{fig:mixreg}~and~\ref{fig:h_inv}.

\subsection{Relation to Causality} \label{sec:causality}
Assumption~\ref{ass:h_inv} is formulated without the notion of causality. 
The following proposition shows that if the data $(\B{Y}, \B{X}, \B{H})$ do come from an SCM, 
the set $S^*$ may be thought of as the set of observable parents of $Y$. 
\begin{prop}[Causal interpretation of $S^*$] \label{prop:causalS}
Consider an SCM over the system of variables $(Y_t, X_t, H^*_t)_{t \in \{1, \dots, n\}}$, where for every $t$, 
$(Y_t, X_t, H_t^*) \in \R^{1} \times \R^d \times \{1, \dots, \ell\}$. Assume that the structural assignment of $Y$ is fixed across time, and for every $t \in \{1, \dots, n\}$ given by
\begin{equation*}
Y_t := f(X_t^{\PA{Y}{0}}, H^*_t, N_t),
\end{equation*}
where $(N_t)_{t \in \{1, \dots, n\}}$ are i.i.d.\ noise variables. Here, $\PA{Y}{0} \subseteq \{1, \dots, d\}$ denotes the set of parents of $Y_t$ among $(X^1_t, \dots, X^d_t)$. The structural assignments for the remaining variables $X^1, \dots, X^d, H^*$ are allowed to change between different time points. 
Then, property~\eqref{eq:SstarH} is satisfied for $S^* = \PA{Y}{0}$. If furthermore the assignment $f(\cdot , h, \cdot)$ is linear for all $h \in \{1, \dots, \ell\}$ and the noise variables $N_t$ are normally distributed, 
then, Assumption~\ref{ass:h_inv} is satisfied for $S^* = \PA{Y}{0}$. 
That is, the set of observable parents of $Y$ is \emph{h}-invariant with respect to $(\B{Y}, \B{X}) = (Y_t, X_t)_{t \in \{1, \dots, n\}}$.
\end{prop}
From a causal perspective, Proposition~\ref{prop:causalS} informs us about the behavior of $P_{Y \vert (X^{S^*} = x)}$ under interventions in the data generating process. The set $S^* = \PA{Y}{0}$ will be $h$-invariant under any type of intervention that does not occur directly on the target variable (except through the latent variable $H^*$). 
The following example demonstrates the $h$-invariance property for an SCM in which 
the assignments of some of the variables change between every time point. 

\begin{figure}
\hspace{1mm}
\begin{minipage}{0.4\linewidth}
\vspace*{-3.5mm}
\begin{mdframed}[roundcorner=10pt, innertopmargin = 17.5pt, innerbottommargin = 17.5pt]
\begin{center}
\begin{tikzpicture}[xscale=1.2, yscale=1.2, shorten >=1pt, shorten <=1pt]
  \draw(1,2) node(e) [observedsmall, draw = red] {\textcolor{red}{$E$}};
  \draw (2,1) node(x3) [observedsmall] {$X^3$};
  \draw (3.2,1) node(x2) [observedsmall] {$X^2$};
  \draw (4.4,1) node(y) [observedsmall, draw = green!60!black] {\textcolor{green!60!black}{$Y$}};
  \draw (3.2,0) node(x1) [observedsmall] {$X^1$};
  \draw (3.2,2) node(h) [observedsmall, dashed] {$H^*$};

  \draw[-arcsq, red] (e) -- (h);
  \draw[-arcsq, red, bend right] (e) to [bend left = 10] (x2);
  \draw[-arcsq, red] (e) to [bend right=35]  (x1);
  \draw[-arcsq, dashed] (h) -- (x3);
  \draw[-arcsq, dashed] (h) -- (y);
  \draw[-arcsq, dashed] (h) -- (x2);
  \draw[-arcsq, green!60!black, thick] (x2) -- (y);
  \draw[-arcsq] (x3) -- (x2);
  \draw[-arcsq] (x3) -- (x1);
  \draw[-arcsq] (y) -- (x1);
\end{tikzpicture}
\end{center}
\end{mdframed}
\end{minipage}
\begin{minipage}{0.6\linewidth}
\begin{center}
\includegraphics[width=.9\linewidth]{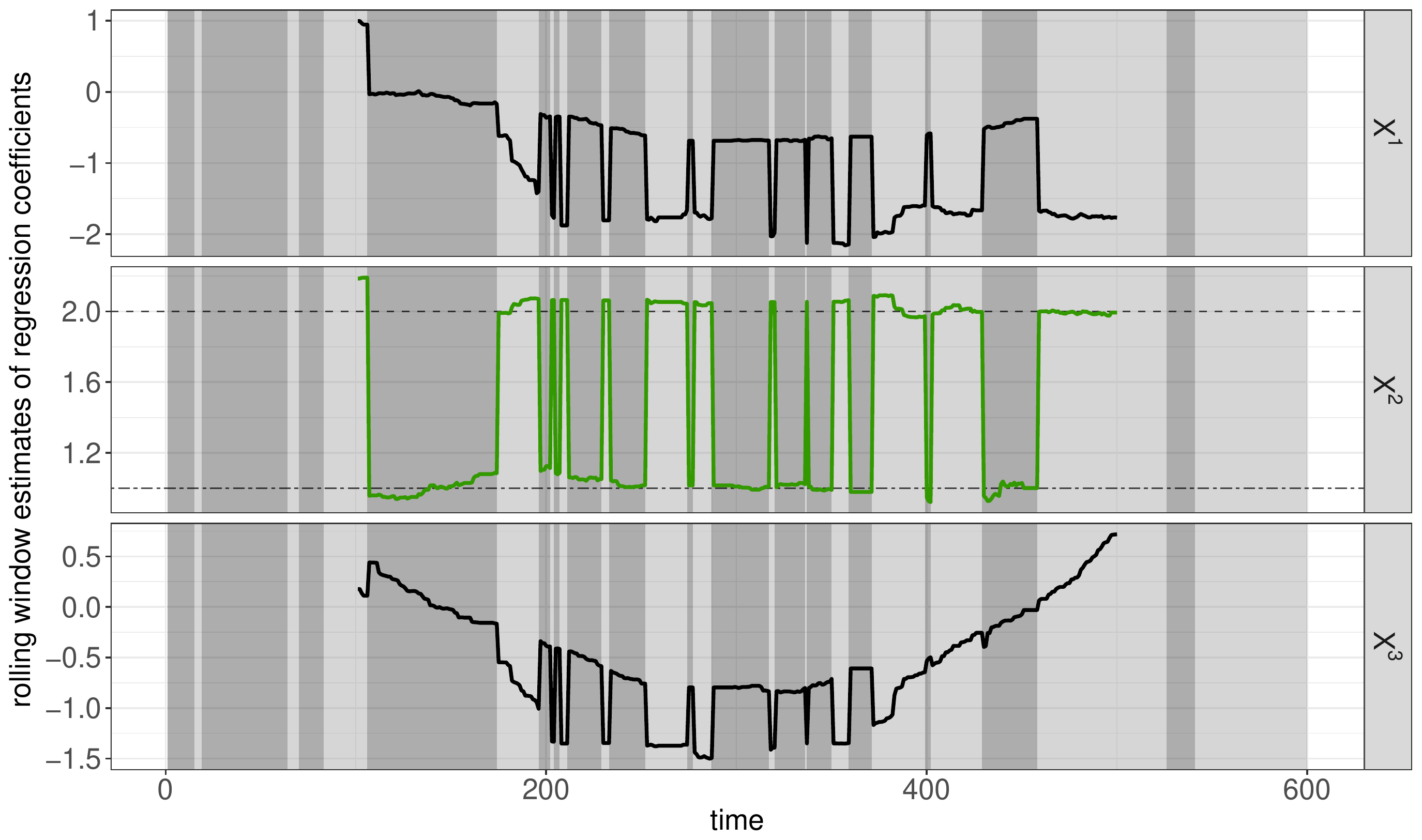}
\end{center}
\end{minipage}
\caption{
An illustration
of the $h$-invariance property based on simulated data from the SCM in Example~\ref{ex:h_inv}. The causal graph (left) and rolling window estimates of regression coefficients in the linear interaction model for the conditional distribution of $ Y $ given $ (X^1, H^*) $, $ (X^2, H^*) $ and $ (X^3, H^*) $, respectively (right).
Within both regimes $H_t^* = 1$ and $H_t^* = 2$ (corresponding to different background colors in the plot), the regression coefficient for $X^2$ (green) is time-homogeneous, and the set $S^* = \{2\}$ is therefore $h$-invariant with respect to $(\B{Y}, \B{X})$. Due to heterogeneity in the data (``the variable $E$ acts on $X^1$, $X^2$ and $H^*$''), neither of the sets $\{1\}$ or $\{3\}$ satisfy $h$-invariance. 
In practice, we test for $h$-invariance using environments, rather than rolling windows, see Section~\ref{sec:tests}.
}
\label{fig:mixreg}
\end{figure}
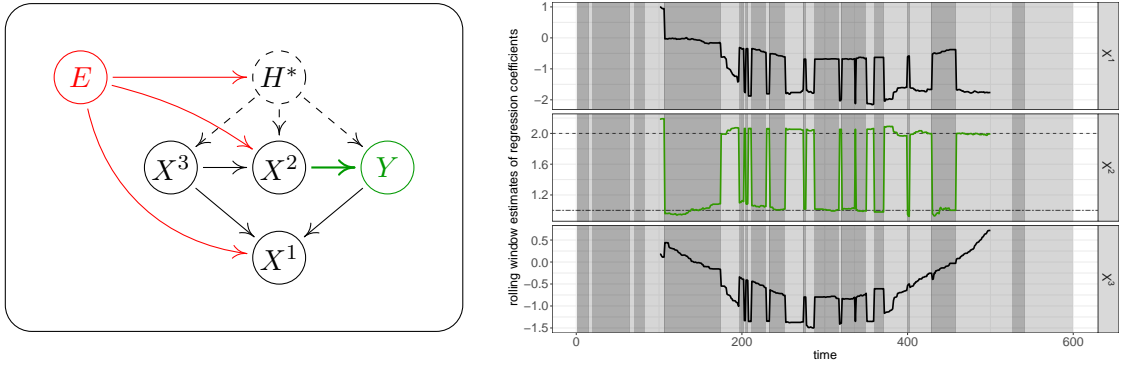

\begin{ex} \label{ex:h_inv}
Consider an SCM over the system of variables $(Y_t, X_t, H_t^*)_{t \in \{1, \dots, n\}}$, where for every $t$, 
the causal graph over 
$(Y_t, X_t, H_t^*) \in \R^{1} \times \R^3 \times \{1,2\}$
is given as in Figure~\ref{fig:mixreg}. 
The node $E$ denotes the ``environment variable'' and the outgoing edges from $E$ to $X^1$, $X^2$ and $H^*$ 
indicate that the structural assignments of these variables change throughout time. The structural 
assignment of $Y$ is fixed across time, and for every $t \in \{1, \dots, n\}$ given by 
\begin{align*}
Y_t &:=  (1+ X^2_t + 0.5 N_t)1_{\{H_t^*=1\}} + (1 + 2 X^2_t + 0.7 N_t) 1_{\{H_t^*=2\}},
\end{align*}
where $(N_t)_{t \in \{1, \dots, n\}}$ are i.i.d.\ standard Gaussian noise variables. 
Then, by Proposition~\ref{prop:causalS}, the set $S^* = \{2\}$ of observable parents of $Y$ is $h$-invariant w.r.t.\ $(\B{Y}, \B{X})$, see Figure~\ref{fig:mixreg}. 
\end{ex}

\subsection{Inference of the $h$-Invariant Set} \label{sec:infh_inv}
In general, Definition~\ref{def:h_inv} does not define a unique set of predictors. 
In analogy to \citet{peters2016causal}, we thus  propose to output the intersection of all $h$-invariant sets. We define 
\begin{equation} \label{eq:H0}
H_{0,S} : S \text{ is } h \text{-invariant with respect to } (\mathbf{Y}, \mathbf{X}), \text{ and}
\end{equation}
\begin{equation} \label{eq:Stilde2}
\tilde S := \bigcap_{S: \, H_{0,S} \text{ true }} S,
\end{equation}
where $S$ runs over subsets ${S \subset \{1, \dots, d\}}$.
In \eqref{eq:Stilde2}, we define the intersection over an empty index set as the empty set. 
In practice, we are given a sample from $(\B{Y}, \B{X})$, and our goal is to estimate $\tilde{S}$. 
Given a family of tests $(\varphi_S)_{S \subset \{1, \dots, d\}}$ of the hypotheses $(H_{0,S})_{S \subset \{1, \dots, d\}}$, we therefore define an empirical version of \eqref{eq:Stilde2} by 
\begin{equation} \label{eq:Shat}
\hat S := \bigcap_{S: \, \varphi_S \text{ accepts } H_{0,S}} S.
\end{equation}
Using that $\{\varphi_{S^*} \text{ accepts } H_{0,S^*}\} \subset \{ \hat{S} \subset S^* \}$, we immediately obtain the following important coverage property. 
\begin{prop}[Coverage property] \label{prop:icp_level}
Under Assumption 1 and given a family of tests $(\varphi_S)_{S \subset \{1, \dots, d\}}$ of $(H_{0,S})_{S \subset \{1, \dots, d\}}$ that are all valid at level $\alpha$, we have that $\P(\hat S \subset S^*) \geq 1 - \alpha$. In words, the (setwise) false discovery rate of \eqref{eq:Shat} is controlled at level $\alpha$.
\end{prop}
The set $S^*$ in Proposition~\ref{prop:icp_level} may not be uniquely determined by the $h$-invariance property. 
But since 
our output is
the \textit{intersection} \eqref{eq:Shat}
of all $h$-invariant sets, this ambiguity does no harm---the coverage guarantee for the inclusion $\hat S \subseteq S^*$ will be valid for \textit{any} choice of $h$-invariant set $S^*$.
The key challenge that remains is the construction of the tests $(\varphi_S)_{S \subseteq \{1, \dots, d\}}$, which we will discuss in Section~\ref{sec:tests}. 

\subsubsection{Tests for non-causality of individual predictors} \label{sec:pvals}
Proposition~\ref{prop:icp_level} proves a level guarantee for the estimator $\hat S$. 
To obtain statements about the significance of individual predictors that could be used for a ranking of all the variables in $X$, for example,
we propose the following construction. 
Whenever at least one hypothesis $H_{0,S}$ is accepted, we define for every $j \in \{1, \dots, d \}$ a $p$-value for the hypothesis $H_{0}^j: j \not \in S^*$ of non-causality of $X^j$ by $p_j := \max \{ p\text{-value for } H_{0,S}:  j \not \in S \}$. 
When all hypotheses $H_{0,S}$, $S \subset \{1, \dots, d\}$, are rejected (corresponding to rejecting the existence of $S^*$), we set all of these $p$-values to 1. The validity of thus defined tests is ensured under the assumptions of Proposition~\ref{prop:icp_level}, and is a direct consequence of $\varphi_{S^*}$ achieving correct level $\alpha$.

\subsection{Tests for the Equality of Switching Regression Models} \label{sec:tests}
\begin{figure}
\begin{minipage}{0.475\linewidth}
\begin{center}
$\mathbf{X^1}$

\vspace{2mm}

\includegraphics[width=.98\linewidth]{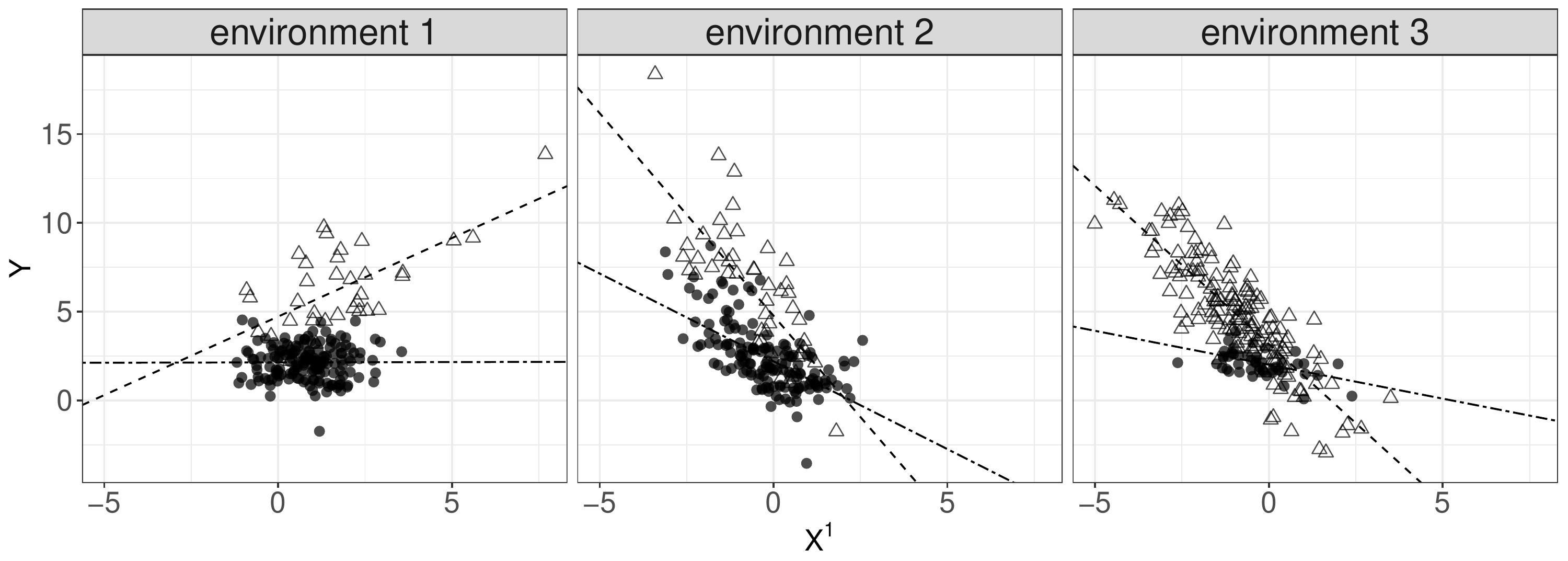}
\end{center}
\end{minipage}
\begin{minipage}{0.025\linewidth}
$ $
\end{minipage}
\begin{minipage}{0.475\linewidth}
\begin{center}
$\mathbf{X^2}$

\vspace{2mm}

\includegraphics[width=.98\linewidth]{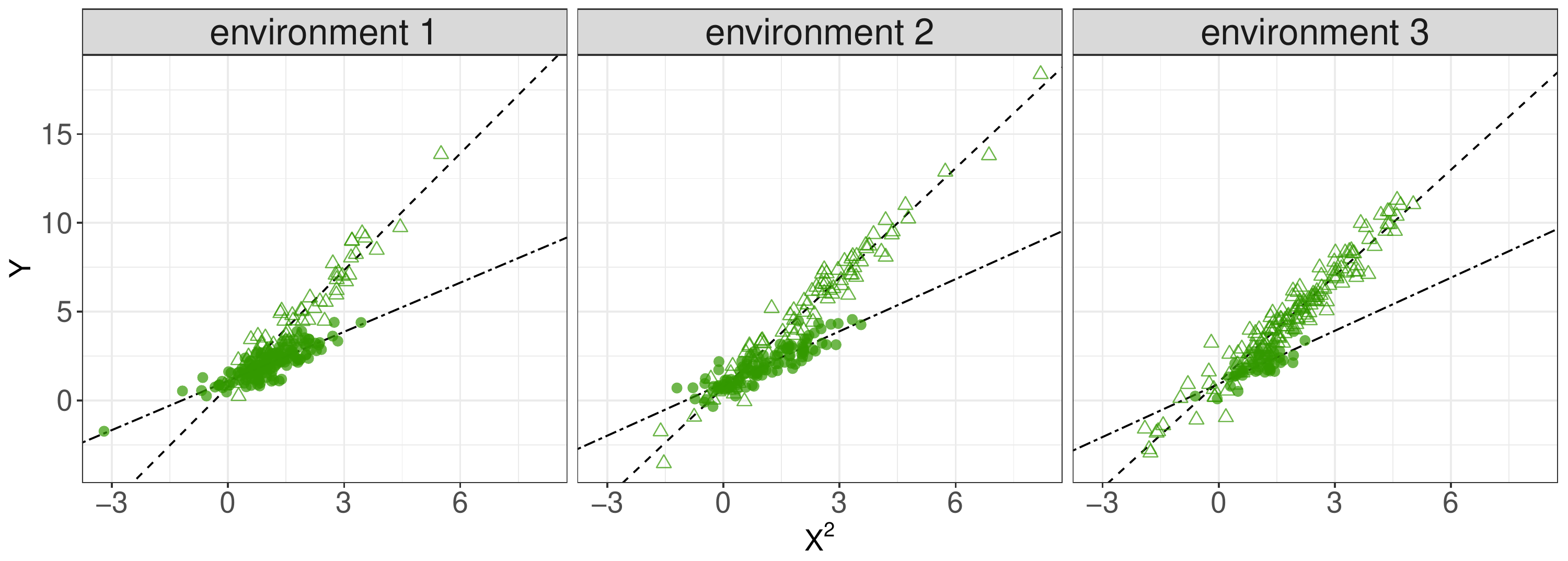}
\end{center}
\end{minipage}
\caption{Testing procedure for $H_{0,S}$, here illustrated for the sets $\{1\}$ (black; not $h$-invariant) and $\{2\}$ (green; $h$-invariant) using the same data that generated Figure~\ref{fig:mixreg}. First, we split data up into several environments, here $e_1 = \{1, \dots, 200\}$, $e_2 = \{201, \dots, 400\}$ and $e_3 = \{401, \dots, 600\}$. Then, we fit an $\mathcal{SR}$ model to each data set $(\B{Y}_e, \B{X}^S_e)$, $e \in \mathcal{E}$, separately, and evaluate whether the mixture components remain invariant across all environments. 
For illustration purposes, we indicate model fits by dashed lines, 
and assign
points to the most likely hidden state ($\bullet: \hat H^*_t = 1$, $\vartriangle: \hat H_t^* = 2$).
(This explicit classification of points is not part of the proposed testing procedure.)}
\label{fig:h_inv}
\end{figure}

We will now focus on the construction of tests for the hypotheses $H_{0,S}$ that are needed to compute the empirical estimator \eqref{eq:Shat}. Let $S \subset \{1, \dots, d\}$ be fixed for the rest of this section. We will make use of the notation $\B{X}^S$ to denote the columns of $\B {X}$ with index in $S$ and $\B{Y}_e = (Y_t)_{t \in e}$ and $\B{X}_e^S = (X^S_t)_{t \in e}$ for the restrictions of $\B{Y}$ and $\B{X}^S$ to environment $e \in \mathcal{E}$. 
For notational convenience, we rewrite $H_{0,S}(\mathcal{E}) := H_{0,S}$ as follows.
\begin{align*}
H_{0,S}(\mathcal{E}): \begin{cases} 
\text{There exist } \lambda_1, \dots, \lambda_n \text{ and } 
(\theta_e)_{e \in \mathcal{E}}, 
\text{ such that, for all } e \in \mathcal{E},  \\ 
P_{(Y_t, X_t^S)} = \mathcal{SR}(\theta_e, \lambda_t \given X_t^S) \text{ if } t \in e,
\text{ and } \text{for all } e,f \in \mathcal{E},  \theta_e = \theta_f.
\end{cases}
\end{align*}
Intuitively, a test $\varphi_S = \varphi_S(\mathcal{E})$ of $H_{0,S}(\mathcal{E})$ should reject whenever the parameters $\theta_e$ and $\theta_f$
differ between at least two environments $e,f \in \mathcal{E}$. This motivates a two-step procedure:
\begin{enumerate}[(i)]
\setlength\itemsep{0em}
\item For every $e \in \mathcal{E}$, fit an $\mathcal{SR}$ model to $(\B{Y}_e, \B{X}^S_e)$ to obtain an estimate $\hat \theta_e$ with confidence intervals, see Section~\ref{sec:inferenceSR}.
\item Based on (i), test if $\theta_e = \theta_f$ for all $e,f \in \mathcal{E}$, see Section~\ref{sec:ci}.
\end{enumerate}
For (i), we use maximum likelihood estimation and construct individual confidence regions for the estimated parameters $\hat \theta_e$ using the asymptotic normality of the MLE.
For (ii), we evaluate the joint overlap of these confidence regions. 
Any other test for the equality of $\mathcal{SR}$ models can be used here, but to the best of our knowledge, we propose the first of such tests.
Figure~\ref{fig:h_inv} illustrates step (i) for the two canditate sets $\{1\}$ and $\{2\}$. Here, we would expect a test to reject the former set, while accepting the truly $h$-invariant set $S^* = \{2\}$.
A generic approach for comparing ordinary linear regression models across different environments can be based on exact resampling of the residuals \citep[e.g.,][]{pfister2017invariant}. This procedure, however, is not applicable to mixture models: 
after fitting the mixture model, 
the states $H_t$ are unobserved,
and thus, there are multiple definitions of the residual
$r^j_t = Y_t - X_t^S \hat \beta_j$, $j \in \{1, \dots, \ell\}$.

\subsection{Intersecting Confidence Regions} \label{sec:ci}
Assume $H_{0,S}(\mathcal{E})$ is true and let $\theta_0$ be the true vector of regression parameters (that is the same for all environments). If for $e \in \mathcal{E}$, $C^{\alpha}_e = C^{\alpha}_e(\B{Y}_e, \B{X}^S_e)$ are valid $(1- \alpha)$--confidence regions for $\theta_e = \theta_0$, we can obtain a $p$-value for $H_{0,S}(\mathcal{E})$ by 
considering their joint overlap. More formally, we construct the test statistic $T_S: \R^{n \times (1+\card{S})} \to [0,1 ]$ by
\begin{equation} \label{eq:T}
T_S(\B{Y}, \B{X}^S) := \max \left \lbrace \alpha \in [0,1] : \bigcap_{e \in \mathcal{E}} C^{\alpha / \card{\mathcal{E}}}_e (\B{Y}_e, \B{X}^S_e) \not = \emptyset \right \rbrace,
\end{equation}
and define a test $\varphi^{\alpha}_S$ by $ \varphi^{\alpha}_S = 1 :\Leftrightarrow T_S < \alpha$.
Due to the Bonferroni correction of the confidence regions, such a test will be conservative. 
The construction of confidence regions is discussed in the following section.

\section{Inference in Switching Regression Models} \label{sec:inferenceSR}
In this section, we discuss maximum likelihood estimation and the construction of confidence regions for the parameters in $\mathcal{SR}$ models. 
In Sections~\ref{sec:time_dep}--\ref{sec:likelihood} we present two different models for time dependencies in the data, introduce the likelihood function for $\mathcal{SR}$ models, and present two types of parameter constraints that ensure the existence of the maximum likelihood estimator. In Section~\ref{sec:fisher}--\ref{sec:labels} we construct confidence regions based on the maximum likelihood estimator, and in Section~\ref{sec:covgarant_cr} we show that these confidence regions attain the correct asymptotic coverage. As a corollary, we obtain that the test defined in~\eqref{eq:T} satisfies asymptotic type I error control. 

Let $S \subseteq \{1, \dots, d\}$ and consider a fixed environment $e$, say $e = \{1, \dots, m\}$. 
Throughout this section, we will omit all indications of $S$ and $e$ and simply write $(Y_t, X_t) \in \R^{1+p}$ for $(Y_t, X_t^S)$ and $(\B{Y}, \B{X})$ for $(\B{Y}_e, \B{X}^S_e)$.

\subsection{Time Dependence and Time Independence} \label{sec:time_dep}
Assume there exist parameters $\theta$ and $\lambda_1, \dots, \lambda_m$ such that, for all $t \in \{1, \dots, m\}$, $(Y_t, X_t) \sim \mathcal{SR}( \theta, \lambda_t \given X_t)$. 
Let $\B{H} = (H_t)_{t \in \{1, \dots, m\}} \in \{1, \dots, \ell\}^m$ be such that for every $t \in \{1, \dots, m\}$, the distributional statement in Definition~\ref{def:mixreg} holds for $(Y_t, X_t, H_t)$. We will now consider two different models for the dependence between observations of different time points:
\begin{itemize}
\setlength\itemsep{0em}
\item Independent observations (``IID''): All observations $(Y_t, X_t, H_t)$ across different time points $t = 1, \dots, m$ are jointly independent and the marginal distribution of $\B{H}$ is time-homogeneous. Furthermore, for every $t \in \{1, \dots, m\}$, the variables $X_t$ and $H_t$ are independent. 
\item A hidden Markov model (``HMM''): The dependence in the data is governed by a first order Markovian dependence structure on the latent variables $\B{H}$ as described in Figure~\ref{fig:HMM}. The Markov chain $\B{H}$ is initiated in its stationary distribution. Furthermore, for every $t \in \{1, \dots, m\}$, the variables $X_t$ and $H_t$ are independent. 
\end{itemize}
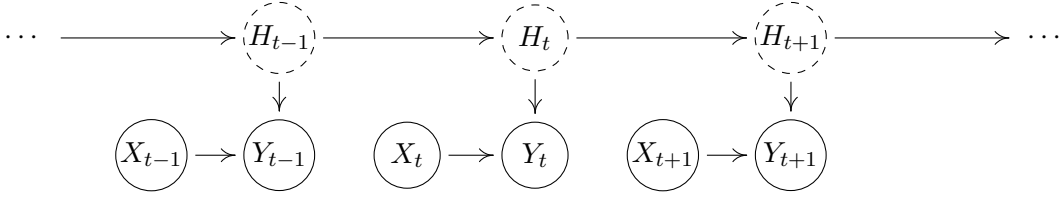
\begin{figure}
\centering
\begin{tikzpicture}[xscale=1.7, yscale=1.2, shorten >=2pt, shorten <=2pt]
  \draw (0,0) node(dots) [] {$\cdots$};
  \draw (2,0) node(h0) [unobserved] {$H_{t-1}$};
  \draw (4,0) node(h1) [unobserved] {$H_{t}$};
  \draw (6,0) node(h2) [unobserved] {$H_{t+1}$};
  \draw (8,0) node(dotss) [] {$\cdots$};
  \draw (2,-1.3) node(y0) [observed, inner sep = .5mm] {$Y_{t-1}$};
  \draw (4,-1.3) node(y1) [observed, inner sep = 1mm] {$Y_{t}$};
  \draw (6,-1.3) node(y2) [observed, inner sep = .5mm] {$Y_{t+1}$};
  \draw (1,-1.3) node(x0) [observed] {$X_{t-1}$};
  \draw (3,-1.3) node(x1) [observed] {$X_{t}$};
  \draw (5,-1.3) node(x2) [observed] {$X_{t+1}$};
  \draw[-arcsq] (dots) -- (h0);
  \draw[-arcsq] (h0) -- (h1);
  \draw[-arcsq] (h1) -- (h2);
  \draw[-arcsq] (h2) -- (dotss);
  \draw[-arcsq] (h0) -- (y0);
  \draw[-arcsq] (h1) -- (y1);
  \draw[-arcsq] (h2) -- (y2);
  \draw[-arcsq] (x0) -- (y0);
  \draw[-arcsq] (x1) -- (y1);
  \draw[-arcsq] (x2) -- (y2);
\end{tikzpicture}
\caption{A hidden Markov model for $(\B{Y}, \B{X})$. All observations (across different $t \in \{1, \dots, m\}$) are conditionally independent given $\B{H}$, and $(Y_t, X_t)$ only depends on $\B{H}$ through the present state $H_t$. Moreover, the variables in $\B{H}$ resemble a first order Markov chain, that is, ($H_1, \dots, H_{t-1}) \indep H_{t+1} \given H_t$ for all $t \in \{2, \dots, m-1\}$.}
\label{fig:HMM}
\end{figure}
We conveniently assume the independence of $X$ and $H$, which allows for likelihood inference 
without explicitly modelling the distribution of $X$. Our robustness analysis in Section~\ref{sec:sensitivity} suggests, however, 
that violations of this assumption do not negatively affect the performance of our causal discovery method.

For $i,j \in \{1, \dots, \ell\}$, let $\Gamma_{ij} = P(H_t =j \given H_{t-1} = i)$ denote the transition probabilities of $\B{H}$. By considering different parametrizations $\gamma \mapsto \B{\Gamma}(\gamma)$, $\gamma \in \mathcal{G}$, where $\mathcal{G}$ is a subset of a Euclidean space, we can encompass 
both of the above models simultaneously. The model IID then simply corresponds to a map $\B{\Gamma}$ satisfying that, for every $\gamma \in \mathcal{G}$, $\B{\Gamma}(\gamma)$ has constant columns. For details on the parametrizations of the models IID and HMM, see Appendix~\ref{app:para}.

\subsubsection{Notation}
The characteristics of the model for the joint distribution of $(\B{Y}, \B{X})$ are determined by the parametrizations $(\B{\Theta}, \mathcal{T})$ and $(\B{\Gamma}, \mathcal{G})$ of the regression matrix $\Theta$ and the transition matrix $\Gamma$, respectively. For every $\gamma \in \mathcal{G}$, let $\lambda (\gamma) = \lambda(\B{\Gamma}(\gamma)) \in \R^{1 \times \ell}$ be the stationary distribution of $\B{\Gamma}(\gamma)$. The stationary distribution $\lambda(\gamma)$ exists (and is unique) if the matrix $\B{\Gamma}(\gamma)$ is irreducible and aperiodic \citep[e.g.,][Propositions~1.31--1.33]{ching2006markov}. In the remainder of this work, we therefore require the image $\B{\Gamma}(\mathcal{G})$ to be a subset of the space of irreducible and aperiodic matrices of dimension $\ell \times \ell$.
We use $\mathcal{SR}_{(\B{\Theta}, \B{\Gamma})}(\theta, \gamma \given \B{X})$ to denote the joint distribution $P$ over $(\B{Y}, \B{X})$ with marginals $(Y_t, X_t) \sim \mathcal{SR}_\B{\Theta}(\theta, \lambda(\gamma) \given X_t)$ and a dependence structure given by $\B{\Gamma}(\gamma)$. Unless explicit parametrizations are referred to, we will usually omit the dependence on $\B{\Theta}$ and $\B{\Gamma}$ and simply write $\mathcal{SR} (\theta, \gamma \given \B{X})$. For every $j \in \{1, \dots, \ell \}$, we use $\beta_j(\cdot)$ and $\sigma_j^2(\cdot)$ to denote the parametrizations of the $j$th regression coefficient and the $j$th error variance, respectively, as induced by $(\B{\Theta}, \mathcal{T})$. Finally, $\phi$ denotes the combined parameter vector $(\theta, \gamma)$ with corresponding parameter space $\mathcal{P} := \mathcal{T} \times \mathcal{G}$.

\subsection{Likelihood} \label{sec:likelihood}
Consider a fixed pair of parametrizations $(\B{\Theta}, \mathcal{T})$ and $(\B{\Gamma}, \mathcal{G})$. For $(\theta, \gamma) \in \mathcal{T} \times \mathcal{G}$, the joint density of $(\B{Y}, \B{X}, \B{H})$ induced by the distribution $\mathcal{SR}(\theta, \gamma \given \B{X})$ is given by
\begin{equation*}
p_{(\B{\Theta}, \B{\Gamma})}( \B{y}, \B{x}, \B{h} \given \theta, \gamma) = p(\B{x}) \lambda(\gamma)_{h_1} \prod_{s=2}^m \B{\Gamma}_{h_{s-1} h_s} (\gamma) \prod_{t=1}^m \mathcal{N}(y_t \given x_t \beta_{h_t}(\theta), \sigma_{h_t}^2(\theta)),
\end{equation*}
where $p(\B{x})$ is the (unspecified) density of $\B{X}$, and where, for $j \in \{1, \dots, \ell\}$, $\mathcal{N}(y_t \given x_t \beta_j, \sigma_j^2)$ is short hand notation for the density of a $\mathcal{N}(x_t \beta_j, \sigma_j^2)$ distribution evaluated at $y_t$. 
Given a sample $(\B{y}, \B{x})$ from $(\B{Y}, \B{X})$, the loglikelihood function for the model $\{ \mathcal{SR}(\theta, \gamma \given \B{X}) : (\theta, \gamma) \in \mathcal{T} \times \mathcal{G} \}$ is then given by
\begin{equation} \label{eq:loglik} 
\ell_{(\B{\Theta}, \B{\Gamma})} (\B{y}, \B{x} \given \theta, \gamma) = \log \sum_{h_1} \cdots \sum_{h_m} p_{(\B{\Theta}, \B{\Gamma})}( \B{y}, \B{x}, \B{h} \given \theta, \gamma), \qquad (\theta, \gamma) \in \mathcal{T} \times \mathcal{G}.
\end{equation}
It is well known that, in general, the loglikelihood function~\eqref{eq:loglik} is non-concave and may have several local maxima. For unconstrained parametrizations $(\B{\Theta}, \mathcal{T})$ and $(\B{\Gamma}, \mathcal{G})$, it is even unbounded. To see this, one may, for example, choose $(\theta, \gamma) \in \mathcal{T} \times \mathcal{G}$ such that all entries of $\B{\Gamma}(\gamma)$ are strictly positive and such that $x_{t_0} \beta_1(\theta) = y_{t_0}$ for a single fixed $t_0$. By letting $\sigma_1^2 (\theta)$ go to zero while keeping all other regression parameters fixed, $p_{(\B{\Theta}, \B{\Gamma})}( \B{y}, \B{x}, \B{h} \given \theta, \gamma)$ approaches infinity for all $\B{h}$ with $h_t = 1 \Leftrightarrow t = t_0$. 

We consider two kinds of parameter constraints: (i) a lower bound on all error variances, and (ii) equality of all error variances. These constraints can be implemented using the 
parametrizations $(\B{\Theta}^c, \mathcal{T}^c)$ and $(\B{\Theta}^=, \mathcal{T}^=)$ given in Appendix~\ref{app:para}. 
In the following theorem, we show that either of these parametrizations ensures the existence of the maximum likelihood estimator. %
\begin{thm}[Existence of the MLE] \label{thm:existence}
	Let $(\B{y}, \B{x})$ be a sample of $(\B{Y}, \B{X}) = (Y_t, X_t)_{t\in \{1, \dots, m\}}$ and assume that the set $\{(y_t,x_t) \given t \in \{1, \dots, m\}\}$ is not contained in a union of $\ell$ hyperplanes of dimension $p$. 
	Let $\mathcal{G}$ be a compact subset of a Euclidean space and let $\B{\Gamma}: \mathcal{G} \to [0,1]^{\ell \times \ell}$ be a continuous parametrization of the transition matrix $\Gamma$. Then, with $(\B{\Theta}, \mathcal{T})$ being either of the parametrizations $(\B{\Theta}^c, \mathcal{T}^c)$ or $(\B{\Theta}^=, \mathcal{T}^=)$ (see Appendix~\ref{app:para}), the loglikelihood function $\ell_{(\B{\Theta}, \B{\Gamma})}$ attains its supremum on $\mathcal{T} \times \mathcal{G}$.
\end{thm}

The assumption involving hyperplanes excludes the possibility of a perfect fit. The conditions on $(\B{\Gamma}, \mathcal{G})$ ensure that the space of possible transition matrices is a compact set. The continuity of all parametrizations 
together with the parameter constraints inherent in $(\B{\Theta}^c, \mathcal{T}^c)$ and $(\B{\Theta}^=, \mathcal{T}^=)$ make for a continuous and bounded likelihood function. 
We use two different methods for likelihood optimization: a numerical optimization routine\footnote{We are grateful to Roland Langrock who 
who shared parts of his code with us.
} and an EM-type algorithm. These methods make use of the \verb|R| packages \verb|nlm| and \verb|mixreg|, respectively, and will be referred to as ``NLM'' and ``EM''; see Appendix~\ref{app:lik_opt} for details.

\subsection{Fisher Confidence Regions} \label{sec:fisher}
Using the asymptotic normality of maximum likelihood estimators, 
we can now construct (approximate) confidence regions for $\theta$. 
Let therefore $\hat \phi = (\hat \theta, \hat \gamma)$ be a global maximizer of the likelihood function and let $ \mathcal{J} (\hat \phi)$
be the observed Fisher information \citep[e.g.,][Chapter~2]{lehmann2006theory} at $\hat \phi$. 
For $\alpha \in (0,1)$, we define the region
\begin{equation} \label{eq:Calpha}
C^{\alpha}(\hat \theta) := \left \lbrace \hat \theta +  \mathcal{J}^{-1/2} (\hat \theta) v \, : \,  \norm{v}_2^2 \leq q_{\chi^2(\text{dim}(\theta))}(\alpha)\right \rbrace,
\end{equation}
where $\text{dim}(\theta)$ is the length of the parameter vector $\theta$, $q_{\chi^2(f)}(\alpha)$ is the $\alpha$-quantile of a $\chi^2(f)$-distribution and $\mathcal{J}^{-1/2}(\hat \theta)$ is the submatrix of $\mathcal{J}(\hat \phi)^{-1/2}$ corresponding to $\hat \theta$. 
For these confidence regions to achieve the correct asymptotic coverage, we need to adjust for the label switching problem described in the following subsection.

\subsection{Label Permutations} \label{sec:labels}
The distribution $\mathcal{SR}(\phi \given \B{X})$ is invariant under certain permutations of the coordinates of the parameter vector $\phi$. For example, when $\ell=2$, the hidden variable has two states. If we exchange all parameters corresponding to the first state with those corresponding to the second state, the induced mixture distribution is unchanged. 
In general, the model $\{ \mathcal{SR}(\phi \given \B{X}) : \phi \in \mathcal{P} \}$ is therefore not identifiable. 
More formally, let
$\Pi$ denote the set of all permutations of elements in $\{1, \dots, \ell \}$. For every permutation $\pi \in \Pi$ with associated permutation matrix $M_\pi$, define the induced mappings $\pi_\mathcal{T} := \B{\Theta}^{-1} \circ (\Theta \mapsto \Theta M_\pi^T) \circ \B{\Theta}$, $\pi_\mathcal{G} := \B{\Gamma}^{-1} \circ (\Gamma \mapsto M_\pi \Gamma M_\pi^T) \circ \B{\Gamma}$ and $\pi_\mathcal{P} := (\pi_\mathcal{T}, \pi_\mathcal{G})$ on $\mathcal{T}$, $\mathcal{G}$ and $\mathcal{P}$, respectively. Then, for every $\phi \in \mathcal{P}$ and every $\pi \in \Pi$, the distributions $\mathcal{SR}(\phi \given \B{X})$ and $\mathcal{SR}(\pi_\mathcal{P} (\phi) \given \B{X})$ coincide (and thus give rise to the same likelihood). 
The likelihood function therefore attains its optimum in a set of different parameter
vectors, all of which correspond to permutations of one another.
Coverage properties of the confidence region \eqref{eq:Calpha}
depend on which particular permutation of the MLE
is output by the optimization routine
(even though each of them parametrizes the exact same distribution). 
To overcome this ambiguity, we introduce the permutation-adjusted confidence regions 
\begin{equation} \label{eq:C_adj}
C_{\text{adjusted}}^{\alpha} (\hat \theta) := \bigcup_{\pi \in \Pi} C^\alpha( \pi_\mathcal{T}(\hat \theta)).
\end{equation}
In the following section, we make precise under which conditions these confidence regions achieve the correct asymptotic coverage.

\subsection{Asymptotic Coverage of Adjusted Confidence Regions} \label{sec:covgarant_cr}
Assume that the distribution of $X_t$ is stationary across  $e = \{1, \dots, m\}$ and has a density $f$ with respect to the Lebesgue measure on $\R^{p}$. 
Consider a fixed pair $(\B{\Theta}, \mathcal{T})$ and $(\B{\Gamma}, \mathcal{G})$ of parametrizations. Let $\phi^0 = (\theta^0, \gamma^0) \in \mathcal{P} := \mathcal{T} \times \mathcal{G}$ be the true parameters and let $\Theta^0 = \B{\Theta}(\theta^0)$ and $\Gamma^0 = \B{\Gamma}(\gamma^0)$ be the associated regression matrix and transition matrix, respectively.

Suppose now that the data within environment $e$ accumulates. For every $m \in \N$, write $(\B{Y}_m, \B{X}_m) = (Y_t, X_t)_{t \in \{1, \dots, m\}}$,
let $\P^m_0 := \mathcal{SR}(\theta^0, \gamma^0 \given \B{X}_m)$ and use 
 $\P_0$ to denote the (infinite-dimensional) limiting distribution of $\P_0^m$. Similarly, $\E_0$ denotes the expectation with respect to $\P_0$. 
We require the following assumptions.
\begin{itemize}
\setlength\itemsep{0em}
\item [(A1)] The maximum likelihood estimator exists. 
\item [(A2)] The true parameter $\phi^0$ is contained in the interior of $\mathcal{P}$.
\item [(A3)] The transition matrix $\Gamma^0$ is irreducible and aperiodic \citep[e.g.,][Section~1]{ching2006markov}. 
\item [(A4)] For every $i \in \{1, \dots, p+1 \}$ and $j,k \in \{1, \dots, \ell\}$, the maps $\theta \mapsto \B{\Theta}_{ij}(\theta)$ and $\gamma \mapsto \B{\Gamma}_{jk}(\gamma)$ have two continuous derivatives. 
\item [(A5)] 
For every $m \in \N$, assume that the joint distribution of $(\B{Y}_m, \B{X}_m)$ 
has a density with respect to the Lebesgue measure that we denote by $f_m$. 
Then, the Fisher information matrix $\mathcal{I}_0$ defined as
\begin{equation*}
\mathcal{I}_0 := \E_{0} [\eta \eta^T], \quad \text{where} \quad \eta = \lim_{m \to \infty} \left. \frac{\partial}{\partial \phi} f_m(Y_m, X_m \given \B{Y}_{m-1}, \B{X}_{m-1}, \phi) \right \vert_{\phi = \phi^0},
\end{equation*}
is strictly positive definite.
\item [(A6)] All coordinates of $X_1$ have finite fourth moment.
\item [(A7)] $\E [\card{\log f(X_1)}] < \infty$.
\end{itemize}
Assumptions~(A1)~and~(A4) are satisfied for the explicit parametrizations of the models IID and HMM given in Appendix~\ref{app:para}, see Theorem~\ref{thm:existence}. 
The irreducibility of $\Gamma^0$ assumed in (A3) guarantees all latent states to be visited infinitely often, such that information on all parameters keeps accumulating. Assumption~(A5) is needed to ensure that, in the limit, the loglikelihood function has, on average, negative curvature and hence a local maximum at $\phi^0$. Finally, (A6) and (A7) are mild regularity conditions on the (otherwise unspecified) distribution of $X_t$.

Essentially, the asymptotic validity of the adjusted confidence regions~\eqref{eq:C_adj} rests on two results: (1) consistency of the MLE and (2) asymptotic normality of the MLE. For every $\phi \in \mathcal{P}$, let $[\phi] := \{\pi_\mathcal{P}(\phi) : \pi \in \Pi \}\subseteq \mathcal{P}$ denote the equivalence class of $\phi$, i.e., the set of parameters in $\mathcal{P}$ that are equal to $\phi$ up to a permutation $\pi_\mathcal{P}$ as defined in Section~\ref{sec:labels}. Consistency in the quotient topology (``$[\hat \phi_m] \rightarrow [\phi^0]$'') then simply means that any open subset of $\mathcal{P}$ that contains the equivalence class of $\phi^0$, must, for large enough $m$, also contain the equivalence class $\hat \phi_m$. 
With this notation, we can now state an asymptotic coverage result for confidence regions~\eqref{eq:C_adj}. 
The main work is contained in Theorems~\ref{thm:consistency}~and~\ref{thm:normality}.
Their proofs 
make use of results given by \citet{leroux1992maximum} and \citet{bickel1998asymptotic}, 
which discuss 
consistency and asymptotic normality, respectively, of the MLE in hidden Markov models with finite state space. 
\begin{thm}[Consistency of the MLE] \label{thm:consistency}
Assume that (A1), (A3), (A4) and (A7) hold true. Then, $\P_0$-almost surely, $[\hat \phi_m] \rightarrow [\phi^0]$ as $m \to \infty$. 
\end{thm}
Theorem~\ref{thm:consistency} says that $(\hat \phi_m)_{m \in \N}$ alternates between one or more subsequences, each of which is convergent to a permutation of $\phi^0$. The following theorem proves a central limit theorem for these subsequences. 
\begin{thm}[Asymptotic normality of the MLE] \label{thm:normality}
Assume that the maximum likelihood estimator is consistent. Then, under (A1)--(A6), it holds that $\mathcal{J}(\hat \phi_m)^{1/2}(\hat \phi_m - \phi^0)  \stackrel{d}{\longrightarrow} \mathcal{N}(0, I)$ under $\P_0$.
\end{thm}
Together, Theorems~\ref{thm:consistency}~and~\ref{thm:normality} imply the following asymptotic coverage guarantee. 
\begin{cor}[Asymptotic coverage of adjusted confidence regions] \label{cor:cr_validity}
Un\-der Assumptions (A1)--(A7), the adjusted confidence regions \eqref{eq:C_adj} achieve the correct asymptotic coverage. That is, for any $\alpha \in (0,1)$,
\begin{equation} \label{eq:ci_asymp}
\liminf_{m \to \infty} \P^m_0 (\theta^0 \in C^\alpha_{\emph{adjusted}} ( \hat \theta_m)) \geq 1 - \alpha.
\end{equation}
\end{cor}
As another corollary, the asymptotic type I error control of the tests defined by~\eqref{eq:T} follows by applying Corollary~\ref{cor:cr_validity} to each environment separately. 

\section{ICPH: Algorithm and False Discovery Control} \label{sec:covgarant}
We can now summarize the above sections into our overall method. In Section~\ref{sec:pseudo} we provide a pseudo code for this procedure, and Section~\ref{sec:covgarant_icph} presents our main theoretical result---an asymptotic version of Proposition~\ref{prop:icp_level}, which states that our procedure is consistent. 
\subsection{Algorithm} \label{sec:pseudo}
Given data $(\B{Y}, \B{X})$ and a collection $\mathcal{E}$ of environments, we run through all $S \subseteq \{1, \dots, d\}$, test the hypothesis $H_{0,S}$ with the test defined by~\eqref{eq:T} using the adjusted confidence regions~\eqref{eq:C_adj}, and output the intersection of all accepted sets. Below, this procedure is formalized in a pseudo code. 

\begin{algorithm}[H] \caption{ICPH (``Invariant Causal Prediction in the presence of Hidden variables'')} \label{alg:icph}
\textbf{Input}: response $\B{Y} \in \R^n$, covariates $\B{X} \in \R^{n \times d}$,  environment indicator $\B{E} \in \{1, \dots, \card{\mathcal{E}} \}^n$ (i.e., $\B{E}_t = k \Leftrightarrow t \in e_k$)\;
\textbf{Options}: $\verb|model| \in \{ \text{``IID'', ``HMM''}\}$, $\verb|method| \in \{\text{``EM'', ``NLM''}\}$, $\verb|variance.constraint| \in \{\text{``lower bound'', ``equality''}\}$, $\verb|number.of.states| \in \N_{\geq 2}$, $\verb|intercept| \in \{ \text{TRUE, FALSE} \}$, $\verb|test.parameters| \subseteq \{\text{``intercept'', ``beta'', ``sigma''} \}$, $\verb|alpha| \in (0,1)$\;
 \For{$S \subset \{1, \dots, d\}$}{
 	\For{$e \in \mathcal{E}$}{
		 Fit an $\mathcal{SR}$ model to $(\B{Y}_e, \B{X}^S_e)$, see Section~\ref{sec:likelihood}\;
		 Construct the permutation-adjusted confidence region \eqref{eq:C_adj}\; 	
 	}
 	Compute a $p$-value $p_S$ for $H_{0,S}$ using the test defined by \eqref{eq:T}\;
 }
 \textbf{Output}: the empirical estimator $\hat S = \bigcap_{S: p_S > \alpha} S$\;
\end{algorithm}

Most of the options in Algorithm~\ref{alg:icph} are self-explanatory. The option \verb|test.parameters| allows the user to specify the ``degree of $h$-invariance'' that is required of the sets $S \subset \{1, \dots, d\}$. If, for example, $\verb|test.parameters| = \{\text{``beta'', ``sigma''}\}$, a set $S$ will be regarded $h$-invariant if the mixture components of $P_{Y_t \vert X_t^S}$ are ``invariant in $\beta$ and $\sigma^2$'', i.e., time-homogeneous up to changes in the intercept between different environments. 
Code is available online (see Section~\ref{sec:organization}). 
To make Algorithm~\ref{alg:icph} scalable to a large number of predictors, it can be combined with a 
variable screening step, e.g., using Lasso \citep{Tibshirani94};
see Section~\ref{sec:covgarant_icph} for more details.

\subsection{Asymptotic False Discovery Control of ICPH} \label{sec:covgarant_icph}
The cornerstone for the false discovery control of ICPH is given in Corollary~\ref{cor:cr_validity}. It proves that if Assumptions~(A1)--(A7) are satisfied for the true set $S^*$, then the test $\varphi_{S^*}$ achieves the correct asymptotic level, which in turn guarantees an asymptotic version of Proposition~\ref{prop:icp_level}. We will now summarize this line of reasoning into out main theoretical result. 

Assume that we are given data $((\B{Y}_n, \B{X}_n))_{n \in \N} = \left((Y_{n,t}, X_{n,t})_{t \in \{1, \dots, n\}}\right)_{n \in \N}$ from a triangular array, where, for every $n$, $(\B{Y}_n, \B{X}_n) \in \R^{n \times (1 + d)}$. Consider a fixed number of $K$ environments and let $(\mathcal{E}_n)_{n \in \N}$ be a sequence of collections $\mathcal{E}_n = \{ e_{n,1}, \dots, e_{n,K} \}$, such that, for all $n$, $e_{n,1}, \dots, e_{n,K}$ are disjoint with $\cup_k e_{n,k} = \{1, \dots, n\}$ and such that, for all $k$, $\card{e_{n,k}} \rightarrow \infty$ as $n \to \infty$. For all $n$ and $k$, write $(\B{Y}_{n,k}, \B{X}_{n,k}) = (Y_t, X_t)_{t \in e_{n,k}}$. Consider a transition parametrization $(\B{\Gamma}, \mathcal{G})$ and a family of regression parametrizations $\{(\B{\Theta}^S, \mathcal{T}^S)\}_{S \subseteq \{1, \dots, d\}}$, i.e., for every $S \subseteq \{1, \dots, d\}$, $\B{\Theta}^S$ maps $\mathcal{T}^S$ into the space of matrices of dimension $(\card{S}+1) \times \ell$ with columns in $\R^{\card{S}} \times \R_{>0}$. For every $n$ and every $S \subseteq \{1, \dots, d\}$, let $H^n_{0,S}$ denote the hypothesis~\eqref{eq:H0} for the data $(\B{Y}_n, \B{X}^S_n)$ and let $\varphi_S^n$ be the corresponding test defined by \eqref{eq:T} with the confidence regions~\eqref{eq:C_adj}. Finally, define for every $n$ the estimator 
\begin{equation} \label{eq:hatSn}
\hat S_n := \bigcap_{S: \, \varphi^n_S \text{ accepts } H^n_{0,S}} S.
\end{equation}
We then have the following result. 
\begin{thm}[Asymptotic false discovery control] \label{thm:asymp_cov}
Assume that Assumption~\ref{ass:h_inv} is satisfied. That is, there exists a set $S^* \subseteq \{1, \dots, d\}$ which, for every $n$, is $h$-invariant with respect to $(\B{Y}_n, \B{X}_n)$. 
Assume furthermore that, for every $k$, (A1)--(A7) hold true for the data $(\B{Y}_{n,k}, \B{X}^{S^*}_{n,k})$ with parametrizations $(\B{\Theta}^{S^*}, \mathcal{T}^{S^*})$ and $(\B{\Gamma}, \mathcal{G})$. Then, the estimator $\hat S_n$ enjoys the following coverage property
\begin{equation} \label{eq:asymp_cov}
\liminf_{n \to \infty} \P^n_0 (\hat S_n \subseteq S^*) \geq 1 - \alpha,
\end{equation}
where $\P^n_0$ is the law of $(\B{Y}_n, \B{X}_n)$.
\end{thm}
If the number of predictor variables is large, our algorithm can be combined with an upfront variable screening step. Given a family $(\hat S^n_{\text{screening}})_{n \in \N}$ of screening estimators, we can for every $n\in \N$ construct an estimator $\bar{S}_n$ of $S^*$ analogously to \eqref{eq:hatSn}, but where the intersection is taken only over those $S$ additionally satisfying that $S \subseteq \hat{S}^n_{\text{screening}}$. Given that $\liminf_{n \to \infty} \P^n_0(S^* \subseteq \hat S^n_{\text{screening}}) \geq 1 - \alpha$, it then follows from 
\begin{align*}
\P^n_0(\bar{S}_n \not \subseteq S^*) 	&= \P^n_0([\bar{S}_n \not \subseteq S^*] \cap [S^* \subseteq \hat S^n_\text{screening}]) + \P^n_0([\bar{S}_n \not \subseteq S^*] \cap [S^* \not \subseteq \hat S^n_\text{screening}]) \\
																	&\leq \P^n_0(\varphi^n_{S^*} \text{ rejects } H^n_{0,S^*}) + \P^n_0 (S^* \not \subseteq \hat S^n_\text{screening}),
\end{align*}
that the estimator $(\bar{S}_n)_{n \in \N}$ satisfies the asymptotic false discovery control \eqref{eq:asymp_cov} at level $1-2\alpha$. In high-dimensional models, assumptions that allow for the screening property have been studied \citep[see, e.g.,][]{Buhlmann2011}.

\section{Experiments} \label{sec:experiments}
In this section, we apply our method to simulated data (Section~\ref{sec:simulation}) and to a real world data set on photosynthetic activity and sun-induced fluorescence (Section~\ref{sec:SIF}). 
We only report results using the NLM optimizer. In all experiments, the results for EM were almost identical to those for NLM, except that the computation time for EM was larger (by approximately a factor of 6). For an experiment that uses the EM-method, see Appendix~\ref{app:em}. 

\subsection{Simulated Data} \label{sec:simulation}
We start by testing the sample properties of the adjusted confidence regions, disregarding the problem of causal discovery, see Section~\ref{sec:coverage_cr}. 
In Section~\ref{sec:SCM}, 
we present the multivariate data generating process that we will use in the subsequent
analyses.
In Section~\ref{sec:level}, we see that, even for sample sizes that are too small for the confidence regions to achieve the correct coverage, our overall method (ICPH) is able to keep the type I error control. Section~\ref{sec:power} contains a power analysis. In Section~\ref{sec:sensitivity}, we test the robustness of ICPH against a range of different model violations, and include a comparison with two alternative causal discovery methods. 
The performance of ICPH for non-binary latent variables, 
for large numbers of predictor variables,
or under violations of the $h$-invariance assumption,
can be found in Appendix~\ref{app:extrasim}. 

\subsubsection{Empirical coverage properties of adjusted confidence regions} \label{sec:coverage_cr}
The finite sample coverage properties of the confidence regions \eqref{eq:C_adj} depend on the true distribution over $(\B{Y}, \B{X})$ (i.e., on the parameters of the $\mathcal{SR}$ model as well on the marginal distribution of $\B{X}$) and on the sample size. We here illustrate this sensitivity in the i.i.d.\ setting. 
Consider a joint distribution $\P$ over $(Y,X,H) \in \R^{1+p} \times \{1, \dots, \ell\}$ which induces an $\mathcal{SR}$ model over $(Y,X)$. For every $j\in \{1, \dots, \ell \}$ let $p_j(y,x) = \P(H = j \given Y = y, X = x)$ denote the posterior probability of state $j$ based on the data $(y,x)$. We then use the geometric mean of expected posterior probabilities 
\begin{equation} \label{eq:pi}
\text{GMEP} :=\left(  \prod_{j = 1}^\ell \E[p_j(Y,X) \given H = j] \right)^{1/ \ell} \in [0,1]
\end{equation}
as a measure of difficulty of fitting the $\mathcal{SR}$ model induced by $\P$.\footnote{If each of the distributions $\P_{(Y, X) \vert H=j}$, $j \in \{1, \dots, \ell\}$ has a density w.r.t.\ the Lebesgue measure on $\R^{1+p}$,
each factor in \eqref{eq:pi} is given as an integral over $\R^{1+p}$. In practice, we approximate these integrals by numerical integration.} We expect smaller values of GMEP to correspond to 
more difficult
estimation problems, which negatively affect the convergence rate of \eqref{eq:ci_asymp} and result in low finite sample coverage.
If the between-states differences in the regression parameters of $X$ are small, 
for example,
we expect the unobserved states to be difficult to infer from the observed data
(i.e., for every $j$, the expected posterior probabilities $\E[p_i(Y,X) \given H = j]$ are close to uniform in $i$),
resulting in small GMEP.

\begin{figure}
\centering
\includegraphics[width = .29\linewidth]{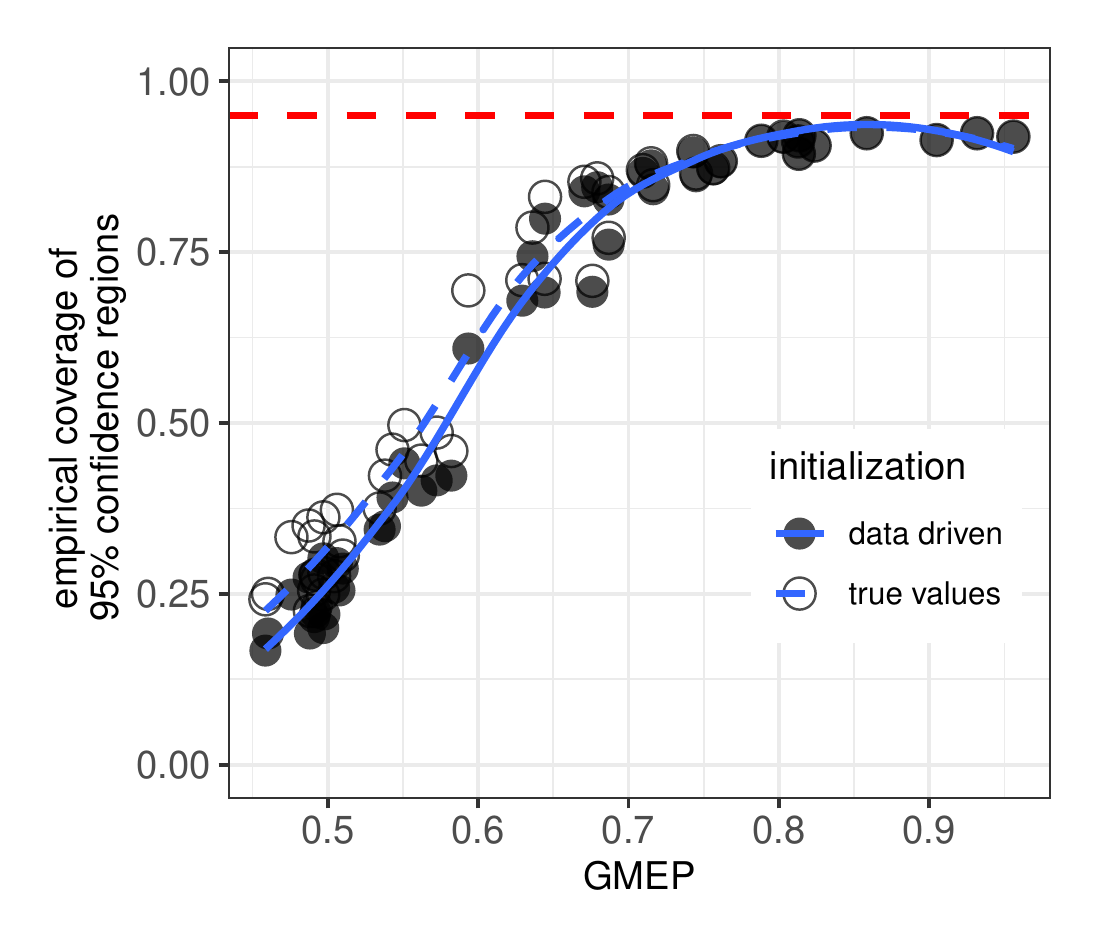}
$\quad$
\includegraphics[width = .31\linewidth]{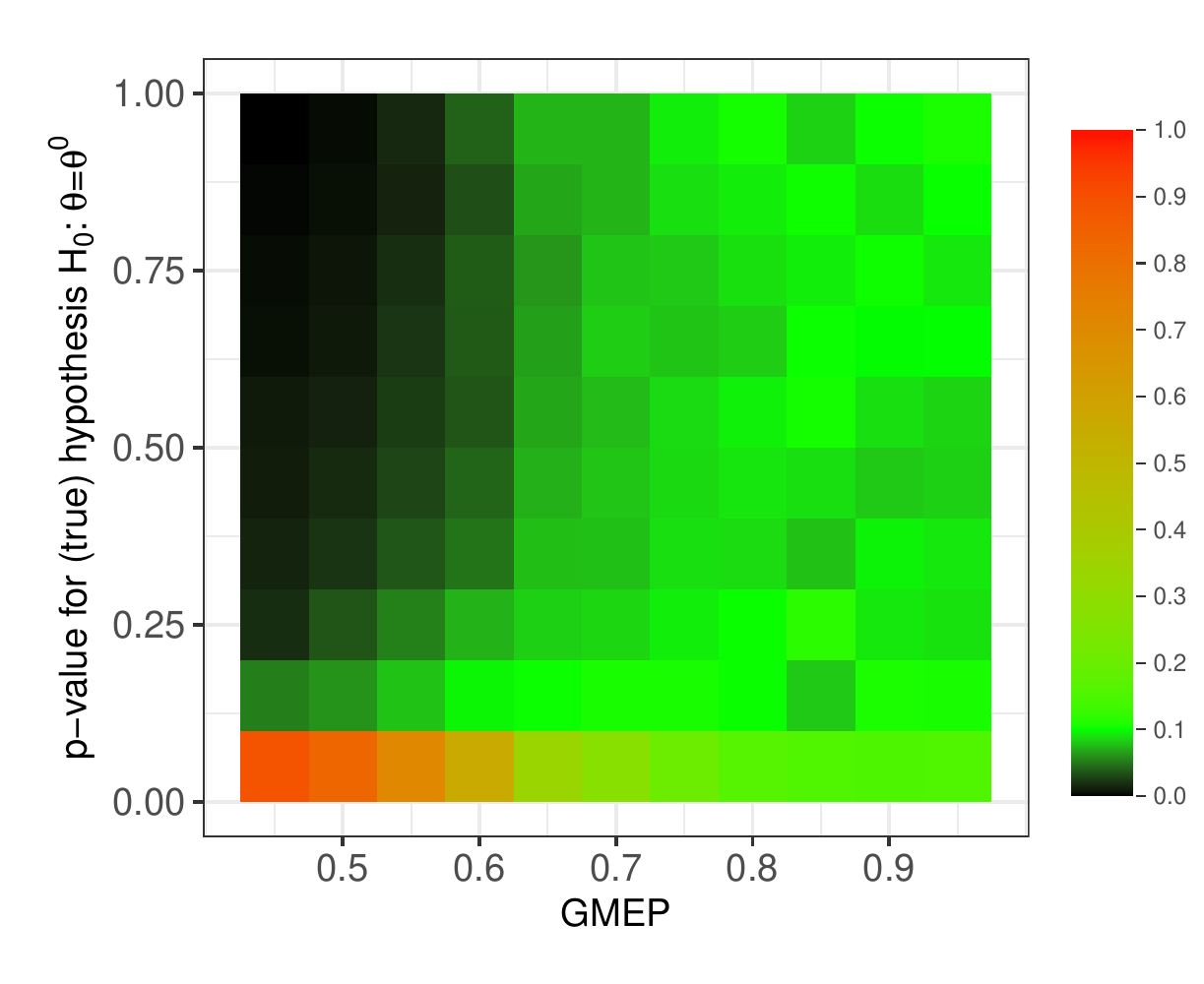}
$\quad$
\includegraphics[width = .29\linewidth]{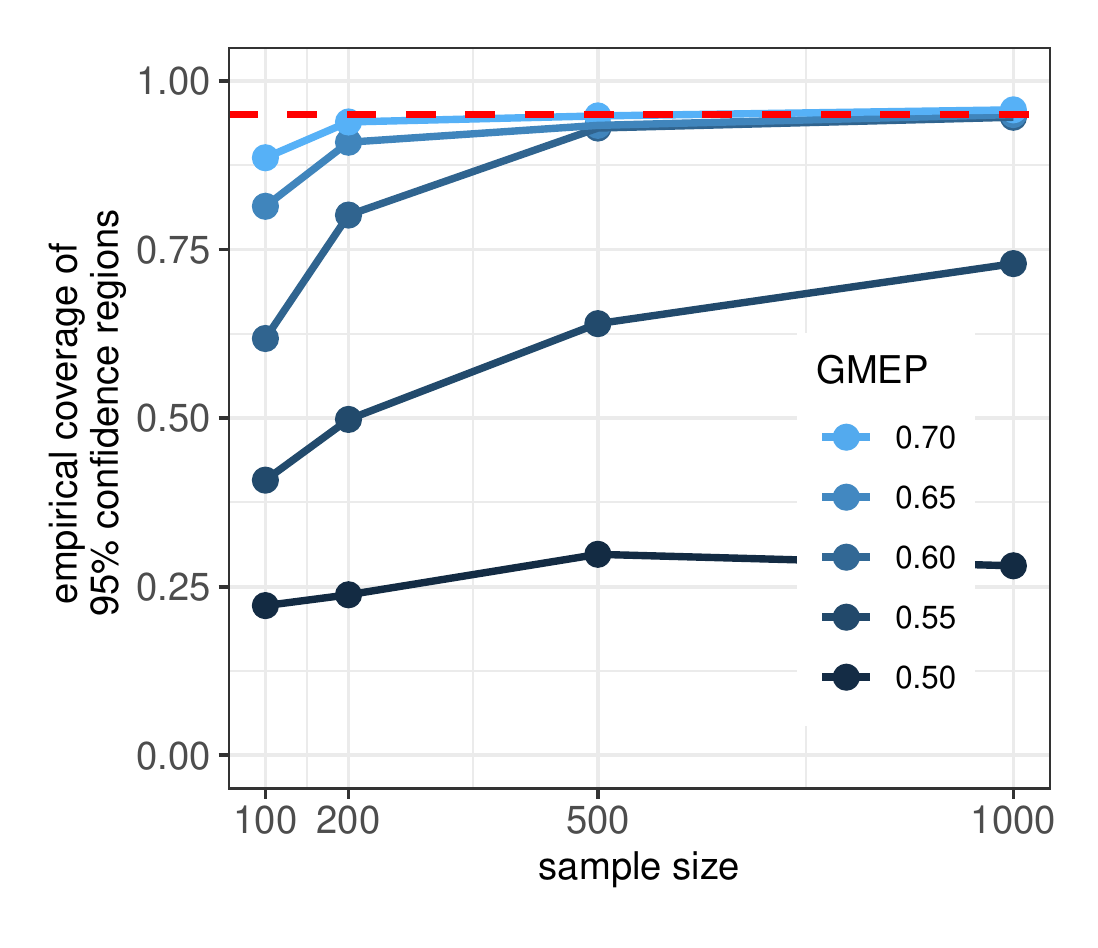}
\caption{
Empirical coverage properties of the adjusted confidence regions~\eqref{eq:C_adj} using data simulated from the model~\eqref{eq:SCM}. %
The left panel shows empirical coverage 
of $95\%$-confidence regions for different model parameters (see Equation \ref{eq:pi} for a definition of GMEP), and a fixed sample size of 100. 
We see that the coverage properties strongly depend on GMEP, 
and that the poor performance for low GMEP is 
not an optimization problem 
(the likelihood scores obtained from starting the algorithm in the true values 
exceed those obtained from data driven initialization in less than $0.2\%$ of simulations). 
In the middle panel, 
we use the same simulations, but only
consider data-driven initialization. 
Each column corresponds to a histogram of $p$-values~\eqref{eq:pval}. 
For increasing GMEP, the $p$-value distribution 
approximates the desired uniform distribution.
For 5 different parameter settings, we further increase the sample size (right). As suggested by Corollary~\ref{cor:cr_validity}, the empirical coverage gradually improves, although very low GMEP demand large amounts of data to obtain satisfactory coverage. 
}
\label{fig:coverage}
\end{figure}

We now perform the following simulation study.
For different model parameters and sample sizes, we generate i.i.d.\ data sets from the SCM
\begin{equation} \label{eq:SCM}
H := N^H_\lambda, \quad X := \mu_X  +\sigma_X N^X,  \quad Y := \mu_Y + \beta_1 X \cdot \mathbbm{1}_{ \{ H = 1 \}} + \beta_2  X \cdot \mathbbm{1}_{ \{ H = 2 \}} + \sigma_Y N^Y,
\end{equation}
where all noise variables are jointly independent with marginal distributions $N^H_\lambda \sim \text{Ber}(\lambda)$, $N^X, N^Y \sim \mathcal{N}(0,1)$. We construct adjusted confidence regions \eqref{eq:C_adj} for the vector of regression parameters $\theta^0 = (\mu_Y, \beta_1, \mu_Y, \beta_2, \sigma_Y^2)$ 
using the likelihood function \eqref{eq:loglik} with parametrizations $(\B{\Theta}^=, \mathcal{T}^=)$ and $(\B{\Gamma}^{\text{IID}}, \mathcal{G}^\text{IID})$ (see Appendix~\ref{app:para}). 
We first sample 50 different sets of parameters independently as $\mu_X, \mu_Y, \beta_1, \beta_2 \sim \text{Uniform}(-1,1)$, $\sigma_X \sim \text{Uniform}(0.1,1)$, $\sigma_Y \sim \text{Uniform}(0.1,0.5)$ and $\lambda \sim \text{Uniform}(0.3,0.7)$. 
For each setting, 
we compute empirical coverage degrees based on 1000 independent data sets, 
each consisting of 100 independent replications from \eqref{eq:SCM}, 
and compare them to the GMEP of the underlying models,
see Figure~\ref{fig:coverage} (left). 
For the same simulations, we also compare the $p$-values 
\begin{equation} \label{eq:pval}
p := \max \{ \alpha \in [0,1] : \theta^0 \not \in C_{\text{adjusted}}^\alpha (\hat \theta) \}
\end{equation}
for the (true) hypotheses $H_0: \theta = \theta^0$ to a uniform distribution (Figure~\ref{fig:coverage} middle). 
For 5 models of different degrees of difficulty ($\text{GMEP} \approx 0.50, 0.55, 0.60, 0.65, 0.70$), we then compute empirical coverage degrees for increasing sample size (Figure~\ref{fig:coverage} right). 

For difficult 
estimation problems (i.e., low GMEP), the finite sample variance of the MLE is inflated, resulting in low empirical coverage and too small $p$-values (Figure~\ref{fig:coverage} left and middle). Although there is no proof that NLM finds the global optimum, it is assuring that 
there is little difference when we start the algorithm at the (usually unknown) true values
(Figure~\ref{fig:coverage} left, hollow circles). Indeed, the thus obtained likelihood scores exceed those obtained from data driven initialization in less than $0.2\%$ of simulations. 
As seen in Figure~\ref{fig:coverage} (right), coverage properties improve with increasing sample size, although in models with low GMEP, we require large amounts of data in order to obtain satisfactory performance. We will see in Section~\ref{sec:level} that even cases where we cannot expect the confidence regions to obtain valid coverage, our overall causal discovery method maintains type I error control.

\subsubsection{Data generating process} \label{sec:SCM}
We now specify the data generating process used in the following sections. We consider an SCM over the system of variables $(Y, X^1, X^2, X^3, H)$ given by the structural assignments
\begin{align*}
X^1	&:= N^1,  	\quad 	H 	:= N^H, \quad		X^2 := \beta^{21} X^1 + N^2 \\
Y		&:= \sum_{j=1}^\ell (\mu^Y_j + \beta^{Y}_{1j} X^1 + \beta^{Y}_{2j} X^2 + \sigma_{Yj} N^Y) \mathbbm{1}_{\{H = j\}} \\
X^3 	&:= \beta^{3Y} Y + N^3, 
\end{align*}
where $N^H \sim \text{Multinomial}(1, \lambda)$, $N^Y \sim \mathcal{N}(0,1)$ and $N^j \sim \mathcal{N}(\mu^j, \sigma_j^2)$. 
In Sections~\ref{sec:level}--\ref{sec:sensitivity}, the latent variable $H$ is assumed to be binary, while Appendix~\ref{app:non_binary} treats the more general case where $\ell \geq 2$.
The different environments are constructed as follows.
We first draw random change points $1 < t_1 < t_2 < n$ and then generate data as described below. 
\begin{figure}
	\centering
	\includegraphics[width = \linewidth]{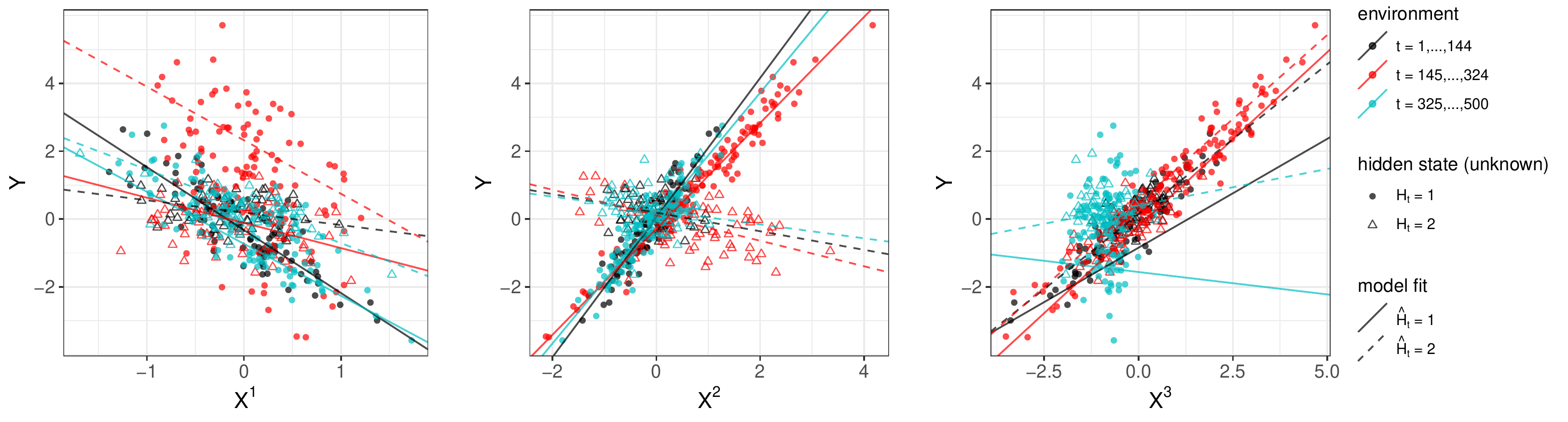}
	\caption{Data generated from the SCM described in Section~\ref{sec:SCM} 
		for each of the three environments (black, red, blue). 
		Here, 
		the only $h$-invariant set is $S^* = \{1, 2\}$ and 
		we would therefore like
		our method to correctly identify the 
		violations of the $h$-invariance of the sets $\{1\}, \{2\}$ and $\{3\}$. 
		These violations are indicated by the different model fits (colored lines), which for none of the three variables are stable across all environments. 
		For numerical results on such data sets, see Sections~\ref{sec:level} and~\ref{sec:power}.
		The issue of label permutations can be seen from the 
		occasional mismatch between the true latent states ($\bullet: H_t = 1$, $\vartriangle: H_t = 2$) and the estimated labels (\sampleline{}: $\hat H_t = 1$, \sampleline{dashed}:~$\hat H_t = 2$).
	}
	\label{fig:simdata}
\end{figure}

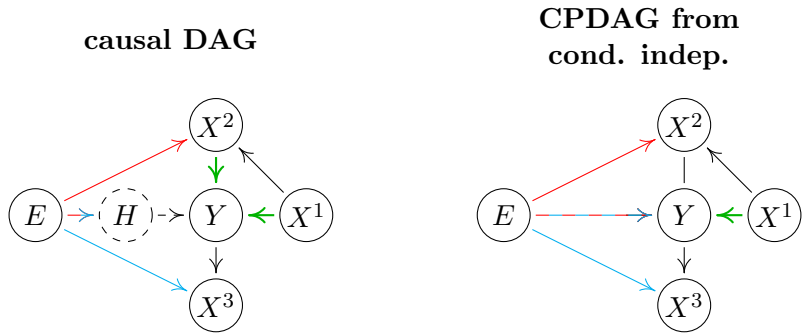
\begin{figure}
\begin{center}
\begin{minipage}{0.4\textwidth}
	\begin{center}
	{\bf causal DAG}
	\end{center}
\end{minipage}
\begin{minipage}{0.4\textwidth}
	\begin{center}
		{\bf CPDAG from \\ cond. indep.}
	\end{center}
\end{minipage}
\end{center}

\vspace{-5mm}

\begin{center}
\begin{minipage}{0.4\textwidth}
	\begin{center}
		\begin{tikzpicture}[xscale=1.2, yscale=1.2, shorten >=1pt, shorten <=1pt]
		\draw (0,0) node(x2) [observedsmall] {$X^2$};
		\draw (-2,-1) node(E) [observedsmall] {$E$};
		\draw (0,-1) node(y) [observedsmall] {$Y$};
		\draw (-1,-1) node(h) [unobservedsmall] {$H$};
		\draw (1,-1) node(x1) [observedsmall] {$X^1$};
		\draw (0,-2) node(x3) [observedsmall] {$X^3$};

		\draw[-arcsq, green!70!black, thick] (x1) -- (y);
		\draw[-arcsq, dashed] (h) -- (y);
		\draw[-arcsq, green!70!black, thick] (x2) -- (y);
		\draw[-arcsq] (x1) -- (x2);
		\draw[-arcsq] (y) -- (x3);
		\draw[-arcsq, dashed,dash pattern=on 5pt off 5pt, red] (E) -- (h);
		\draw[-arcsq, dashed, dash pattern=on 5pt off 5pt, cyan, dash phase=5pt] (E) -- (h);
		\draw[-arcsq, red] (E) -- (x2);
		\draw[-arcsq, cyan] (E) -- (x3);
		\end{tikzpicture}
	\end{center}
\end{minipage}
\begin{minipage}{0.4\textwidth}
	\begin{center}
		\begin{tikzpicture}[xscale=1.2, yscale=1.2, shorten >=1pt, shorten <=1pt]
		\draw (0,0) node(x2) [observedsmall] {$X^2$};
		\draw (-2,-1) node(E) [observedsmall] {$E$};
		\draw (0,-1) node(y) [observedsmall] {$Y$};
		\draw (1,-1) node(x1) [observedsmall] {$X^1$};
		\draw (0,-2) node(x3) [observedsmall] {$X^3$};

		\draw[-arcsq, green!70!black, thick] (x1) -- (y);
		\draw[-arcsq] (h) -- (y);
		\draw (x2) -- (y);
		\draw[-arcsq] (x1) -- (x2);
		\draw[-arcsq] (y) -- (x3);
		\draw[-arcsq, dashed,dash pattern=on 5pt off 5pt, red] (E) -- (y);
		\draw[-arcsq, dashed, dash pattern=on 5pt off 5pt, cyan, dash phase=5pt] (E) -- (y);
		\draw[-arcsq, red] (E) -- (x2);
		\draw[-arcsq, cyan] (E) -- (x3);
		\end{tikzpicture}
	\end{center}
\end{minipage}
\end{center}
\caption{
Left: the causal graph induced by the SCM in Section~\ref{sec:SCM}. The node $E$ represents the different environments 
($E$ points into variables that have been intervened on, the color corresponds to the environments shown in Figure~\ref{fig:simdata}). 
Right: the CPDAG representing the Markov equivalence class of the graph where $H$ is marginalized out. 
Since the edge $X^2 - Y$ is not oriented, the full set of causal parents $S^* = \{1,2\}$ cannot be identified 
only from conditional independence statements. 
Our method exploits the simple form of the influence of $H$ on $Y$.
Note that in the case of an additional edge $E \to X^1$, 
none of the edges among the variables $(Y, X^1, X^2, X^3)$ would be oriented in the CPDAG. 
}
\label{fig:graphs}
\end{figure}

\begin{itemize}
	\setlength\itemsep{0em}
	\item $e_1 = \{1, \dots, t_1\}$: Here, we sample from the observational distribution.
	\item $e_2 = \{t_1+1, \dots, t_2\}$: Here, we set $X^2 := \beta^{21} X^1 +\tilde{N}^2$, where $\tilde{N}^2$ is a Gaussian random variable with mean sampled uniformly between 1 and 1.5 and variance
	sampled uniformly between 1 and 1.5. Also, the mixing proportions $\lambda$ are resampled.
	\item $e_3 = \{t_2+1, \dots, n \}$: We again sample data from the above SCM, but this time we intervene on $X^3$. The structural assignment is replaced by $X^3 := \tilde{N}^3$, where $\tilde{N}_3$ is a Gaussian random variable with mean sampled uniformly between $-1$ and $-0.5$ and the same variance as the noise $N^3$ from the observational setting. The mixing proportions $\lambda$ are again resampled.
\end{itemize}
A sample data set can be seen in Figure~\ref{fig:simdata}, where points have been colored according to the above environments (black, red and blue for $e_1$, $e_2$ and $e_3$, respectively). 
The only $h$-invariant set is the set $S^* = \{1,2\}$ of observable parents of $Y$. In the population case, our method therefore correctly infers $\tilde{S} = \{1,2\}$, see Equation~\eqref{eq:Stilde2}. 
The causal graph induced by the above data generating system can be seen in Figure~\ref{fig:graphs} (left).
Here, the environment is drawn as a random variable.\footnote{To view the data set as i.i.d.\ realizations from such a model one formally adds a random permutation of the data set, which breaks the dependence of the realizations of the environment variable (this has no effect on the causal discovery algorithm, of course).
\citet{Constantinou2017} discuss a non-stochastic treatment of such nodes.}
We also display the CPDAG representing the Markov equivalence class of the induced graph over the observed variables
(right),
showing that the full set of causal parents $S^* = \{1,2\}$ cannot be identified only from conditional independence statements.

\subsubsection{Level analysis} \label{sec:level}
\begin{figure}
\centering
\includegraphics[width = .7\linewidth]{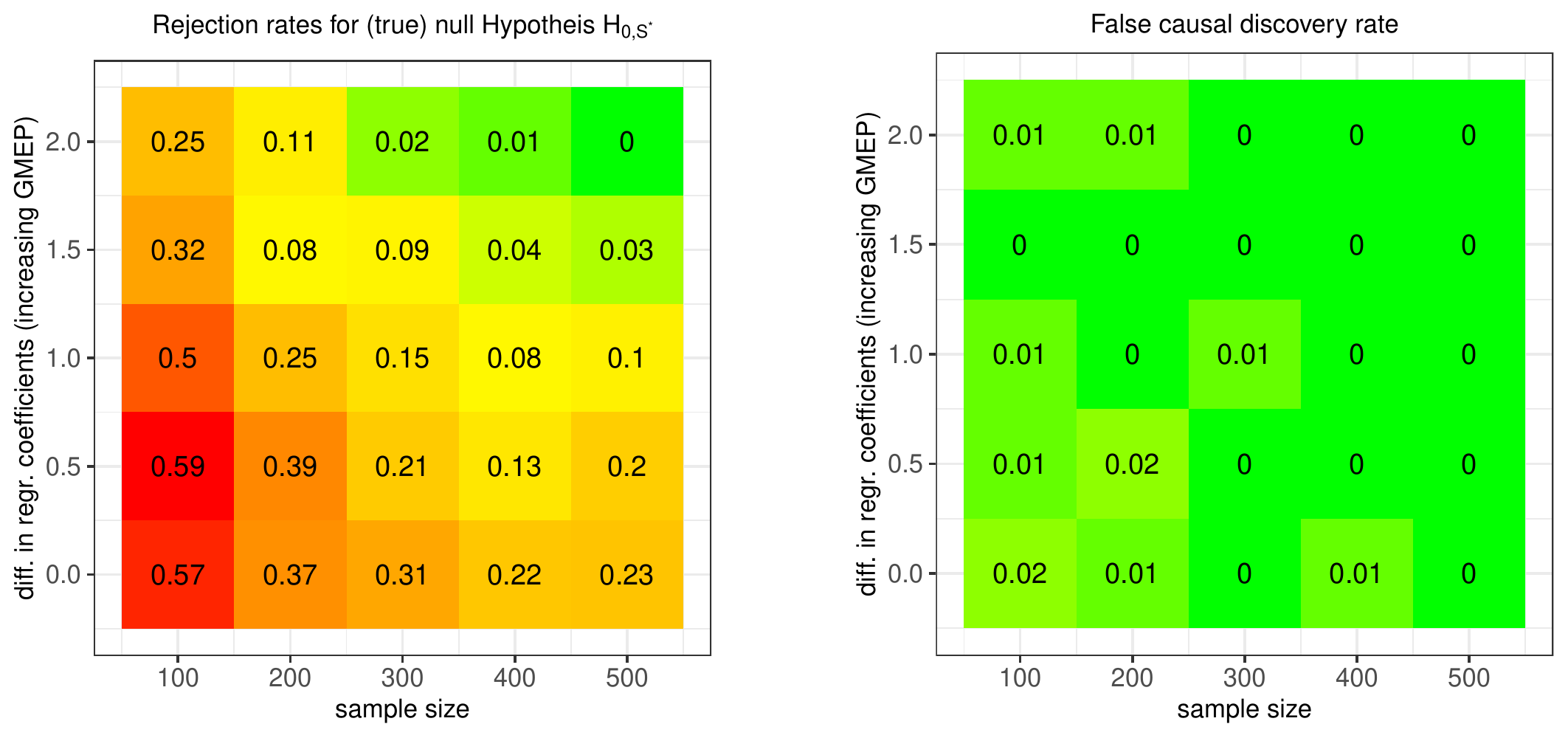}
\caption{Estimates $\hat \P(\varphi_{S^*} \text{ rejects } H_{0, S^*})$ (left) and $\hat \P(\hat S \not \subseteq S^*)$ (right) of the type I error rates of the test $\varphi_{S^*}$ and the overall method ICPH, respectively, based on the experiment described in Section~\ref{sec:level} and 100 repetitions. 
The desired level is $\alpha = 0.05$. 
We have used NLM with parametrizations $\B{\Theta}^=$ and $\B{\Gamma}^{\text{IID}}$ (see Appendix~\ref{app:para}). 
The average GMEP values are $0.51, 0.56, 0.64, 0.66, 0.78$ (ordered in accordance to the vertical axis). 
For small sample sizes, and in particular for low GMEP, the type I error control of the test $\varphi_{S^*}$ is violated. 
Even in these cases, however,
the false causal discovery control of ICPH is satisfied.
}
\label{fig:robust_beta}
\end{figure}

Given that the theoretical coverage guarantees are only asymptotic, we cannot expect the tests~\eqref{eq:T} to satisfy type I error control for small sample sizes---especially if GMEP is low, see also Section~\ref{sec:coverage_cr}. 
The following empirical experiments suggest, however, that 
even if the test level of the true hypothesis $H_{0,S^*}$ is violated, ICPH may still keep the overall false discovery control.
We use data sets $ (Y_t, X_t^1, X_t^2, X_t^3)_{t \in \{1, \dots, n\}} $ generated 
as described in 
Section~\ref{sec:SCM}, and 
analyse the performance of ICPH for different sample sizes and different GMEP. Since the latter is difficult to control directly, we vary the between-states difference in regression coefficients for $X^1$ and $X^2$ in the structural assignment for $Y$, and report the average GMEP for each setting.
For every $n \in \{100, 200, 300, 400, 500\}$ and every $\Delta \beta \in \{0, 0.5, 1, 1.5, 2\}$, we simulate 100 independent data sets by drawing model parameters $\mu \sim^{iid} %
\text{Uniform}(-0.2, 0.2)$, 
$\sigma^2 \sim^{iid} %
\text{Uniform}(0.1, 0.3)$ (with the restriction that $\sigma_{Y1}^2 = \sigma_{Y2}^2$), $\beta \sim^{iid} %
\text{Uniform}([-1.5, -0.5] \cup [0.5,1.5])$ and $\lambda \sim \text{Uniform}(0.3, 0.7)$. For $j \in \{1,2\}$ we then assign $\beta_{j,2}^Y := \beta_{j,1}^Y + \text{sign}(\beta_{j,1}^Y) {\Delta \beta}$.
The results are summarized in Figure~\ref{fig:robust_beta}. We see that even in settings for which the true hypothesis $H_{0,S^*}$ is rejected for about every other simulation, ICPH stays conservative. 

\subsubsection{Power analysis} \label{sec:power}
\begin{figure}
\begin{center}
\begin{minipage}{0.47\textwidth}
\centering
\includegraphics[width = \linewidth]{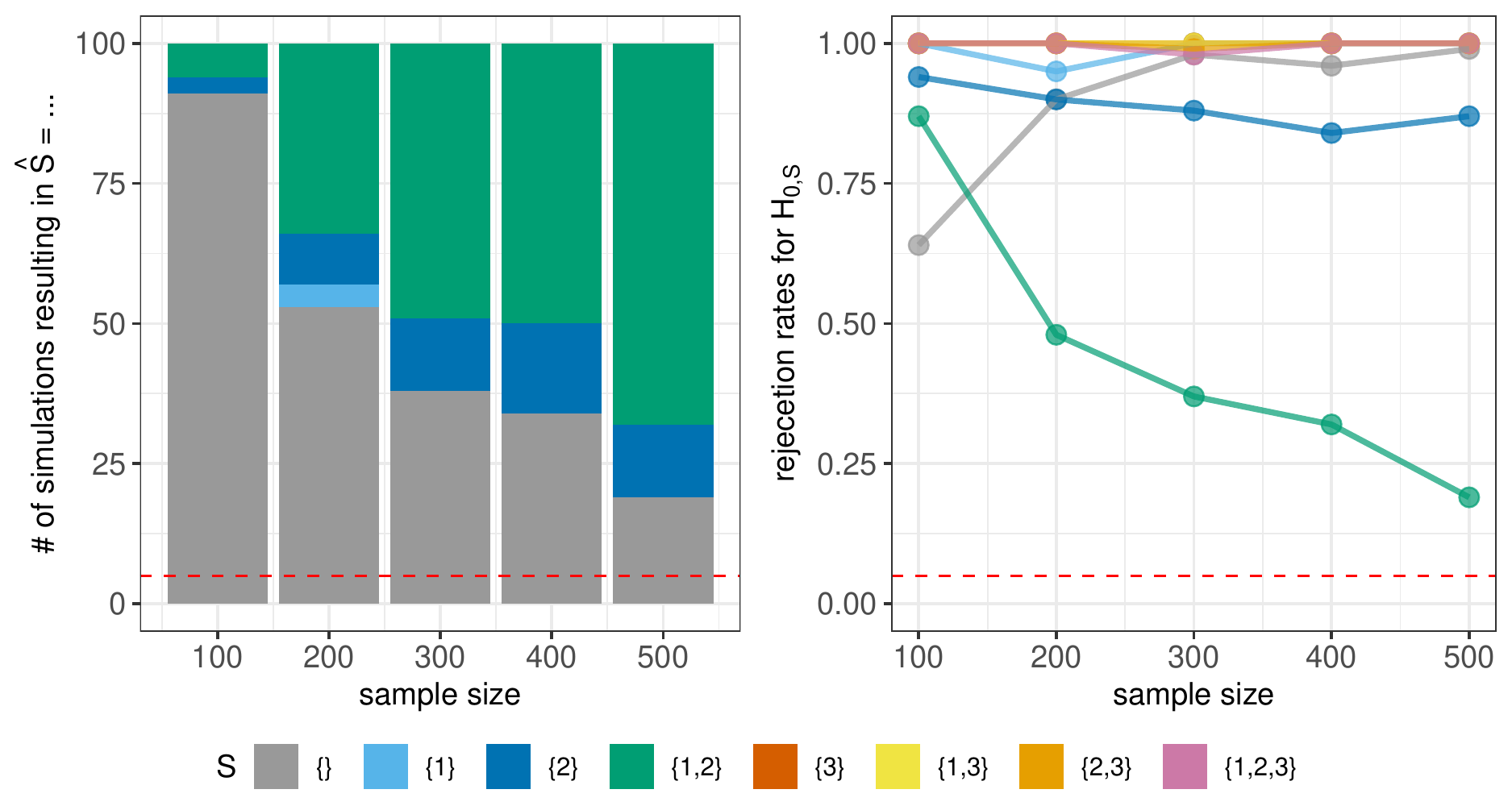}
\end{minipage}
\quad
\begin{minipage}{0.47\textwidth}
\centering
\includegraphics[width = \linewidth]{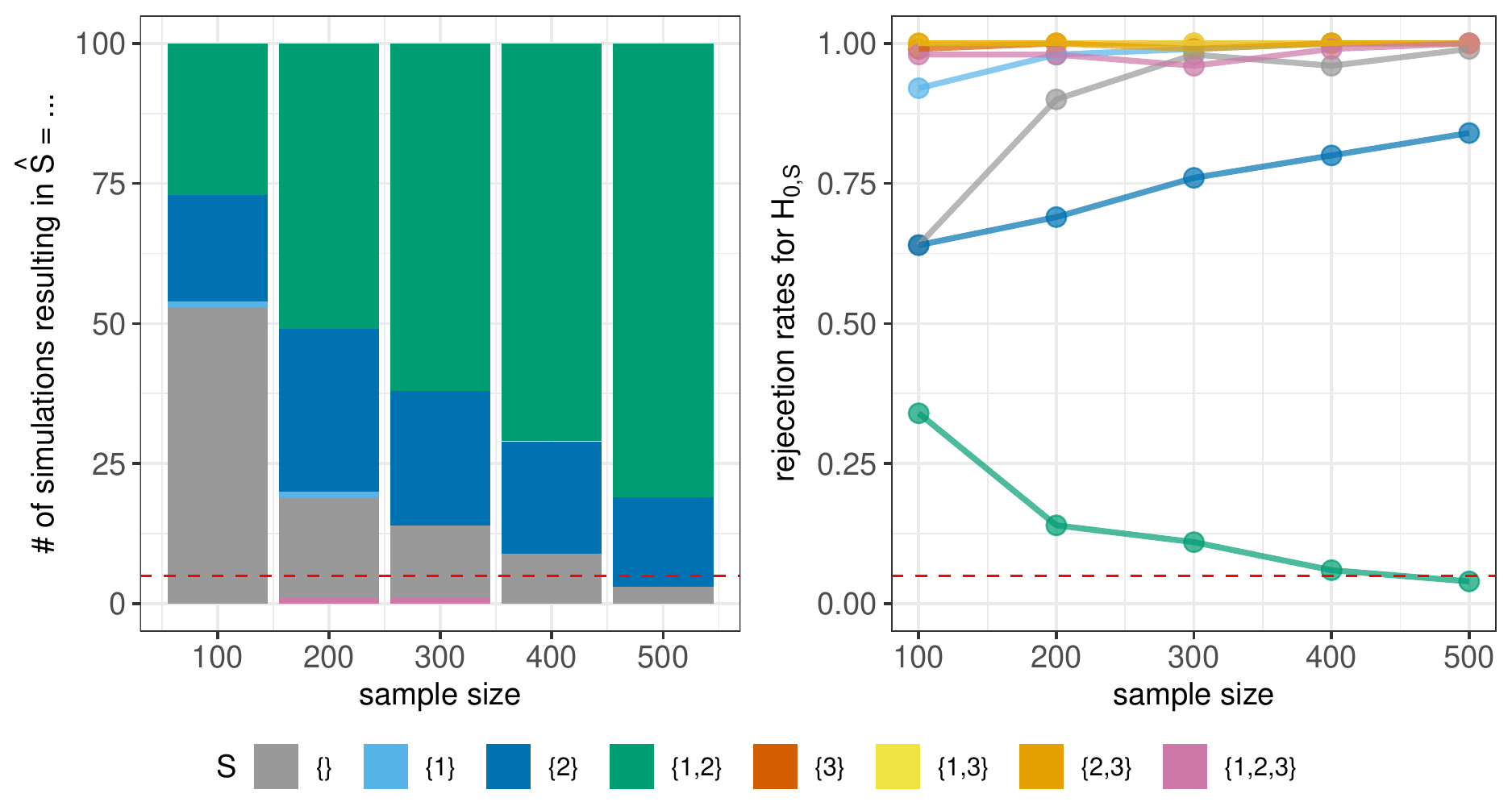}
\end{minipage}
\end{center}

\caption{
Output of ICPH (bar plots) and rejection rates for individual hypotheses (curve plots) for the experiment in Section~\ref{sec:power} with parameter constraint $\sigma_{Y1}^2, \sigma_{Y2}^2 \geq 10^{-4}$ (left) and $\sigma_{Y1}^2 = \sigma_{Y2}^2$ (right). 
The larger the proportion of blue and green colors in the bar plots, the more power our method has. Simulations are ordered such that, within each bar, the bottom colors 
(yellow, light orange, dark orange, purple) correspond to false positives, i.e., cases where $\hat S \not \subseteq S^*$. Even though the level of the test for $H_{0,S^*}$ is violated in the finite sample case, 
ICPH controls the empirical type I error rate at $\alpha = 0.05$
(indicated by a dashed horizontal line).
Enforcing equality on error variances is beneficial, especially for small data sets. 
For both settings, the identification of $S^*$ improves with increasing sample size.
}
\label{fig:shat}
\end{figure}

Since the only $ h $-invariant set is the set $S^* = \{1,2\}$ of causal parents of $ Y $, the population version of our method correctly infers $\tilde{S} = \{1,2\}$, see Equation~\eqref{eq:Stilde2}. For finite samples, identifiability of $ S^* $ is determined by the power of the tests for the hypotheses  $ H_{0,S} $. For a fixed value of $\Delta \beta = 1.5$ (average GMEP of 0.66) and increasing sample size, we generate i.i.d.\ data sets as described in Section~\ref{sec:level} and analyze the performance of ICPH for two different variance constraints $\sigma_{Y1}^2 = \sigma_{Y2}^2$ and $\sigma_{Y1}^2, \sigma_{Y2}^2 \geq 10^{-4}$. The results in Figure~\ref{fig:shat} suggest that the former constraint results in higher performance, and it will therefore be our default setting for the rest of this section. As the sample size increases, ICPH tends to identify the set $S^*$ (larges  shares of green in bar plots).

For the same data that generated Figure~\ref{fig:shat}, we compute rejection rates for non-causality (i.e., empirical proportions of not being contained in $\hat S$) for each of the predictors 
$X^1$, $X^2$ and $X^3$. 
Here, we also add a comparison to other methods. 
We are not aware of any other method that is 
suitable for inferring $S^*$,
but we nevertheless add two approaches as baseline.
\begin{itemize}
	\setlength\itemsep{0em}
	\item ``$k$-means ICP'': Pool data points from all environments and infer estimates $\hat H$ of the hidden states using $2$-means clustering. Run the ordinary ICP algorithm \citep{peters2016causal} on each of the data sets $\{(Y_t, X_t) \, : \, \hat H_t = j \}$, $j \in \{1,2\}$, testing all hypotheses at level $\alpha/2$, and obtain $\hat S_1$ and $\hat S_2$. Output the final estimate $\hat S = S_1 \cup S_2$.
	\item ``JCI-PC'': We use a modified version of the PC algorithm \citep{Spirtes2000}, which exploits our background knowledge of $E$ being a source node: in between skeleton search and edge orientation, we orient all edges connecting $E$ to another node. The resulting algorithm may be viewed as a variant of the
	of JCI algorithm \citep{mooij2019joint}. We apply it to the full system of observed variables $(E,Y, X^1, X^2, X^3)$, and output the set of variables (among $\{X^1, X^2, X^3\}$) which have a directed edge to $Y$ in the resulting PDAG.\footnote{
Note that $H$ can be marginalized out, so it is not necessary to use FCI. Furthermore, since we do not assume the intervention targets to be known, search algorithms for interventional data such as the GIES algorithm \citep{hauser2012characterization} are not applicable.}
\end{itemize}

In the JCI-PC algorithm, we use conditional independence tests based on partial correlations. 
Since we apply it to a system of mixed variables (i.e., continuous as well as discrete), 
the assumptions underlying some of the involved tests will necessarily be violated. 
We are not aware of any family of tests which is more suitable. 
However, even in the population case, we cannot expect constraint-based methods such as 
JCI-PC to infer the set full $S^*$, see Figure~\ref{fig:graphs}.
ICPH solves a specific problem
and is the only method which exploits the simple influence
of $H$ on $Y$.
The results in Figure~\ref{fig:sensitivity} (black curves) confirm our previous findings: 
causal discovery improves with increasing sample size, and our method stays conservative. ICPH
outperforms both other methods in terms of level and power.

\subsubsection{Robustness analysis} \label{sec:sensitivity}
Our results are based on various assumptions, and we now investigate the robustness of ICPH
against different kinds of model violations. We use data generated from the following modified versions of the SCM in Section~\ref{sec:SCM}. Unless mentioned otherwise, parameters are sampled as described in Section~\ref{sec:level}.

\begin{figure}
\centering
\includegraphics[width =  \linewidth]{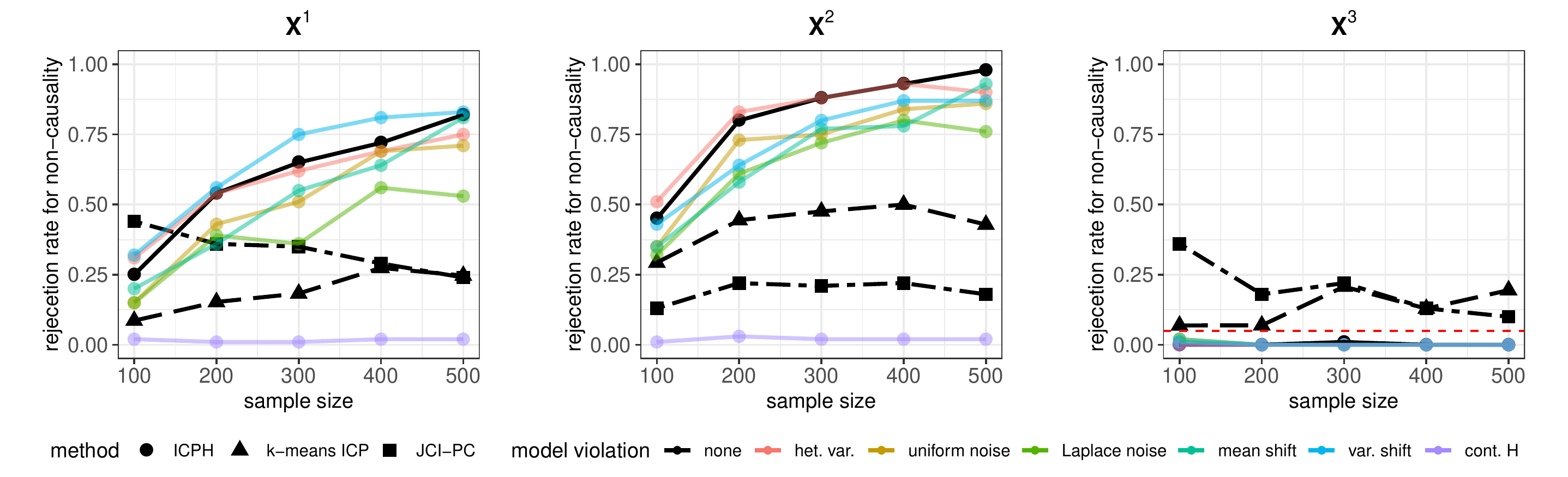}
\caption{Rejection rates for non-causality. This figure contains two comparisons (one among all black curves, and another among all curves with round points). For data generated from the SCM in Section~\ref{sec:SCM} (black), we compare the performance of ICPH ({\Large$\bullet$}) against the two alternative methods $k$-means ICP ($\blacktriangle$) and JCI-PC ({\footnotesize$\blacksquare$}) described in Section~\ref{sec:power}. For increasing sample size, ICPH 
outperforms both
methods in terms level and power. As a robustness analysis, we further we apply ICPH to simulated data sets from the modified SCMs described in Section~\ref{sec:sensitivity} (colored). Each of the modified SCMs yields a misspecification of the model for $P_{Y \vert X^{S^*}}$ that is assumed by our method. Most of these model misspecifications do not qualitatively affect the results: for increasing sample size, both causal parents $X^1$ and $X^2$ tend to be identified. For a continuous hidden variable, none of the variables is identified as causal (which is not incorrect, but uninformative). In all scenarios, ICPH maintains empirical type I error control.
}
\label{fig:sensitivity}
\end{figure}

\begin{itemize}
	\setlength\itemsep{0em}
	\item Heterogeneous variances: The error variances $\sigma_{Y1}^2$ and $\sigma_{Y2}^2$ are sampled independently.
	\item Non-Gaussian noise: We generate errors $N^Y$ from (i) a uniform distribution and (ii) a Laplace distribution. 
	\item A direct effect $H \to X^1$: We allow for an influence of $H$ on $X^1$ through binary shifts in (i) the mean value 
	and (ii) the error variance. 
	Parameters are sampled independently as $\mu_1^1, \mu_2^1 \sim \text{Uniform}(-1,1)$ and $\sigma_{11}^2, \sigma_{12}^2 \sim \text{Uniform}(0.1,1)$.
	\item A continuous hidden variable: We substitute the structural assignment for $Y$ by $Y:= (\mu^Y + \beta_1^Y X^1 + \beta_2^Y X^2) H + \sigma_Y N^Y$, where $H \sim \mathcal{N}(0,1)$. The distribution of $H$ does not change across environments.
\end{itemize}
We now repeat the power analysis from Section~\ref{sec:power} for data sets generated in the above way (Figure~\ref{fig:sensitivity}, colored curves). 
Most model violations do not qualitatively affect the results. Only the assumption on the state space of $H$ is crucial for the power (not the level) of our method; for a continuous hidden variable, we mostly 
output the empty set.

\subsection{Sun-Induced Fluorescence and Land Cover Classification} \label{sec:SIF}
We now consider a real world data set for which we can compare our method's output against a plausible causal model constructed from background knowledge. 
The data set is related to the study of global carbon cycles, which are determined by the movement of carbon between land, atmosphere and ocean. 
Carbon is emitted, e.g., during fossil fuel combustion, land use change or cellular respiration, and assimilated back into the Earth's surface by processes of carbon fixation. A major component hereof is photosynthesis, where inorganic carbon is converted into organic compounds by terrestrial ecosystems.
Direct measurements of carbon fluxes can be obtained from fluxtowers (\url{http://fluxnet.fluxdata.org}), but
are only available at single locations.
Constructing reliable global models for predicting photosynthesis using satellite data is an active line of research. 
While most of the commonly used models \citep[e.g.,][]{jung2009towards, running2015daily} use sunlight as the predominant driver, recent work \citep[e.g.,][]{guanter2012retrieval, zhang2016consistency} explores the predictive potential of sun-induced fluorescence (SIF), a (remotely sensible) electromagnetic radiation that is emitted by plants during the photosynthetic process. 
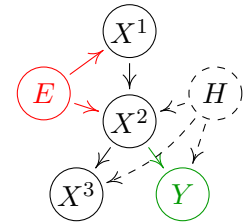
\begin{figure}
\begin{minipage}{0.7\textwidth}
\begin{center}
\begin{tabular}{c|l}
variable & description \\ \hhline{=|=}
$Y$ & sun-induced fluorescence (SIF) \\
$X^1$ & incident shortwave radiation (SW) \\ 
$X^2$ & absorbed photosynthetically active radiation ($\text{APAR}_{\text{chl}}$) \\
$X^3$ & gross primary productivity (GPP) \\
$H$ & vegetation type \\
\end{tabular}
\end{center}
\end{minipage}
\begin{minipage}{0.29\textwidth}
\begin{center}
	\hspace{10mm}
\begin{tikzpicture}[xscale=1, yscale=1, shorten >=1pt, shorten <=1pt]
  \draw (-1.141,0.371) node(e) [observedsmall, draw = red] {\textcolor{red}{$E$}};
  \draw (0,1.2) node(x1) [observedsmall] {$X^1$};
  \draw (1.141,0.371) node(h) [observedsmall, dashed] {$H$};
  \draw (0,0) node(x2) [observedsmall] {$X^2$};
  \draw (-0.706,-0.971) node(x3) [observedsmall] {$X^3$};
  \draw (0.706,-0.971) node(y) [observedsmall, draw=green!60!black] {\textcolor{green!60!black}{$Y$}};
  
  \draw[-arcsq, red] (e) -- (x1);
  \draw[-arcsq, red] (e) -- (x2);
  \draw[-arcsq] (x1) -- (x2);
  \draw[-arcsq] (x2) -- (x3);
  \draw[-arcsq, draw=green!60!black] (x2) -- (y);
  \draw[-arcsq, dashed] (h) [bend left = 0] to (x2);
  \draw[-arcsq, dashed] (h) [bend left = 0] to (y);
  \draw[-arcsq, dashed] (h) [bend left = 15] to (x3);
\end{tikzpicture}
\end{center}
\end{minipage}
\caption{Variable descriptions (left) and causal graph constructed from background knowledge (right). 
In our analysis, we use the temporal ordering of data to construct the environment variable $E$. 
Due to seasonal cycles of aggradation and degradation of chlorophyll, $\text{APAR}_{\text{chl}}$ is not a constant fraction of SW (which itself is time-heterogeneous). The environment therefore
``acts'' on the variables $X^1$ and $X^2$. Furthermore, different vegetation types differ not only in their chlorophyll composition (and thus in $\text{APAR}_{\text{chl}}$), but also in their respective efficiencies of converting $\text{APAR}_{\text{chl}}$ into GPP and SIF---hence the arrows from $H$ to $X^2$, $X^3$ and $Y$. }
\label{fig:data}
\end{figure}

Here, we take SIF as the target variable. 
As predictors, we include the incident shortwave radiation (SW), the photosynthetically active radiation absorbed 
by the plants' chlorophyll cells ($\text{APAR}_{\text{chl}}$), and the gross primary productivity (GPP), the latter of which is a measure of photosynthesis. Since GPP cannot be directly measured, %
we use spatially upscaled measurements from a network of fluxtowers \citep{jung2009towards}. 
Background knowledge %
suggests that out of these three variables, 
only $\text{APAR}_{\text{chl}}$ is a direct causal parent of the target SIF. 
\citet{zhang2016consistency} suggest evidence for a linear relationship between SIF and $\text{APAR}_{\text{chl}}$, 
and show that this relationship strongly depends on the type of vegetation. 
Estimates of the vegetation type can be obtained from the IGBP global land cover data base \citep{loveland2000development}.
We use the IGBP classification to select data coming from two different vegetation types only. 
In the resulting data set, we thus expect the causal influence of SIF on $\text{APAR}_{\text{chl}}$ to be confounded by a binary variable. 
When applying our method to these data, we remove information on vegetation type, so that this binary variable becomes latent. 
The data and the ground truth we consider is shown in Figure~\ref{fig:data}.

In Section~\ref{sec:realcd}, we use our causal discovery method to identify the causal predictor of SIF. 
In Section~\ref{sec:reconstruction},
we explore the possibility to
reconstruct the vegetation type from the observed data $(Y, X^1, X^2, X^3)$
when assuming that we have 
inferred the correct causal model.
We believe that such estimates
may be used to complement conventional vegetation type classifications.

\subsubsection{Causal discovery} \label{sec:realcd}
We denote the observed variables by $(Y, X^1, X^2, X^3)$ as described in Figure~\ref{fig:data} (left). The data are observed along a spatio-temporal grid with a temporal resolution of 1 month (Jan 2010 -- Dec 2010), and a spatial resolution of $0.5^\circ \times 0.5^\circ$ covering the North American continent. The setup is directly taken from \citet{zhang2016consistency}, and we refer to their work for a precise description of the data preprocessing for the variables $(Y, X^2, X^3)$. The data for $X^1$ is publicly available at \url{https://search.earthdata.nasa.gov}. 
We select
pixels classified as either \textit{Cropland (CRO)} or \textit{Evergreen Needleleaf Forest (ENF)}.
These vegetation types are expected to differ in their respective relationships $X^2 \rightarrow Y$ \citep{zhang2016consistency}. 
As environments we use the periods Feb -- Jul and Aug -- Jan.\footnote{
We also conducted the experiments with 
alternative constructions of the environments. 
Since switching regression models are hard to fit 
if the distribution of the predictors strongly differs
between states, 
some choices of environments 
make our method output the empty set---a result that is not incorrect, but uninformative.
}

\begin{figure}
	\centering 
	\includegraphics[width=\linewidth]{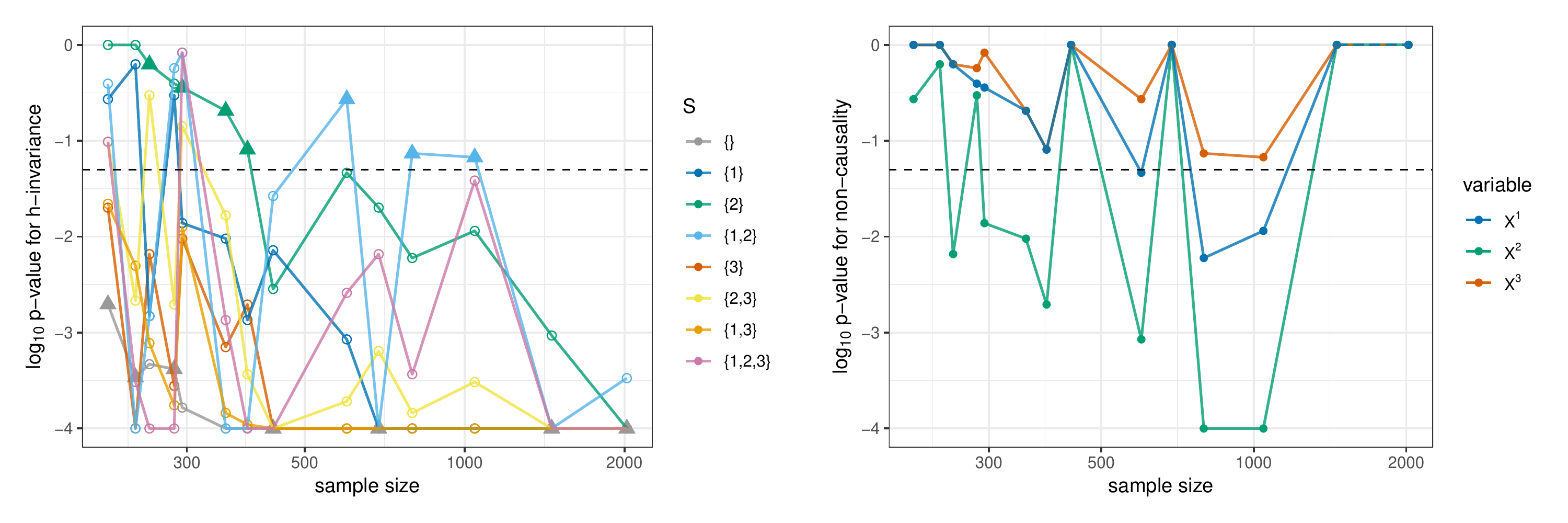}
	\caption{\textit{P}-values for $h$-invariance of different sets $S \subseteq \{1,2,3\}$ (left) and \textit{p}-values for non-causality (see Section~\ref{sec:pvals}) of the individual variables $X^1, X^2$ and $X^3$ (right). 
		For every experiment, the estimated set $\hat S$ in the left plot is indicated by a triangle. For several sample sizes, our method correctly infers $\hat S = \{2\}$ (left), and the causal parent $X^2$ consistently obtains the lowest $p$-value for non-causality (right). 
		Experiments for which all $p$-values for non-causality are equal to 1 correspond to instances in which all sets have been rejected.
		For large amounts of data, this is always the case (the two largest sample sizes are not shown here).
		At sample sizes 436, 797 and 1045, our method infers the set $\hat S = \{1,2\}$. This finding may be due to imperfect measurements of the variable $X^2$, that do not contain all information from $X^1$ that is relevant for $Y$. 
	}
	\label{fig:pval_vs_n}
\end{figure}

The goal of the statistical analysis is to identify the set $S^* = \{2\}$ of causal parents of $Y$ among the vector $(X^1, X^2, X^3)$. Since the variables $X^1$ and $X^2$ are closely related, we regard distinguishing between their respective causal relevance for $Y$ as a difficult problem.
We analyze the data for different sample sizes. To do so, we gradually lower the spatial resolution %
in the following way. For every $c \in \{1, \dots, 16\}$, we 
construct a new data set by 
increasing the pixel size of the original data set by a factor of $c^2$, and then averaging observations within each pixel. Grid cells that do not purely contain observations from either of the two vegetation types are discarded. We then apply our causal discovery method to each of the generated data sets, allowing for a binary hidden variable. The results are illustrated in Figure~\ref{fig:pval_vs_n}.\footnote{We omit all intercept terms, impose an equality constraint on the error variances, and assume an i.i.d.\ structure on the hidden variables. For estimation, we use the NLM optimizer. In our implementation of the test \eqref{eq:T}, the lowest attainable \textit{p}-value is $10^{-4}$.}
Indeed, for several sample sizes 
($n \leq 390$),
the true hypothesis $H_{0, S^*}$ is accepted, and our method mostly correctly infers $\hat S = \{2\}$ (left plot). In all experiments, the variable $X^2$ is attributed the highest significance as a causal parent of $Y$ (right plot). Also, we consistently do not
reject the only non-ancestrial variable $X^3$, and the causal ordering implied by the right hand plot is in line with the assumed causal structure from Figure~\ref{fig:data}. 
As the sample size grows, the power of our tests of the hypotheses $H_{0,S}$ increases, and even small differences in regression coefficients are detected. For sample sizes above 1459 (the two largest sample sizes are not shown here), all hypotheses $H_{0,S}$ are rejected, and our method returns the uninformative output $\hat S = \emptyset$. At sample sizes 436, 797 and 1045, our method infers the set $\hat S = \{1,2\}$, that is, the two predictors $\text{APAR}_{\text{chl}}$ and SW. A possible explanation is that the true chlorophyll content is unknown, and that $\text{APAR}_{\text{chl}}$ therefore itself is estimated (on the basis of the greenness index EVI \citep{huete2002overview}). Due to these imperfect measurements, $X^1$ may still contain information about $Y$ that cannot be explained by $X^2$. 

\subsubsection{Reconstruction of the Vegetation Type} \label{sec:reconstruction}
\begin{figure}
	\centering 
	\includegraphics[width=\linewidth]{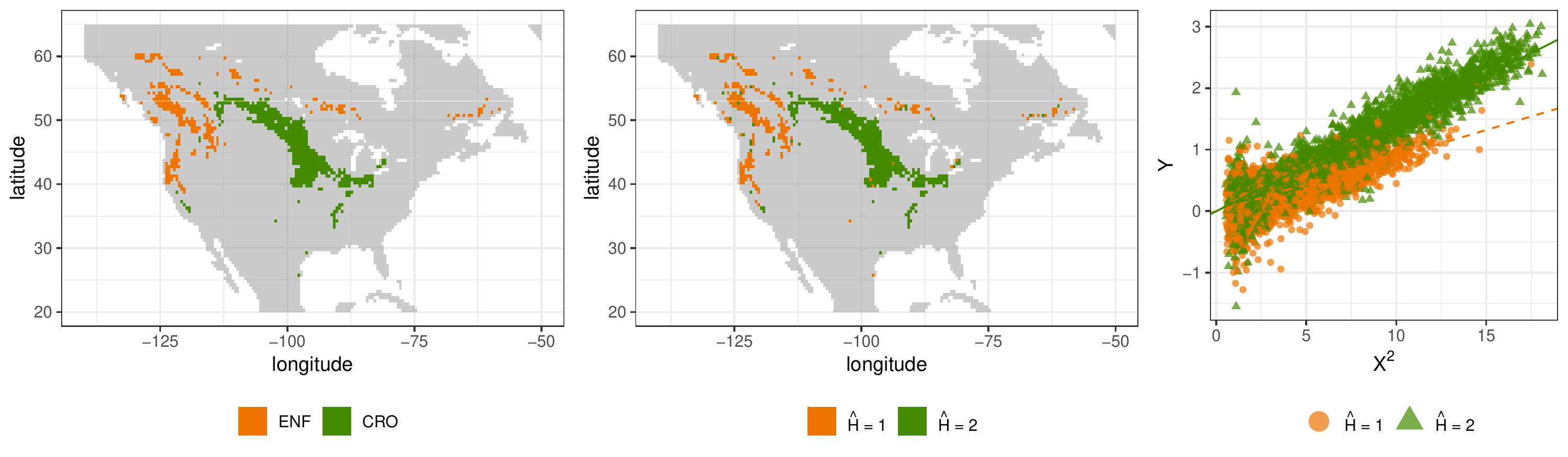}
	\caption{
Vegetation type by IGBP (left) and estimates obtained from reconstructing the values of the hidden variable, as described in Section~\ref{sec:reconstruction} (middle).  
We correctly classify more than 95\% of the pixels. 	The right hand plot illustrates the vegetation-dependent linear relationship between $Y$ and $X^2$. Switching regression model fits are indicated by straight lines, and points are colored according the reconstructed value of $\hat H$. Since the data are not well-clustered in the $X^2$-$Y$ space, classifying observations based on data from $(Y, X^2)$ is generally not a straight-forward task.}
	\label{fig:reconstruction}
\end{figure}

We know that $(Y,X^2)$ follows a switching regression model (see Figure~\ref{fig:data}), 
and that the hidden variable in this model corresponds to the true vegetation type. 
We can thus obtain estimates of the vegetation type by reconstructing the values of the hidden variable in the fitted model. 
We use the data set at its highest resolution, %
and exploit
the background knowledge that $H$ does not change throughout the considered time span. All observations obtained from one spatial grid cell are therefore assumed to stem from the same underlying regime. 
Let $\mathcal{S} \subset \R^2$ and $\mathcal{T} = \{1, \dots, 12\}$ be the spatial and the temporal grid, respectively, along which data are observed.
We then classify each grid cell $s \in \mathcal{S}$ as 
$\hat{H}_s := \argmax_{j \in \{1,2\}} \sum_{t \in \mathcal{T}} \hat \P(H_{st} = j \given Y_{st}, X_{st})$, where $\hat \P$ refers to the fitted model. 
Our method correctly reconstructs the hidden variable in more than 95\% of the grid cells (Figure~\ref{fig:reconstruction}, left and middle). 
As seen in Figure~\ref{fig:reconstruction} (right), reconstructing $H$ based on data from $(Y, X^2)$ is not an easy classification problem. 

So far, we have assumed that the IGBP classification corresponds to the true vegetation type. 
In reality, it is an estimate based on greenness indices that are constructed from remotely sensed radiation reflected from the Earth's surface. 
The outcome of our method may be viewed as an alternative ecosystem classification scheme, 
which additionally comes with a process-based interpretation: each cluster corresponds to a different 
slope parameter in the linear regression of SIF on $\text{APAR}_{\text{chl}}$. This parameter represents the efficiency 
at which absorbed energy is quenched as fluorescence, and is referred to as \textit{fluorescence yield}.

\section{Conclusions and Future Work} \label{sec:concl}
This paper discusses methodology for causal discovery that is applicable in the presence of discrete hidden variables. If the data set is time-ordered, 
the hidden variables may follow a Markov structure.
The method is formulated in the framework of invariant causal prediction. It aims at inferring causal predictors of a target variable and comes with the following coverage guarantee: whenever the method's output is non-empty, it is correct with large probability.
Our algorithm allows 
for several user choices 
and is tested on 
a wide range of simulations.
We see that also in small sample regimes and under a variety of different model violations,
the coverage is not negatively affected. 
Our implementation allows for using either the EM-algorithm or a numerical maximization technique. 
In our experiments, we find that the two options yield very similar results, 
but that the latter is computationally faster and more suitable for handling parameter constraints. 
The power of both methods decreases with an increasing number of hidden states. 
This conforms to the theoretical result that, in general, identifiability of causal predictors cannot be achieved if the hidden variable may take arbitrarily many states, for example.

As part of the method, 
we propose a test for the equality of two switching regression models; 
to the best of our knowledge this is the first example of such a test and may be of interest in itself. 
We prove the asymptotic validity of this test by providing sufficient conditions for the existence, the consistency and the asymptotic normality of the maximum likelihood estimator in switching regression models. 

On the real world data, the true causal parent is consistently attributed the highest significance as a causal predictor of the target variable. 
Switching regression models can also be used for classifying data points based on reconstructed values of the hidden variables.

For large sample sizes, most goodness of fits test are usually rejected in real data.
Since the $h$-invariance assumption may not hold exactly either, 
it may be interesting to explore %
relaxations of this assumption. %
For example, \citet{pfister2018identifying} propose a causal ranking, and \citet{rothenhausler2018anchor} interpolate between prediction and invariance. 
Our robustness analysis in Section~\ref{sec:sensitivity} suggests that 
the performance of our method is not negatively affected when allowing for a dependence
between $X$ and $H$, and we believe that our theoretical results could be extended 
to 
such scenarios (possibly adding mild assumptions).
To widen the range of applicability of our method, it might also be worthwhile to consider  
non-linear models. In particular, it would be interesting to construct conditional 
independence tests that are able to take into account a mixture model structure.

\acks{}%
We thank Roland Langrock 
for insightful comments and 
providing parts of the code;  
Jens Ledet Jensen, 
Miguel Mahecha and
Markus Reichstein
for helpful discussions;
and Yao Zhang and Xiangming Xiao for providing parts of the data used in Section~\ref{sec:SIF}. 
We thank two anonymous referees and the AE for many helpful and constructive comments.
This research was supported by a research grant (18968) from VILLUM FONDEN.

\appendix

\section{Structural Causal Models} \label{app:SCM}

Below, we formally define structural causal models \citep{Pearl2009, bollen2014structural}, and use a presentation similar to \citet[Chapter~6]{Peters2017book}. 

\begin{defi}[Structural causal model] \label{defi:SCM}
A \emph{structural causal model} (SCM) over variables $(Z_1, \dots, Z_p)$ consists of a family of structural assignments
\begin{equation*}
Z_j := f_j(\B{PA}_j, N_j), \qquad j = 1, \dots, p,
\end{equation*}
where for each $j \in \{1, \dots, p\}$, $\B{PA}_j \subset \{Z_1, \dots, Z_p \}  \setminus \{ Z_j\}$ denotes the \emph{parent set} of variable $Z_j$, and a product distribution over the noise variables $(N_1, \dots, N_p)$. Every SCM induces a graph over the nodes in $\{Z_1, \dots, Z_p\}$: for every $j$, one draws an arrow from each of the variables in $\B{PA}_j$ to $Z_j$. We here require this graph to be acyclic. A variable $Z_i$ is a \emph{cause} of $Z_j$, if there exists a directed path from $Z_i$ to $Z_j$. The variables in $\B{PA}_j$ are said to be the \emph{direct causes} of $Z_j$. 
\end{defi}
Due to the acyclicity of the graph, an SCM induces a distribution over the variables $Z_1, \ldots, Z_p$.
An \emph{intervention} on $Z_j$ corresponds to replacing the corresponding assignment. 
(We still require joint independence of all noise variables, as well as the acyclicity of the induced graph to be preserved under interventions.)
This yields another SCM and another distribution, the intervention distribution.

\section{Parametrizations of the Models IID and HMM} \label{app:para}
Define $\mathcal{G}^{\text{IID}} := [0,1]^{\ell-1}$ and 
$\mathcal{G}^{\text{HMM}} := \{ \gamma \in [0,1]^{(\ell-1)\ell} \given \text{ for all } j \in \{1, \dots, \ell \} : \sum_{k = 1}^{\ell-1} \gamma_{j \ell + k} \leq 1 \}$ and parametrize the transition matrix via the maps $\B{\Gamma}^{\text{IID}} : \mathcal{G}^{\text{IID}} \to [0,1]^{\ell \times \ell}$ and $\B{\Gamma}^{\text{HMM}} : \mathcal{G}^{\text{HMM}} \to [0,1]^{\ell \times \ell}$, for all $i,j \in \{1, \dots, \ell\}$ given by 
\begin{equation*}
\B{\Gamma}_{ij}^{\text{IID}}(\gamma) = \begin{cases} \gamma_j & j < \ell \\ 1-\sum_{k=1}^{\ell-1} \gamma_k & j = \ell \end{cases}
\qquad \text{and} \qquad
\B{\Gamma}_{ij}^{\text{HMM}}( \gamma) = \begin{cases} \gamma_{i \ell + j} & j < \ell \\ 1-\sum_{k=1}^{\ell-1} \gamma_{i \ell+k} & j = \ell. \end{cases}
\end{equation*}
For the regression matrix $\Theta$, we consider the two types of parameter constraints discussed in Section~\ref{sec:likelihood}. For $c>0$, let $\mathcal{T}^c := (\R^p \times [c, \infty))^\ell$ and $\mathcal{T}^{=} := \R^{p \ell} \times (0, \infty)$ and parametrize the regression matrix via the maps $\B{\Theta}^c : \mathcal{T}^c \to \R^{p \times \ell}$ and $\B{\Theta}^= : \mathcal{T}^= \to \R^{p \times \ell}$, for all $i \in \{1, \dots, p+1\}$ and $j \in \{1, \dots, \ell\}$ given by
\begin{equation*}
\B{\Theta}^c_{ij}(\theta) = \theta_{(j-1) (p+1) + i} \qquad \text{ and } \qquad
\B{\Theta}^=_{ij}(\theta) = \begin{cases} \theta_{(j-1)p + i} & i \leq p \\ \theta_{p \ell +1} & i = p+1. \end{cases}
\end{equation*}
Both of the parameter constraints induced by $(\B{\Theta}^c, \mathcal{T}^c)$ and $(\B{\Theta}^=, \mathcal{T}^=)$ ensure the existence of the maximum likelihood estimator, see Theorem~\ref{thm:existence}. Since all of the above coordinate mappings are linear in $\theta$ and $\gamma$, Assumption~(A4) in Section~\ref{sec:covgarant_cr} is satisfied for any pair $(\B{\Theta}, \B{\Gamma})$ with $\B{\Theta} \in \{\B{\Theta}^c, \B{\Theta}^=\}$ and $\B{\Gamma} \in \{\B{\Gamma}^{\text{IID}}, \B{\Gamma}^{\text{HMM}}\}$.

\section{Proofs} \label{app:proofs}

\subsection{Proof of Proposition~\ref{prop:causalS}}
Recall that by Definition~\ref{defi:SCM}, we require the underlying causal graph to be acyclic. For every $t \in \{1, \dots, n\}$, 
we can therefore recursively substitute structural assignments to express $(X_t^{\PA{Y}{0}}, H_t^*)$ as a function of all
noise variables appearing in the structural assignments of the ancestors of $Y_t$. Using the joint independence of all noise variables
(see~Definition~\ref{defi:SCM}), it follows that $(X_t^{\PA{Y}{0}}, H_t^*) \indep N_t$. Using the i.i.d.\ assumption on $(N_t)_{t \in \{1, \dots, n\}}$,
we have that for all $t$ and for all $x,h$, the distribution of $Y_t \given (X_t^{\PA{Y}{0}} = x, H^*_t = h) \stackrel{d}{=} f(x, h, N_t)$ 
does not depend on $t$, which shows that $S^*={\PA{Y}{0}}$ satisfies \eqref{eq:SstarH}. By writing 
$Y_t = \sum_{h = 1}^\ell f(X_t^{\PA{Y}{0}}, h, N_t) \mathbbm{1}_{\{H_t^*=h\}}$ and using the linearity of 
the functions $f(\cdot, h , \cdot)$, it follows that $S^*={\PA{Y}{0}}$ is $h$-invariant with respect to $(\B{Y}, \B{X})$. $\hfill \blacksquare$

\subsection{Proof of Theorem~\ref{thm:existence}}
We first introduce some notation. Since neither of the parametrizations in question impose any constraints on the regression coefficients, we will throughout this proof write $\theta = (\beta, \delta)$, where $\beta = (\beta_1, \dots, \beta_\ell) \in \mathcal{B} := \R^{p \times \ell}$ and $\delta \in \mathcal{D}$ is the part of $\theta$ that parametrizes the error variances, i.e., $\mathcal{D}^= = (0, \infty)$ and $\mathcal{D}^c = [c, \infty)^\ell$. Also, we will use $\bar{\mathcal{D}}^= = [0, \infty]$, $\bar{\mathcal{D}}^c = [c, \infty]^\ell$, $\bar{\mathcal{B}} = (\R \cup \{- \infty, +\infty\})^{p \times \ell}$ to denote the ``compactifications'' of $\mathcal{D}^c$, $\mathcal{D}^=$ and $\mathcal{B}$, respectively. 
For every $h \in \{1, \dots, \ell\}^m$ and every $j \in \{1, \dots, \ell\}$ define $T_{h=j} := \{t \in \{1, \dots, m\} : h_t = j\}$ 
and write the likelihood function as $G = \sum_{h \in \{1, \dots, \ell\}^m} g_h$, where
\begin{equation*}
g_h(\phi) = p(\B{x}) \lambda(\gamma)_{h_1} \prod_{s=2}^m \B{\Gamma}_{h_{s-1} h_s} (\gamma) \prod_{j = 1}^\ell \prod_{t \in T_{h=j}} \mathcal{N}(y_t \given x_t \beta_{h_t}, \sigma_{h_t}^2(\delta)),
\end{equation*}
where the product over an empty index set is defined to be 1. 

Let $G^* := \sup_{\phi \in \mathcal{P}} G(\phi) \in (0, \infty]$. 
We want to show that there exists $\phi^* \in \mathcal{P}$ with $G(\phi^*) = G^*$ (which in particular shows that $G^* < \infty$). The idea of the proof is as follows. 
We first show that given an arbitrary point $\bar{\phi}$ in the compactification $\bar{\mathcal{P}}$  
and an arbitrary sequence $(\phi^n)_{n \in \N}$ in $\mathcal{P}$ that converges to $\bar{\phi}$, we can construct a sequence $(\tilde{\phi}^n)_{n \in \N}$ with limit point $\tilde{\phi} \in \mathcal{P}$, such that $\lim_{n \to \infty} G(\tilde{\phi}^n) \geq \lim_{n \to \infty} G(\phi^n)$. We then let $(\phi^{*n})_{n \in \N}$ be a sequence with $\lim_{n \to \infty} G(\phi^{*n}) = G^*$. By compactness of $\bar{\mathcal{P}}$, we can wlog assume that $(\phi^{*n})_{n \in \N}$ is convergent in $\bar{\mathcal{P}}$ (otherwise we may choose a convergent subsequence). By the first part of the proof, there exists a sequence $(\tilde{\phi}^{*n})_{n \in \N}$ that is convergent to some $\phi^* \in \mathcal{P}$, and with $\lim_{n \to \infty} G(\tilde{\phi}^{*n}) = G^*$. By continuity of $G$, $G(\phi^*) = G^*$. 

Let $\bar{\phi} = (\bar{\beta}, \bar{\delta}, \gamma) \in \bar{\mathcal{P}}$ and let $(\phi^n)_{n \in \N} = (\beta^n, \delta^n, \gamma^n)_{n \in \N}$ be such that $\lim_{n \to \infty} \phi^n = \bar{\phi}$. If $\bar{\phi} \in \mathcal{P}$, there is nothing to prove. 
Assume therefore $\bar{\phi} \in \bar{\mathcal{P}} \setminus \mathcal{P}$. Since $\mathcal{G}$ was assumed to be compact, $\bar{\mathcal{P}} = \bar{\mathcal{B}} \times \bar{\mathcal{D}} \times \mathcal{G}$. 
The problem can therefore be divided into the two cases $\bar{\delta} \in \bar{\mathcal{D}} \setminus \mathcal{D}$ and $\bar \beta \in \bar{\mathcal{B}} \setminus \mathcal{B}$, which are treated in Lemma~\ref{lem:delta} and Lemma~\ref{lem:beta}, respectively. 
Together, they imply the existence of a sequence $(\tilde{\phi}^n)_{n \in \N}$ with $\lim_{n \to \infty} \tilde{\phi}^n \in \mathcal{P}$ and $\lim_{n \to \infty} G(\tilde{\phi}^n) \geq \lim_{n \to \infty} G(\phi^n)$, thereby completing the proof of Theorem~\ref{thm:existence}.

We first consider the case where $\bar{\delta} \in \bar{\mathcal{D}} \setminus \mathcal{D}$. 
\begin{lem} \label{lem:delta}
	Let $(\phi^n)_{n \in \N}$ be a sequence in $\mathcal{P}$ that converges to a point $\bar{\phi} = (\bar{\beta}, \bar{\delta}, \gamma) \in \bar{\mathcal{B}} \times (\bar{\mathcal{D}} \setminus \mathcal{D}) \times \mathcal{G}$ and assume that the limit $\lim_{n \to \infty} G(\phi^n)$ exists in $[0,\infty]$. Then, there exists a sequence $(\tilde{\phi}^n)_{n \in \N}$ with limit point $(\bar{\beta}, \delta, \gamma) \in \bar{\mathcal{B}} \times  \mathcal{D} \times \mathcal{G}$, such that $\limsup_{n \to \infty} G(\tilde{\phi}^n) \geq \lim_{n \to \infty} G(\phi^n)$.
\end{lem}
\proof
We treat the two parametrizations $(\B{\Theta}^c, \mathcal{T}^c)$ and $(\B{\Theta}^=, \mathcal{T}^=)$ separately. 

If $\mathcal{D} = \mathcal{D}^c$, then $\bar{\mathcal{D}} \setminus \mathcal{D} = \{ (\bar \delta_1, \dots, \bar \delta_\ell) \in [c, \infty]^\ell : \bar \delta_j = \infty \text{ for at least one } j \}$. 
Let $j$ be such that $\bar{\delta}_j = \infty$. Since for every $h \in \{1, \dots, \ell\}^m$, 
\begin{equation}
g_h(\phi^n) \begin{cases} \to 0 \text{ as } n \to \infty & \text{if } T_{h=j} \not = \emptyset  \\ \text{does not depend on } \delta_j^n & \text{otherwise}, \end{cases}
\end{equation}
we can simply substitute $(\delta_j^n)_{n \in \N}$ by the sequence $(\tilde{\delta}_j^n)_{n \in \N}$ that is constantly equal to $c$, to obtain $(\tilde{\phi}^n)_{n \in \N}$ with $\limsup_{n \to \infty} G (\tilde{\phi}^n) \geq \lim_{n \to \infty} G (\phi^n)$. By repeating this procedure for all $j$ with $\bar{\delta}_j = \infty$, we obtain a sequence $(\tilde{\phi}^n)_{n \in \N}$ with $\limsup_{n \to \infty} G(\tilde{\phi}^n) \geq \lim_{n \to \infty} G(\phi^n)$ and such that $\delta = \lim_{n \to \infty} \delta^n \in \mathcal{D}$. 

If $\mathcal{D} = \mathcal{D}^=$, then $\bar{\mathcal{D}} \setminus \mathcal{D} = \{0, \infty\}$. If $\bar{\delta} = \infty$, then $\lim_{n \to \infty} G(\phi^n) = 0$ and the result is trivial. 
Assume therefore that $\bar{\delta} = 0$. Let $h \in \{1, \dots, \ell\}^m$ be fixed. By the assumption on the sample $(\B{y}, \B{x})$, there exists no set of parameters that yield a perfect fit. We may therefore find a sequence $(s(n))_{n \in \N}$ of elements in $\{1, \dots, m\}$ such that $y_{s(n)} - x_{s(n)} \beta_{h_{s(n)}}^n$ is bounded away from zero for all $n$ large enough. For every $n \in \N$ we have
\begin{equation*}
g_h(\phi^n) \leq p(\B{x}) (2 \pi \sigma_1^2(\delta^n))^{-m/2} \exp \left( -\dfrac{1}{2 \sigma_1^2(\delta^n)} (y_{s(n)} - x_{s(n)} \beta^n_{h_{s(n)}})^2 \right).
\end{equation*}
Since the last factor on the right hand side goes to zero exponentially fast in $\sigma_1^2(\delta^n)$, it follows that $\lim_{n \to \infty} g_h(\phi^n) = 0$. Since $h$ was arbitrary, we have that $\lim_{n \to \infty} G(\phi^n) = 0$, and the result follows. %
\endproof

We now turn to the case where $\bar{\beta} \in \bar{\mathcal{B}} \setminus \mathcal{B}$.
\begin{lem} \label{lem:beta}
	Let $(\phi^n)_{n \in \N}$ be a sequence in $\mathcal{P}$ that converges to a point $\bar{\phi} = (\bar{\beta}, \delta, \gamma) \in (\bar{\mathcal{B}} \setminus \mathcal{B}) \times \mathcal{D} \times \mathcal{G}$. Then, there exists a sequence $(\tilde{\phi}^n)_{n \in \N}$ with limit point $(\beta, \delta, \gamma) \in \mathcal{B} \times \mathcal{D} \times \mathcal{G}$, such that $\lim_{n \to \infty} G(\tilde{\phi}^n) \geq \limsup_{n \to \infty} G(\phi^n)$.
\end{lem}
\proof
The idea of the proof is as follows. We construct a bounded sequence $(\tilde{\beta}^n)_{n \in \N}$, such that the sequence $(\tilde{\phi}^n)_{n \in \N}$ obtained from $(\phi^n)_{n \in \N}$ by substituting 
$(\beta^n)_{n \in \N}$ by $(\tilde{\beta}^n)_{n \in \N}$ satisfies that $\lim_{n \to \infty} G(\tilde{\phi}^n) \geq \limsup_{n \to \infty} G(\phi^n)$. 
Since $(\delta^n)_{n \in \N}$ was assumed to be convergent in $\mathcal{D}$ (and hence bounded) and by compactness of $\mathcal{G}$, the whole sequence $(\tilde{\phi}^n)_{n \in \N}$ is bounded. We can therefore find a compact set $\mathcal{K} \subseteq \mathcal{P}$, such that $\{\tilde{\phi}^n : n \in \N \} \subseteq \mathcal{K}$. Consequently, we can wlog assume that $(\tilde{\phi}^n)_{n \in \N}$ is convergent in $\mathcal{K}$ (otherwise we may choose a convergent subsequence). The sequence $(\tilde{\phi}^n)_{n \in \N}$ then fulfills the requirements in Lemma~\ref{lem:beta}, thereby completing the proof. 

The crucial part that remains is the construction of the sequence $(\tilde{\beta}^n)_{n \in \N}$. This is done by induction. Let $(\phi^n)_{n \in \N} = (\beta_1^n, \dots, \beta_\ell^n, \delta^n, \gamma^n)$ be as stated in Lemma~\ref{lem:beta} and let $K^\infty$ be the set of states $k$, for which $\norm{\beta_k^n} \to \infty$ as $n \to \infty$. We then construct $(\tilde{\beta}^n)_{n \in \N}$ in the following way. Pick an arbitrary
$k \in K^\infty$ 
and construct a bounded sequence $(\tilde{\beta}^n_k)_{n \in \N}$ (this construction is described below), such that the sequence $(\tilde{\phi}_{(k)}^n)_{n \in \N}$ obtained from $(\phi^n)_{n \in \N}$ by substituting 
$(\beta^n_k)_{n \in \N}$ by $(\tilde{\beta}^n_k)_{n \in \N}$ satisfies that $\limsup_{n \to \infty} G(\tilde{\phi}_{(k)}^n) \geq \limsup_{n \to \infty} G(\phi^n)$.
We then take $k^\prime \in K^\infty \setminus \{k\}$ and similarly construct 
$(\tilde{\phi}^n_{(k, k^\prime)})_{n \in \N}$ from $(\tilde{\phi}_{(k)}^n)_{n \in \N}$ such that $\limsup_{n \to \infty} G(\tilde{\phi}^n_{(k, k^\prime)}) \geq \limsup_{n \to \infty} G(\tilde{\phi}_{(k)}^n)$. 
By inductively repeating this procedure for all elements of $K^\infty$, we obtain a bounded sequence $(\tilde{\beta}^n)_{n \in \N}$, such that $(\tilde{\phi}^n)_{n \in \N} = (\tilde{\beta}^n, \delta^n, \gamma^n)_{n \in \N}$ satisfies that $\limsup_{n \to \infty} G(\tilde{\phi}^n) \geq \lim_{n \to \infty} G(\phi^n)$. Once again, we can wlog assume that $(G(\tilde{\phi}^n))_{n \in \N}$ converges, since otherwise we can choose a convergent subsequence $(G(\tilde{\phi}^{n_i}))_{i \in \N}$ with $\lim_{i \to \infty} G(\tilde{\phi}^{n_i}) = \limsup_{n \to \infty} G(\tilde{\phi}^n)$. 

We now prove the induction step. 
Assume that we have iteratively constructed sequences for $k_1, \dots, k_j \in K^\infty$ (if $j=0$, this corresponds to the base case). For simplicity write $(\check{\phi}^n)_{n \in \N} = (\tilde{\phi}^n_{(k_1, \dots, k_j)})_{n \in \N} $. Pick an arbitrary $k \in K^\infty \setminus \{k_1, \dots, k_j\}$. If for all $t \in \{1, \dots, m\}$, $\card{x_t \beta_k^n} \to \infty$ as $n \to \infty$, we could (similar to the proof of Lemma~\ref{lem:delta}) take $(\tilde{\beta}_k^n)_{n \in \N}$ to be a constant sequence. Since in general, there might exist $s$ such that $\card{x_s \beta_k^n} \not \to \infty$ as $n \to \infty$, we divide the problem as follows. 
Define $\mathcal{S}_1 := \{s \in \{1, \dots, m\} : \card{x_s \beta_k^n} \to \infty \text{ as } n \to \infty \}$, 
$\mathcal{S}_2 := \{1, \dots, m\} \setminus \mathcal{S}_1$, $\mathcal{H}_1 := \{h \in \{1, \dots, \ell \}^m: T_{h=k} \cap \mathcal{S}_1 \not = \emptyset \}$ and $\mathcal{H}_2 := \{1, \dots, \ell\}^m \setminus \mathcal{H}_1$, and write the likelihood function as $G = G_1 + G_2$, where $G_1 := \sum_{h \in \mathcal{H}_1} g_h$ and $G_2 := \sum_{h  \in \mathcal{H}_2} g_h$. We now show that $\lim_{n \to \infty} G_1(\check{\phi}^n) = 0$. We formulate a slightly more general result, which we will also make use of later in the proof:
\begin{itemize}
	\item [(*)] Let $h \in \{1, \dots, \ell\}^m$ and assume there exists a sequence $(s(n))_{n \in \N}$ of elements in $T_{h=k}$, such that $\card{x_{s(n)} \beta_k^n} \to \infty$ as $n \to \infty$. Then, $\lim_{n \to \infty} g_h(\phi^n) = 0$. 
\end{itemize}
\begin{proof}[*] 
	Since $(\delta^n)_{n \in \N}$ was assumed to be convergent in $\mathcal{D}$, all sequences $\{\sigma_j^2(\delta^n)\}_{n \in \N}$, $j \in \{1, \dots, \ell\}$, are bounded from above and bounded away from $0$. Since for all $n \in \N$,
	\begin{equation*}
	g_h(\phi^n) \leq p(x) (2 \pi)^{-n/2} \prod_{t=1}^m (\sigma_{h_t}^2(\delta^n))^{-1/2} \exp \underbrace{\left(- \frac{1}{2 \sigma_{k}^2(\delta^n)} (y_{s(n)} - x_{s(n)} \beta_k^n)^2 \right)}_{\to - \infty},
	\end{equation*}
	we are done. %
\end{proof}

For $h \in \mathcal{H}_1$, we can simply pick $s_0 \in T_{h=k} \cap \mathcal{S}_1$ and consider the sequence $(s(n))_{n \in \N}$ that is constantly equal to $s_0$. The result (*) therefore shows that $\lim_{n \to \infty} G_1(\check{\phi}^n) = 0$. 
It thus suffices to construct $(\tilde{\phi}^n_k)_{n \in \N}$ from $(\check{\phi}^n)_{n \in \N}$ such that $\limsup_{n \to \infty} G_2(\tilde{\phi}^n_k) \geq \limsup_{n \to \infty} G_2(\check{\phi}^n)$. Since for every $h \in \mathcal{H}_2$ we have $T_{h=k} \subseteq \mathcal{S}_2$, we take a closer look at $\mathcal{S}_2$. For every $s \in \mathcal{S}_2$, the sequence $(\card{x_s \beta_k^n})_{n \in \N}$ is either bounded or can be decomposed into two sequences, one of which is bounded and one of which converges to infinity. For every $s \in \mathcal{S}_2$, let therefore $I_s^b$ and $I_s^\infty$ be disjoint subsets of $\N$ with $I_s^b \cup I_s^\infty = \N$, such that $(\card{x_s \beta_k^n})_{n \in I_s^b}$ is bounded and such that either $I_s^\infty = \emptyset$ or $\card{I_s^\infty} = \infty$ with $(\card{x_s \beta_k^n})_{n \in I_s^\infty}$ converging to infinity. 
Let $I^b := \cup_{s \in \mathcal{S}_2} I_s^b$ and define a sequence $(\tilde{\beta}^n_k)_{n \in \N}$ by 
\begin{equation*}
\tilde{\beta}^n_k := \begin{cases} \text{ the projection of } \beta^n_k \text{ onto } \text{span}_\R(\{x_s : s \text{ satisfies } n \in I_s^b\}) & \text{ if } n \in I^b \\ 
0 & \text{ otherwise}. \end{cases}
\end{equation*}
We now show that the above defines a bounded sequence.
\begin{itemize}
	\item [($\circ$)] The sequence $(\tilde{\beta}_k^n)_{n \in \N}$ is bounded.
\end{itemize}
\begin{proof}[$\circ$]
	For every $\mathcal{S} \subseteq \mathcal{S}_2$, define $I^b_\mathcal{S} := \{n \in \N: n \in I_s^b \Leftrightarrow s \in \mathcal{S} \}$ (where $I^b_\emptyset := \N \setminus I^b$). We can then decompose $(\tilde{\beta}^n_k)_{n \in \N}$ into the subsequences $(\tilde{\beta}^n_k)_{n \in I^b_\mathcal{S}}$, $\mathcal{S} \subseteq \mathcal{S}_2$, and prove that each of these sequences is bounded.
	Let $\mathcal{S} \subseteq \mathcal{S}_2$ and let $\{u_1, \dots, u_d \}$ be an orthonormal basis for $\text{span}_\R (\{ x_s : s \in \mathcal{S}\})$. 
	Since all sequences in $\{ (\card{x_s \tilde{\beta}^n_k})_{n \in I^b_\mathcal{S}} : s \in \mathcal{S} \}$ are bounded, then so are the sequences $(\card{u_1 \tilde{\beta}^n_k})_{n \in I^b_\mathcal{S}}, \dots, (\card{u_d \tilde{\beta}^n_k})_{n \in I^b_\mathcal{S}}$ (this follows by expressing each of the $u_i$s as a linear combination of elements in $\{x_s : s \in \mathcal{S}\}$). The result now follows from the identities $\norm{\tilde{\beta}^n_k}^2 = \sum_{j=1}^d \card{u_j \tilde{\beta}^n_k}^2$, $n \in I^b_\mathcal{S}$. 
\end{proof}

Let $(\tilde{\phi}^n_k)_{n \in \N}$ be the sequence obtained from $(\check{\phi}^n)_{n \in \N}$ by substituting $(\beta_k^n)_{n \in \N}$ by $(\tilde{\beta}_k^n)_{n \in \N}$. Finally, we show the following result. 
\begin{itemize}
	\item [($\triangle$)] $\limsup_{n \to \infty} G(\tilde{\phi}_k^n) \geq \limsup_{n \to \infty} G(\check{\phi}^n)$. 
\end{itemize}
\begin{proof}[$\triangle$]
	Let $h \in \mathcal{H}_2$ and define $I_h^\infty := \bigcup_{s \in T_{h=k}} I_s^\infty$ (if $T_{h=k} = \emptyset$, we define $I_h^\infty:= \emptyset$). 
	The idea is to decompose $(\check{\phi}^n)_{n \in \N}$ into $(\check{\phi}^n)_{n \in I_h^\infty}$ and $(\check{\phi}^n)_{n \not \in I_h^\infty}$ and to treat both sequences separately. 
	
	We start by considering $(\check{\phi}^n)_{n \not \in I_h^\infty}$. 
	First, observe that for every $s$, $\mathcal{N}(y_s \given x_s \beta_k, \sigma_k^2 (\delta))$ only depends on $\beta_k$ via the inner product $x_s \beta_k$. 
	By construction of $I_h^\infty$ and $(\tilde{\beta}_k^n)_{n \in \N}$, we thus have that for all $n \not \in I_h^\infty$ and for all $s \in T_{h=k}$, 
	the function values $\mathcal{N}(y_s \given x_s \tilde{\beta}_k^n, \sigma_k^2 (\delta^n))$ and $\mathcal{N}(y_s \given x_s \beta_k^n, \sigma_k^2 (\delta^n))$ coincide. Consequently, we have that for all $n \not \in I_h^\infty$, $g_h(\tilde{\phi}_k^n) = g_h(\check{\phi}^n)$. In particular, the sequences $(\check{g}^n_{h, b})_{n \in \N}$ and $(\tilde{g}^n_{h, b})$, for every $n \in \N$ defined by $\check{g}^n_{h, b} := g_h(\check{\phi}^n) \mathbbm{1}_{\{n \not \in I_h^\infty\}}$ and $\tilde{g}^n_{h,b} := g_h(\tilde{\phi}^n) \mathbbm{1}_{\{n \not \in I_h^\infty\}}$, coincide. 
	
	We now consider $(\check{\phi}^n)_{n \in I_h^\infty}$. By construction of the sets $I_s^\infty, s \in T_{h=k}$, either $I_h^\infty = \emptyset$ or $\card{I_h^\infty} = \infty$. If $\card{I_h^\infty} = \infty$, then for every $n \in \N$, there exists $\check{s}(n) \in T_{h=k}$ such that $n \in I_{\check{s}(n)}^\infty$. By applying (*) to the sequence $(\check{\phi}^n)_{n \in I^\infty}$ with $(s(n))_{n \in I_h^\infty} = (\check{s}(n))_{n \in I_h^\infty}$, it follows that $\lim_{n \to \infty, n \in I_h^\infty} g_h(\check{\phi}^n) = 0$. In particular, the sequences $(\check{g}^n_{h, \infty})_{n \in \N}$ and $(\tilde{g}^n_{h, \infty})_{n \in \N}$, for every $n \in \N$ defined by $\check{g}^n_{h, \infty} := g_h(\check{\phi}^n) \mathbbm{1}_{\{n \in I_h^\infty\}}$ and $\tilde{g}^n_{h,\infty} := g_h(\tilde{\phi}^n) \mathbbm{1}_{\{n \in I_h^\infty\}}$, converge to 0 as $n \to \infty$ (this holds also if $I^\infty = \emptyset$). 
	
	By combing the above results for all $h \in \mathcal{H}_2$, we finally have 
	\begin{align*}
	\limsup_{n \to \infty} G_2(\check{\phi}^n) 
	&= \limsup_{n \to \infty} \left( \sum_{h \in \mathcal{H}_2} \check{g}^n_{h,b} + \sum_{h \in \mathcal{H}_2} \check{g}^n_{h,\infty} \right) = \limsup_{n \to \infty} \left( \sum_{h \in \mathcal{H}_2} \check{g}^n_{h,b} \right) \\ &= \limsup_{n \to \infty} \left( \sum_{h \in \mathcal{H}_2} \tilde{g}^n_{h,b} \right) 
	\leq \limsup_{n \to \infty} \left( \sum_{h \in \mathcal{H}_2} \tilde{g}^n_{h,b} +\sum_{h \in \mathcal{H}_2} \tilde{g}^n_{h,\infty} \right) \\
	&= \limsup_{n \to \infty} G_2(\tilde{\phi}_k^n).
	\end{align*}
	Since $\limsup_{n \to \infty} G_1 (\tilde{\phi}^n_k) \geq 0 = \limsup_{n \to \infty} G_1(\check{\phi}^n)$, the result follows. %
\end{proof}

This completes the proof of Lemma~\ref{lem:beta}. 
\endproof

\subsection{Proof of Theorem~\ref{thm:consistency}} \label{sec:proof_cons}
We start by introducing some notation to be used in the proofs of Theorem~\ref{thm:consistency} and Theorem~\ref{thm:normality}. 
Let $\mathcal{K} := \R^p \times (0, \infty)$ be the full parameter space for a single pair $\kappa = (\beta^T, \sigma^2)^T$ of regression parameters. In analogy to previous notation, we will use $\kappa_j(\theta)$ to denote the $j$th pair of regression parameters of a parameter vector $\theta \in \mathcal{T}$. If the conditional distribution of $Y_t \given (X_t = x, H_t = j)$ is a normal distribution with regression parameters $\kappa$, we will denote the conditional density of $(X_t, Y_t) \given (H_t=j)$ by $f(x,y \given \kappa)$. We use $\P_{0}$ for the distribution $\mathcal{SR}(\phi^0 \given X_1)$ and $\E_0$ for the expectation with respect to $\P_0$. Finally, for every $k \in \N$, let $\mathcal{SR}^k ( \cdot \given X_1)$ denote the unconstrained class of mixture distributions of degree $k$ (i.e., all parameters can vary independently within their range). 

Theorem~\ref{thm:consistency} now follows from \citet[Theorem~3]{leroux1992maximum}. To prove the applicability of their result, we first state slightly adapted versions of their conditions (L1)--(L6) and prove afterwards that they are satisfied.
(L1) $\Gamma^0$ is irreducible, (L2) for each $(x,y)$, $\kappa \mapsto f(x,y \given \kappa)$ is continuous and vanishes at infinity (see the last paragraph of Section 2 in \citet{leroux1992maximum}), (L3) for all $j,k \in \{1, \dots, \ell\}$, the maps $\theta \mapsto \kappa_j (\theta)$ and $\gamma \mapsto \B{\Gamma}_{jk} (\gamma)$ are continuous, (L4) for all $j \in \{0, \dots, \ell \}$, $\E_{0} [\card{\log f(X_1, Y_1 \given \kappa_{j}(\theta^0))} ] < \infty$, (L5) for all $\kappa \in \mathcal{K}$, there exists a $\delta > 0$ such that $\E_{0}[\sup_{\kappa^{\prime}: \norm{\kappa - \kappa^{\prime}} < \delta} (\log f(X_1, Y_1 \given \kappa^{\prime}))^+] < \infty$, and (L6) for every $k \in \{1, \dots, \ell\}$, the class $\mathcal{SR}^k (\cdot \given X_1)$ satisfies the following identifiability property. Define
\begin{align*}
\Lambda^k &:= \{(\lambda_1, \dots, \lambda_k) : \sum_{j=1}^k \lambda_j = 1 \}, \text{ and} \\
\mathcal{Q}^k &:= \left\lbrace \{ (\lambda_1, \kappa_1), \dots, (\lambda_k, \kappa_k) \} \; : \;
\begin{tabular}{@{}l@{}}
$(\lambda_1, \dots, \lambda_k) \in \Lambda^k \text{ and }  \kappa_1, \dots, \kappa_k \in \mathcal{K}$ \\  
$\text{with all } \kappa_j \text{s being distinct}$
\end{tabular}
\right\rbrace
\end{align*}
and consider the mapping $\varphi^k : \mathcal{Q}^k \to \mathcal{SR}(\cdot \given X_1)$ that sends $q = \{(\lambda_1, \kappa_1), \dots, (\lambda_k, \kappa_k) \}$ into the mixture distribution $P_q \in \mathcal{SR}(\cdot \given X_1)$ with density
\begin{equation*}
f_q(x,y) := \sum_{j=1}^k \lambda_j f(x,y \given \kappa_j) = f(x) \sum_{j=1}^k \lambda_j f(y \given x, \kappa_j).
\end{equation*}
Then, for every $k \in \{1, \dots, \ell\}$, $\varphi^k$ is a one-to-one map of $\mathcal{Q}^k$ onto $\mathcal{SR}^k (\cdot \given X_1)$. It is therefore the \textit{set} $\{(\lambda_1, \kappa_1), \dots, (\lambda_k, \kappa_k)\}$, rather than the parameters $(\kappa_1, \dots, \kappa_k)$ and $(\lambda_1, \dots, \lambda_k)$ themselves, that is required to be identifiable. 

We now show that (L1)--(L6) are satisfied.
Condition~(L1) is implied by (A3). Condition~(L2) follows by the continuity of $\kappa \mapsto \mathcal{N}(y \given x, \kappa)$ and (L3) is implied by (A4). For (L4), we see that for all $j \in \{ 0, \dots, \ell \}$, 
\begin{equation*}
\log f(X_1, Y_1 \given \kappa_j(\theta^0)) = \log(2 \pi \sigma_{j}^2(\theta^0)) - \frac{1}{2\sigma_{j}^2(\theta^0)} (Y_1 - X_1 \beta_j(\theta^0))^2 + \log f(X_1) \in \mathcal{L}^1(\P_0),
\end{equation*}
by (A7) and by moment-properties of the normal distribution. For (L5), let $\kappa = (\beta, \sigma^2) \in \mathcal{K}$ and choose $\delta := \sigma^2/2$. We then have
\begin{align*}
\E_{0} \left[ \sup_{\kappa^{\prime}: \norm{\kappa^{\prime} - \kappa} < \delta} (\log f(X_1, Y_1 \vert \kappa^{\prime}))^+\right] 
&\leq \E_{0} \left[ \sup_{\kappa^{\prime}: \norm{\kappa^{\prime} - \kappa} < \delta} (\log f(Y_1 \vert X_1, \kappa^{\prime}))^+ + \card{\log f(X_1)} \right] \\
&\leq  \E_{0} \left[ \sup_{\sigma^{\prime}:\norm{\sigma^{\prime 2} - \sigma^2} < \delta} (-\dfrac{1}{2} \log(2\pi \sigma^{\prime 2}))^+ + \card{\log f(X_1)} \right] \\
&\leq \E_{0} \left[ \dfrac{1}{2} \card{\log(\pi \sigma^2)} + \card{\log f(X_1)} \right] < \infty.
\end{align*}
It is left to prove (L6), the identifiability of the classes $\mathcal{SR}^k (\cdot \given X_1)$. \citet[Proposition~1]{teicher1963identifiability} shows an analogous result for mixtures of univariate normal distributions, that are parametrized by their mean and variance. His result will be the cornerstone for our argument. Consider a fixed $k \in \{1, \dots, \ell\}$, let $q = \{(\lambda_1, \beta_1, \sigma_1^2), \dots, (\lambda_k, \beta_k, \sigma_k^2) \}, q^{\prime} = \{(\lambda^{\prime}_1, \beta^{\prime}_1, {\sigma_1^{\prime}}^2), \dots, (\lambda^{\prime}_k, \beta^{\prime}_k, {\sigma_k^{\prime}}^2) \} \in \mathcal{Q}^k$ and assume that the induced mixtures $P_q$ and $P_{q^\prime}$ are identical. Collect $q$ and $q^{\prime}$ into two matrices $Q, Q^\prime$ with columns $Q_{\cdot j} = (\lambda_j, \sigma_j^2, \beta_j^T)^T$ and $Q^\prime_{\cdot j} = (\lambda^\prime_j, {\sigma_j^\prime}^2, {\beta^\prime_j}^T)^T$ for $j \in \{ 1, \dots, k \}$. We wish to show that $Q$ and $Q^\prime$ are equal up to a permutation of their columns. Because the densities $f_q$ and $f_{q^\prime}$ coincide Lebesgue-almost everywhere, it holds that, for all $x \in \text{int}(\text{supp}(X_1))$,
\begin{equation*}
f_q(y \given x) = \sum_{j=0}^k \lambda_j f(y \given x, \kappa_j) = \sum_{j=0}^k \lambda^\prime_j f(y \given x, \kappa_j^\prime) = f_{q^\prime}(y \given x) \qquad \text{ for almost all } y.
\end{equation*}
It now follows from \citet[Proposition~1]{teicher1963identifiability} that, for all $x \in \text{int}(\text{supp}(X_1))$,
\begin{equation} 
\label{eq:sets}
\{(\lambda_1, \sigma_1^2, x \beta_1), \dots, (\lambda_k, \sigma_k^2, x \beta_k)\} = \{(\lambda^\prime_1, {\sigma^\prime}^2_1, x \beta^\prime_1), \dots, (\lambda^\prime_k, {\sigma^\prime}^2_k, x \beta^\prime_k)\}.
\end{equation}
In the remainder of the proof, we will consider
several $x$ simultaneously (rather than a fixed $x$). This will help us to draw conclusions about the betas.
Equation~\eqref{eq:sets} means that for every $z \in \mathcal{Z} := \R^{2} \times \text{int}(\text{supp}(X_1))$, the vectors $zQ$ and $zQ^\prime$ are equal up to a permutation of their entries (this permutation may depend on $z$). Let $\Sigma$ denote the (finite) family of permutation matrices of size $k \times k$ and consider the partition
\begin{equation*}
\mathcal{Z} = \bigcup_{M \in \Sigma} \mathcal{Z}_M, \quad \text{ where } \quad \mathcal{Z}_M = \{z \in \mathcal{Z} : zQ = z Q^\prime M^T \}.
\end{equation*}
Since  $\mathcal{Z}$ is an open subset of $\R^{p+2}$, there exists an element $M_0 \in \Sigma$, such that $\mathcal{Z}_{M_0}$ contains an open subset of $\R^{p+2}$. We can therefore choose $p+2$ linearly independent elements $z_1, \dots, z_{p+2} \in \mathcal{Z}_{M_0}$ and construct the invertible matrix $\mathbf{Z} = [z_1^T, \dots, z_{p+2}^T]^T$. Since $\mathbf{Z} Q = \mathbf{Z} Q^\prime M_0^T$, it follows that $Q = Q^\prime M_0^T$. $\hfill \blacksquare$

\subsection{Proof of Theorem~\ref{thm:normality}}
Throughout the proof, we make use of the notation introduced in the first paragraph of Appendix~\ref{sec:proof_cons}. Theorem~\ref{thm:normality} follows if both the below statements hold true. 
\begin{itemize}
	\setlength\itemsep{0em}
	\item [(i)] $m^{-1} \mathcal{J}(\hat \phi_m) \to \mathcal{I}_0$ as $m \to \infty$ in $\P_0$-probability. 
	\item [(ii)] $\sqrt{m} (\hat \phi_m - \phi^0) \mathcal{I}_0^{1/2} \stackrel{d}{\longrightarrow} \mathcal{N}(0,I)$ as $m \to \infty$ under $\P_0$. 
\end{itemize}
These results correspond to slightly adapted versions of Lemma~2 and Theorem~1, respectively, in \citet{bickel1998asymptotic} (here referred to as L2 and T1). L2 builds on assumptions (B1)--(B4) to be stated below. T1 additionally assumes that $\phi^0 \in \text{int}(\mathcal{P})$ and that the Fisher information matrix $\mathcal{I}_0$ is positive definite, i.e., our (A2) and (A5). Assumptions~(B1)--(B4) state local regularity conditions for a neighborhood of the true parameter $\phi^0$. We therefore need to verify that there exists an open neighborhood $\mathcal{T}_0$ of $\theta^0$, such that
the following conditions are satisfied. 
\begin{itemize}
	\setlength\itemsep{0em}
	\item [(B1)] The transition matrix $\Gamma^0$ is irreducible and aperiodic.
	\item [(B2)] For all $j,k \in \{1, \dots, \ell\}$ and for all $(x,y)$, the maps $\gamma \mapsto \B{\Gamma}_{jk}(\gamma)$ and $\theta \mapsto f(x,y \vert \kappa_j(\theta))$ 
	(for $\theta \in \mathcal{T}_0$) 
	have two continuous derivatives. 
	\item [(B3)] Write $\theta = (\theta_1, \dots, \theta_K)$. For all $n \in \{1,2\}$, $i_1, \dots, i_n \in \{1, \dots, K\}$ and $j \in \{1, \dots, \ell\}$, it holds that
	\begin{itemize}
		\setlength\itemsep{0em}
		\item [(i)] 
		\begin{equation*}
		\int\sup_{\theta \in \mathcal{T}_0} \left \vert  \frac{\partial^n}{\partial \theta_{i_1} \cdots \partial \theta_{i_n}} f(x, y \vert \kappa_j(\theta)) \right \vert d(x,y) < \infty, \quad \text{and}
		\end{equation*}
		\item [(ii)] 
		\begin{equation*}
		\E_{0} \left[ \sup_{\theta \in \mathcal{T}_0} \left \vert  \frac{\partial^n}{\partial \theta_{i_1} \cdots \partial \theta_{i_n}} \log f(X_1, Y_1 \vert \kappa_j(\theta)) \right \vert^{3-n} \right] < \infty.
		\end{equation*}
	\end{itemize}
	\item [(B4)] For all $(x,y)$, define 
	\begin{equation*}
	\rho(x,y) = \sup_{\theta \in \mathcal{T}_0} \max_{0 \leq i,j \leq \ell} \frac{f(x,y \vert \kappa_i(\theta))}{f(x,y \vert \kappa_j(\theta))}. 
	\end{equation*}
	Then for all $j \in \{1, \dots, \ell\}$, $\P_{0}(\rho(X_1, Y_1) = \infty \given H_1 = j) < 1$.
\end{itemize}

We first construct the set $\mathcal{T}_0$.
Let therefore $\epsilon > 0$ and choose $\mathcal{T}_0$ so small that there exists $c>0$, such that for all $\theta \in \mathcal{T}_0$ and for all $j \in \{1, \dots, \ell\}$ and $k \in \{1, \dots, d\}$, it holds that $\beta_{jk} (\theta) \in (\beta_{jk}(\theta^0) - \epsilon, \beta_{jk}(\theta^0) + \epsilon)$ and $\sigma_j^2(\theta) \geq c$. We can now verify the conditions (B1)--(B4).

Assumption~(B1) is satisfied by (A3). For every $(x,y)$, the maps $\kappa \mapsto f(x,y \vert \kappa)$ are two times continuously differentiable on $\R^p \times (0, \infty)$. Together with (A4), this implies (B2), independently of the choice of $\mathcal{T}_0$. 

For the proof of (B3)(i)--(ii) we will make use of the following result. Let $g$ be a polynomial of $(x,y)$ of degree at most 4, i.e., a sum of functions on the form $b x_i^r x_k^s y^t$ for some $i,k \in \{1, \dots, p\}$ and $r,s,t \in \{0, \dots, 4\}$ with $r+s+t \leq 4$. 
Then, for every $\kappa \in \mathcal{K}$, $\int g(\card{x}, \card{y}) f(x,y \given \kappa) d(x,y) < \infty$, where $\card{x} = (\card{x_1}, \dots, \card{x_p})$. This result follows from the fact that for every $x$, $\int \card{y}^t f(y \given x, \kappa) dy$ is a polynomial of $\card{x}$ of degree $t$, and the assumption that, for all $j \in \{1, \dots, p\}$, $\E[\card{X_1^j}^4] < \infty$.

For (B3)(i), we treat all derivatives simultaneously. Let $n \in \{1, 2\}$, $i_1, \dots, i_n \in \{1, \dots, K\}$ and $j \in \{1, \dots, \ell\}$ be fixed. Let $\{g_\theta\}_{\theta \in \mathcal{T}_0}$ be the functions, for all $(x,y)$ and for all $\theta \in \mathcal{T}_0$ defined by
\begin{equation*}
\frac{\partial^{n}}{\partial \theta_{i_1} \cdots \partial \theta_{i_{n}}} f(x, y \given \kappa_{j}(\theta)) = g_\theta(x, y) \exp \left(- \dfrac{1}{2 \sigma_{j}^2(\theta)} (y - x \beta_{j}(\theta))^2 \right) f(x),
\end{equation*}
(note that $f(x) = 0$ implies $f(x,y \vert \kappa_{j}(\theta)) = 0$).
Then, for all $(x,y)$, $\theta \mapsto g_\theta(x, y)$ is continuous, and for all $\theta \in \mathcal{T}_0$, $(x,y) \mapsto g_\theta(x, y)$ is a polynomial of degree at most 4.
By the compactness of $\bar{\mathcal{T}_0}$, the closure of $\mathcal{T}_0$, and by the continuity of $\theta \mapsto g_\theta(x, y)$, there exists a polynomial $g$ of degree 4, such that, for all $(x,y)$, $\sup_{\theta \in \mathcal{T}_0} \card{g_\theta(x, y)} \leq g(\card{x},\card{y})$.

Consider now a fixed $k \in \{1, \dots, p\}$. By choice of $\mathcal{T}_0$, we have that for all $x_k$ and for all $\theta \in \mathcal{T}_0$, it holds that $x_k (\beta_{jk}(\theta^0) - \text{sign}(x_k) \epsilon) \leq x_k \beta_{jk}(\theta) \leq x_k (\beta_{jk}(\theta^0) + \text{sign}(x_k) \epsilon)$. By writing $s(x) = (\text{sign}(x_1), \dots, \text{sign}(x_p))$ it follows that for all $(x,y)$ and all $\theta \in \mathcal{T}_0$, $y- x (\beta_{j}(\theta^0) - \text{diag}(s(x)) \epsilon) \leq y - x \beta_{j}(\theta) \leq y - x (\beta_{j}(\theta^0) + \text{diag}(s(x)) \epsilon)$. Consequently, we may for every $(x,y)$ find $s(x,y) \in \{-1,1\}^p$ (either $s(x)$ or $-s(x)$) such that for all $\theta\in \mathcal{T}_0$, 
\begin{equation*}
-(y - x \beta_j(\theta))^2 \leq -(y - x ( \underbrace{\beta_j(\theta^0) + \text{diag}(s(x,y)) \epsilon}_{=: \beta_s}))^2.
\end{equation*}
By choosing $C>0$ small enough, it follows that for all $(x,y)$ and for all $\theta \in \mathcal{T}_0$ it holds that
\begin{equation*}
\exp \left( -\dfrac{1}{2 \sigma_j^2(\theta)} (y - x \beta_j(\theta))^2 \right) \leq \exp \left( -C(y - x \beta_j(\theta))^2 \right) \leq 
\sum_{s \in \{-1, 1\}^p} \exp \left( -C(y - x \beta_s)^2 \right).
\end{equation*}
Since all integrals $\int g(\card{x},\card{y}) \exp ( -C(y - x \beta_s)^2 ) f(x) d(x,y)$, $s \in \{-1,1\}^p$, are finite, this completes the proof of (B3)(i). 

The proof of (B3)(ii) is similar to that of (B3)(i). Fix $n \in \{1, 2\}$, $i_1, \dots, i_n \in \{1, \dots, K\}$ and $j \in \{1, \dots, \ell\}$. Let $\{h_\theta\}_{\theta \in \mathcal{T}_0}$ be the functions, for all $(x,y)$ and for all $\theta \in \mathcal{T}_0$ defined by
\begin{equation*}
\frac{\partial^n}{\partial \theta_{i_1} \cdots \partial \theta_{i_n}} \log f(x,y \given \kappa_j(\theta)) = h_\theta(x,y).
\end{equation*}
Then, for all $(x,y)$, $\theta \mapsto h_\theta(x, y)$ is continuous, and for all $\theta \in \mathcal{T}_0$, $(x,y) \mapsto h_\theta(x, y)$ is a polynomial of degree at most 2. We can therefore find a dominating polynomial $h$ of degree 2, such that, for all $(x,y)$, $\sup_{\theta \in \mathcal{T}_0} \card{h_\theta (x,y)} \leq h(\card{x}, \card{y})$. Since $h(\card{X_1}, \card{Y_1}) \in \mathcal{L}^2(\P_0)$, this completes the proof of (B3)(ii). 

(B4) is easily verified. Since the support $\mathcal{S}$ of the functions $f(\cdot \given \kappa)$ does not depend on $\kappa$, it is enough to consider $(x,y) \in \text{int}(\mathcal{S})$. For all $(x,y) \in \text{int}(\mathcal{S})$ and for all $j \in \{1, \dots, \ell \}$, $\theta \mapsto f(x,y \given \kappa_j(\theta))$ is bounded from above and bounded away from zero (by choice of $\mathcal{T}_0$). The function $\rho$ is therefore finite everywhere. $\hfill \blacksquare$

\subsection{Proof of Corollary~\ref{cor:cr_validity}}
Let (A1)--(A7) hold true. By Theorem~\ref{thm:consistency}, we can decompose $(\hat \phi_m)_{m \in \N} = ((\hat \theta_m, \hat \gamma_m))_{m \in \N}$ into one or more subsequences, each of which is convergent to a permutation of $\phi^0$. We can therefore find a sequence $(\pi^m_\mathcal{P})_{m \in \N} = ((\pi^m_\mathcal{T}, \pi^m_\mathcal{G}))_{m \in \N}$ of permutations on $\mathcal{P}$, such that, $\P_0$-almost surely, the sequence of maximum likelihood estimators $( \pi^m_\mathcal{P}(\hat \phi_m) )_{m \in \N}$ converges to $\phi^0$ as $m \to \infty$. For $\alpha \in (0,1)$ and for every $m \in \N$, we then have 
\begin{align*}
\P_0^m (\theta^0 \in C_{\text{adjusted}}^\alpha(\hat \theta_m)) 
\geq \P_0^m (\theta^0 \in C^\alpha (\pi^m_\mathcal{T} (\hat \theta_m))) 
= \P_0^m (\phi^0 \in C^\alpha (\pi^m_\mathcal{T} (\hat \theta_m)) \times \mathcal{G}).
\end{align*}
By Theorem~\ref{thm:normality}, the right hand side converges to $1 - \alpha$ as $m \to \infty$. $\hfill \blacksquare$

\subsection{Proof of Theorem~\ref{thm:asymp_cov}}
By Corollary~\ref{cor:cr_validity}, the adjusted confidence regions within each environment all achieve the correct asymptotic coverage, ensuring the asymptotic validity of the test $\varphi_{S^*}$ of $H_{0,S^*}$. Since, for every $n$, $\P^n_0( \hat S_n \subseteq S^*) \geq \P^n_0(\varphi^n_{S^*} \text{ accepts } H^n_{0,S^*})$, the result follows. $\hfill \blacksquare$

\section{Further Details on Likelihood Optimization} \label{app:lik_opt}
Below, we describe the two optimization methods NLM and EM. Since the loglikelihood function \eqref{eq:loglik} is non-convex, 
the performance of these routines depend on the initialization. 
In practice, we restart the algorithms in 5 different sets of starting values 
(using the \verb|regmix.init| function from the \verb|R| package \verb|mixtools|).

\subsection{Method I (``NLM''): Non-Linear Maximization}
This method maximizes the loglikelihood function \eqref{eq:loglik} numerically. We use the \verb|R| optimizer \verb|nlm|, which is a non-linear maximizer based on a Newton-type optimization routine \citep[e.g.,][]{schnabel1985modular}. 
The method also outputs an estimate of the observed Fisher information, which is used for the construction of the confidence regions~\eqref{eq:Calpha}. An equality constraint on the error variances can be enforced directly by using the parametrization $(\B{\Theta}^=, \mathcal{T}^=)$ described in Appendix~\ref{app:para}. A lower bound (we use $10^{-4}$ as a default value) can be imposed by suitable reparametrization of all error variances \citep[e.g.,][Section~3.3.1]{Zucchini2016}. 

\subsection{Method II (``EM''): The EM-algorithm} \label{app:em}
Given starting values $\phi^{(0)} \in \mathcal{P}$, the EM-algorithm operates by alternating between the following two steps until a convergence criterion is met. (1) The E-step: Compute the posterior distribution $P^{(t)}_{(\B{y}, \B{x})}$ of $\B{H} \given (\B{Y} = \B{y}, \B{X} = \B{x}, \phi^{(t)})$ given the current parameters $\phi^{(t)}$. (2) The M-step: Maximize the expected complete data loglikelihood
\begin{equation} \label{eq:Q}
Q( \phi \given \phi^{(t)}) := \E_{P^{(t)}_{(\B{y}, \B{x})}} \left[ \ell_{\text{complete}} (\B{y}, \B{x}, \B{H} \given \phi) \right]
\end{equation}
to obtain updates $\phi^{(t+1)} \in \argmax_{\phi \in \mathcal{P}} Q( \phi \given \phi^{(t)})$. Here, $\ell_{\text{complete}}$ is the loglikelihood function of the complete data $(\B{y}, \B{x}, \B{h})$.
The explicit forms of $P^{(t)}_{(\B{y}, \B{x})}$ and $Q$ depend on the choice of model. In model IID, $P^{(t)}_{(\B{y}, \B{x})}$ is a product distribution which can be computed by simple applications of Bayes' theorem. In model HMM, the posterior distribution is obtained by the forward-backward algorithm. In both cases, \eqref{eq:Q} can be maximized analytically \citep[e.g.,][Chapters 9 and 13]{bishop2006machine}. 
The observed Fisher information $\mathcal{J}(\hat \phi)$ can be computed analytically from the derivatives of \eqref{eq:Q}, see \citet{oakes1999direct}. 
In our \verb|R| package, the EM-algorithm is only implemented for model IID and makes use of the package \verb|mixreg|. 
An equality constraint on the error variances can be accommodated using the parametrization $(\B{\Theta}^=, \mathcal{T}^=)$ from Appendix~\ref{app:para}. 
A lower bound on the error variances is enforced by restarting the algorithm 
whenever an update $\phi^{(t)}$ contains a variance component that deceeds the lower bound
(\verb|mixreg| uses the threshold $10^{-16}$). 

\begin{figure}
\begin{center}
\begin{minipage}{0.47\textwidth}
\centering
\includegraphics[width = \linewidth]{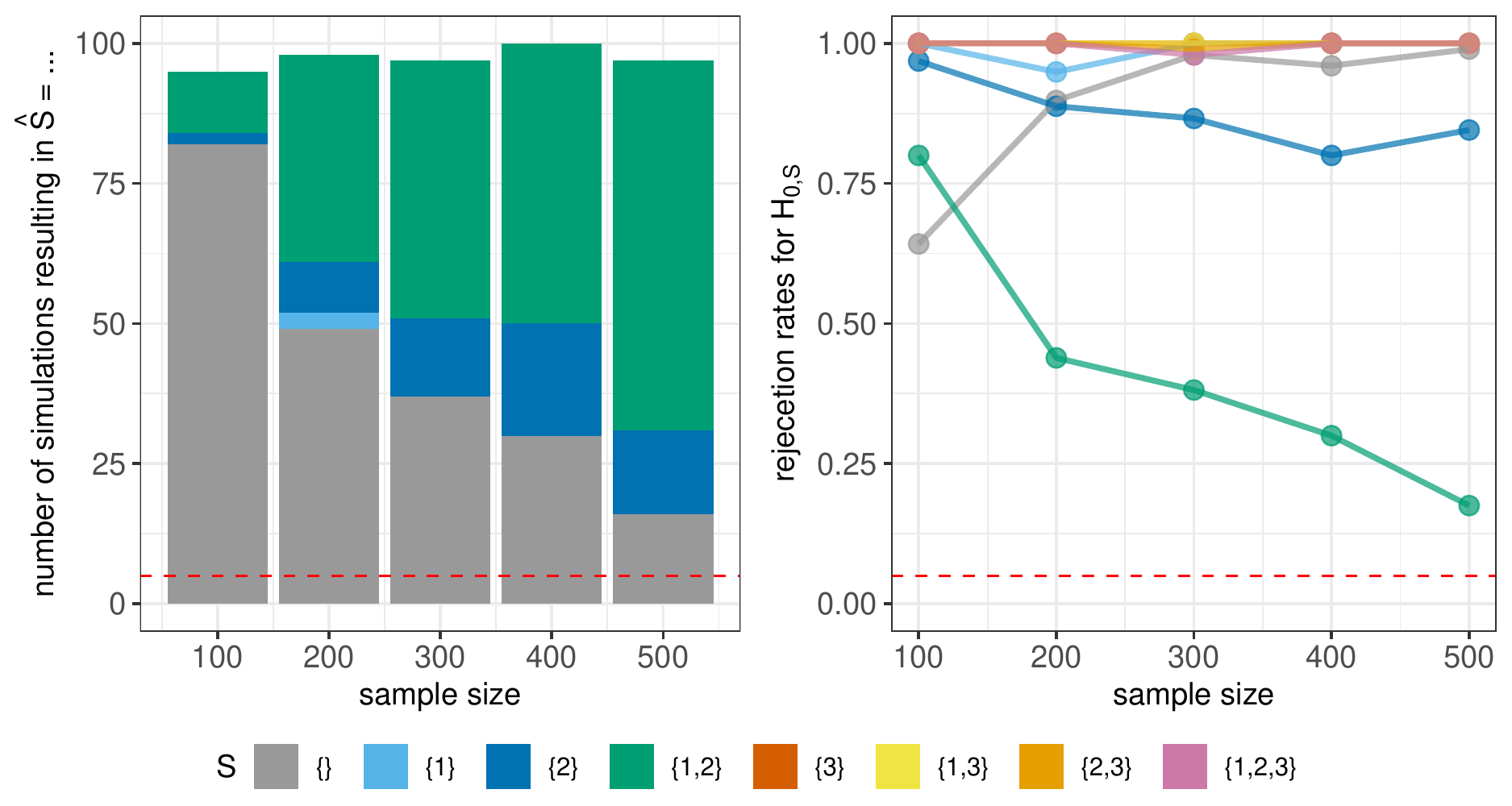}
\end{minipage}
\quad
\begin{minipage}{0.47\textwidth}
\centering
\includegraphics[width = \linewidth]{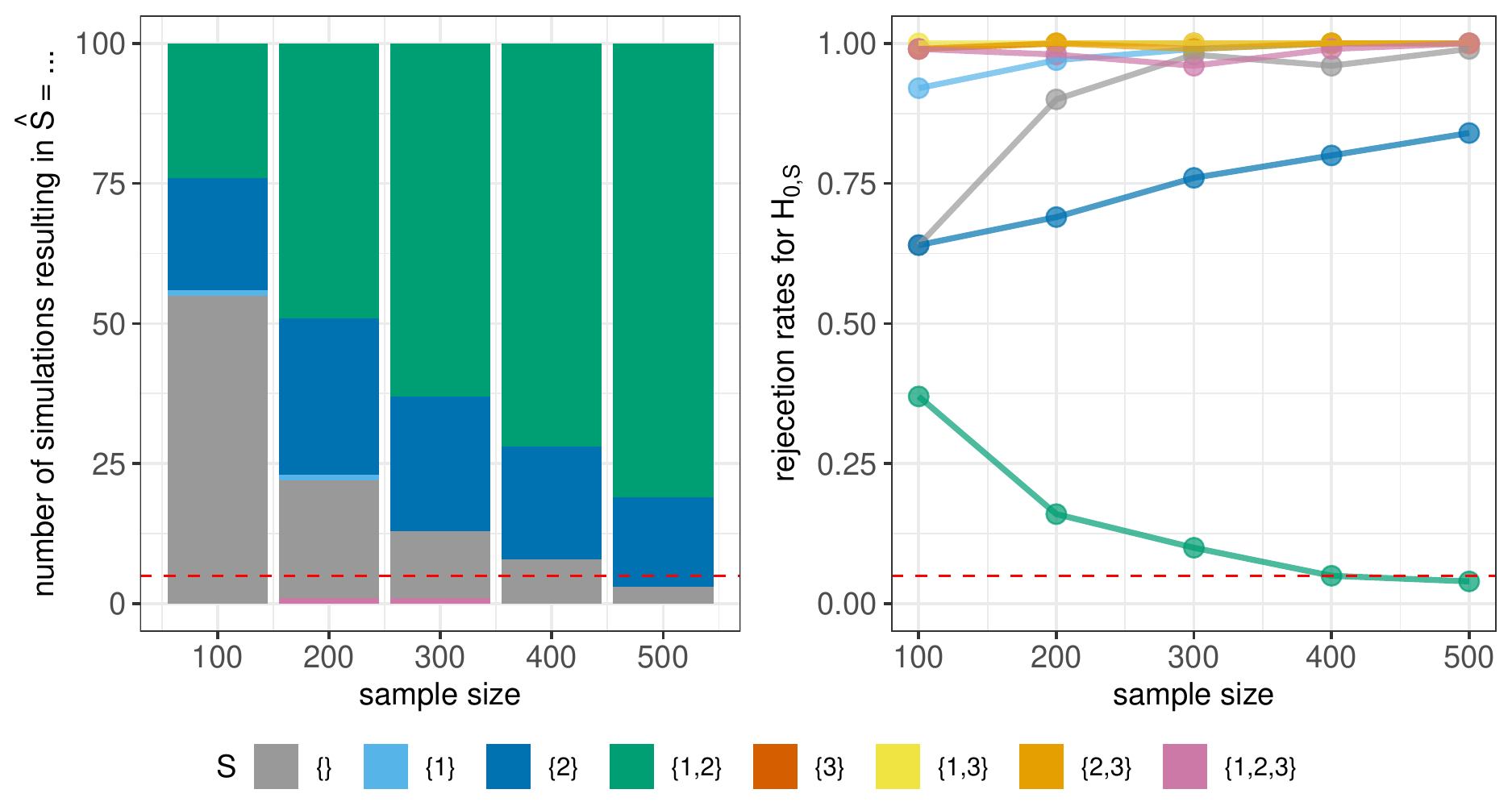}
\end{minipage}
\end{center}

\caption{
Output of ICPH (bar plots) and rejection rates for individual hypotheses (curve plots) for the experiment in Section~\ref{sec:power} with parameter constraint $\sigma_{Y1}^2, \sigma_{Y2}^2 \geq 10^{-16}$ (left) and $\sigma_{Y1}^2 = \sigma_{Y2}^2$ (right), using the EM-algorithm as optimization routine. The results are very similar to those presented in Figure~\ref{fig:shat}, where NLM is applied to the same data. The only notable differences are the missing values in the bar plots (left). These simulations correspond to instances in which the EM-algorithm, after trying several different starting values, failed to converge to a solution which satisfies the variance constraints. 
}
\label{fig:em}
\end{figure}

Figure~\ref{fig:em} shows numerical results for ICPH when using the EM-algorithm as optimization routine. 
The results should be compared to Figure~\ref{fig:shat}, where NLM has been applied to the same data. 
The two methods perform very similarly, 
although NLM is computationally faster (by approximately a factor of 6),
and better suited for handling the lower bound constraint on the error variances. 

\section{Additional Numerical Experiments} \label{app:extrasim}
In this section, we present additional experimental results. In all simulations, we use slight adaptations of the SCM in Section~\ref{sec:SCM}, and measure the performance of ICPH using rejection rates for non-causality (similar to Figure~\ref{fig:sensitivity}). All results are summarized in Figure~\ref{fig:extrasim}.

\subsection{Non-Binary Latent Variables and Unknown Number of States} \label{app:non_binary}
ICPH requires the number of states as an input parameter---we test for $h$-invariance of degree $\ell$ in line 8 of Algorithm~\ref{alg:icph}. If $\ell$ is unknown, we propose the following modification. Let $K \geq 3$ be some predefined integer (e.g., $K=5$), and let for every $S \subset \{1, \dots, d\}$ and every $k \in \{2, \dots, K\}$, $p_S^k$ be a $p$-value for the hypothesis $H_{0,S}^k$ of $h$-invariance of degree $k$ of the set $S$, obtained from the test \eqref{eq:T}. We then substitute the $p$-value $p_S$ in line 8 of Algorithm~\ref{alg:icph} by $p^\prime_S := \max \{p_S^k  \, : \, 2 \leq k \leq K\}$. By construction, the test defined by $p_S^\prime$ is a valid test of $H_{0,S}^\ell$ for any (unknown) $\ell \in \{2, \dots, K\}$. Our code package automatically performs this procedure when the supplied argument
\verb|number.of.states| 
is a vector of length greater than one. 
We now investigate this procedure numerically. %
For a fixed sample size of $n=500$ and for every $\ell \in \{2,3,4,5\}$, we generate $100$ i.i.d.\ data sets from the SCM in Section~\ref{sec:SCM} with parameters sampled as in Section~\ref{sec:level}. The probabilities $\lambda_j = P(H = j)$, $j \in \{0, \dots, \ell\}$ are sampled uniformly between $0.1$ and $1/(\ell+1)$ and standardized correctly. 
In Figure~\ref{fig:extrasim} (left), we 
compare three different approaches: (i) we always test for $h$-invariance of degree 2 (circles), (ii) we always test for $h$-invariance of degree less than or equal to 5, using the approach described above (triangles), and (iii) we test for $h$-invariance using the true number of states $\ell$ (squares). For all methods, ICPH maintains the type I error control, but drops in power as the number of latent states increases. Even if the number of latent states is unknown (but small), ICPH often recovers the causal parents $X^1$ and $X^2$. In general, we propose to limit the application of ICPH to cases where the hidden variables is expected to take only a few different values. 

\subsection{Systems with Large Numbers of Variables}
\begin{figure}
\centering
\hspace{2mm}
\includegraphics[width = .32\linewidth]{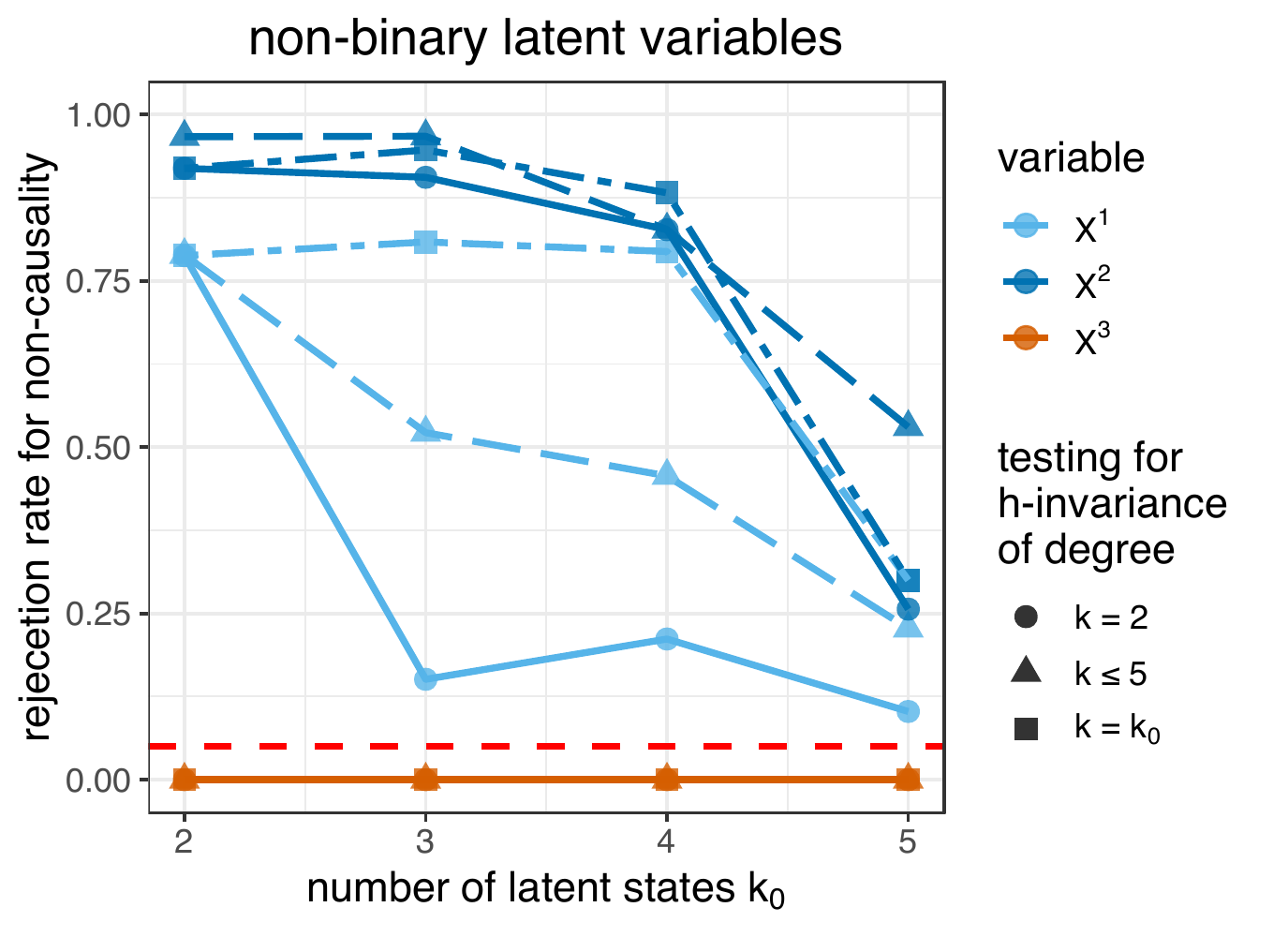}
\includegraphics[width = .32\linewidth]{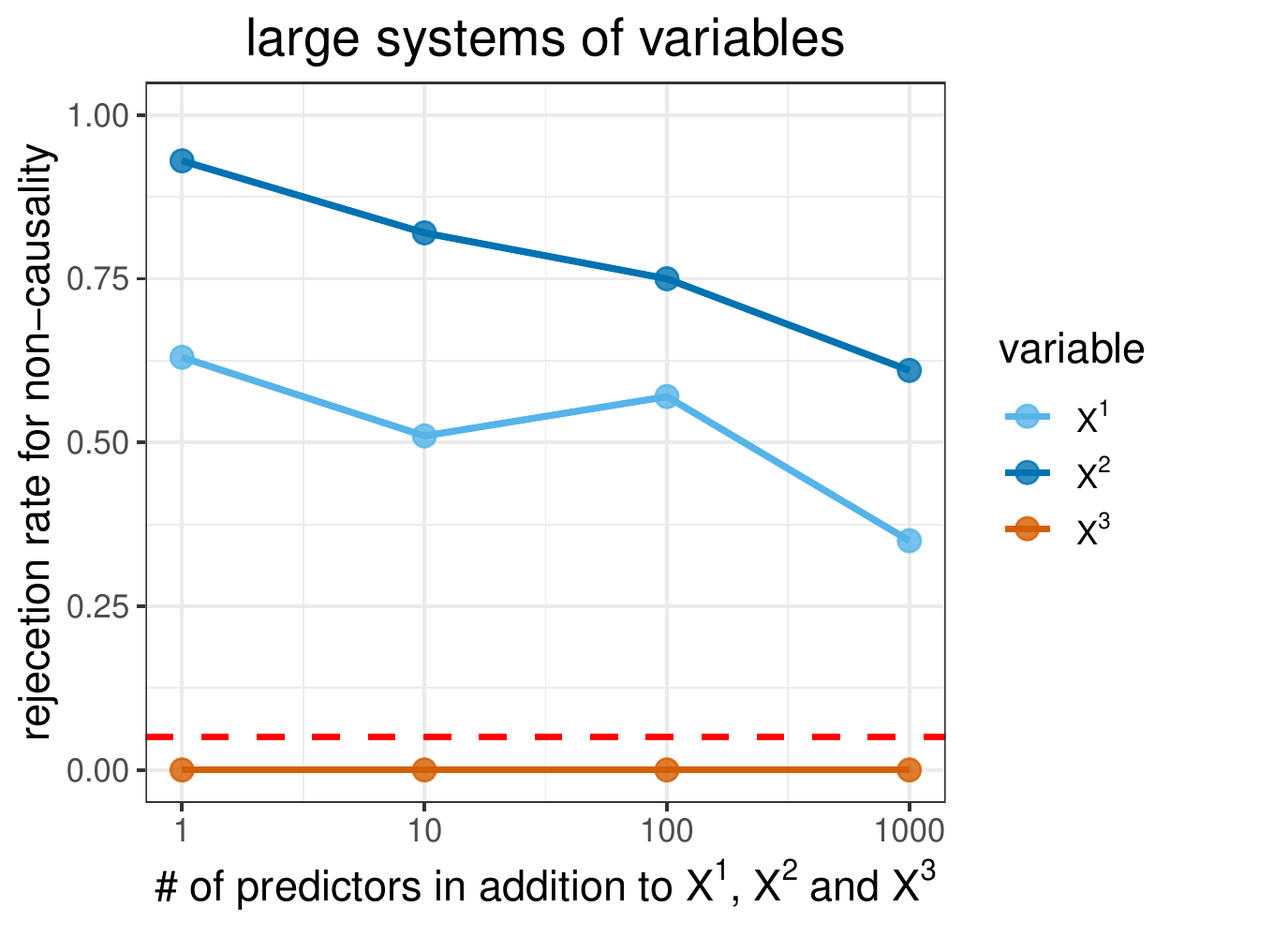}
\includegraphics[width = .32\linewidth]{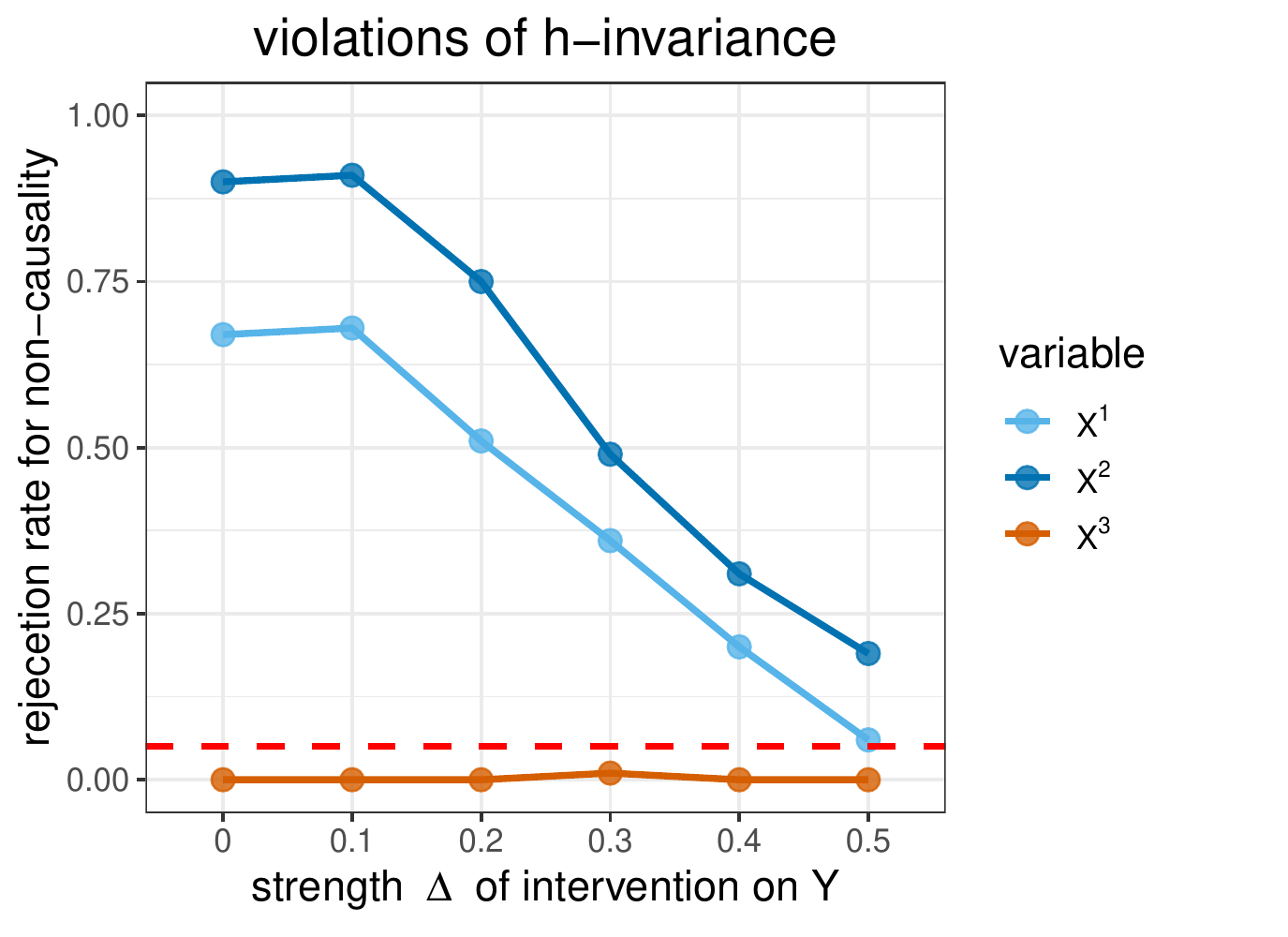}
\caption{
Rejection rates for non-causality of the variables $X^1$, $X^2$ and $X^3$ for the experiments described in Appendix~\ref{app:extrasim}. We investigate the performance of ICPH for non-binary variables (left), for large numbers of predictors (middle), and under violations of the $h$-invariance assumption (right). By simultaneously testing for $h$-invariance of different degrees (see Appendix~\ref{app:non_binary} for details), we can recover $X^1$ and $X^2$ even if the true number of latent states is unknown (left figure, triangles). Our algorithm can be combined with an upfront variable screening (here using Lasso), which results in satisfactory performance even for large number of predictor variables (middle). Under violations of Assumption~\ref{ass:h_inv}, the population version of ICPH is not able to infer $S^* = \{1,2\}$. In the finite sample case we still identify $X^1$ and $X^2$ if $H_{0,S^*}$ is only mildly violated (right). 
}
\label{fig:extrasim}
\end{figure}

For a fixed sample size of $n=300$, we simulate data $(Y, X^1, X^2, X^3, H)$ as described in Section~\ref{sec:SCM}. For increasing $m \in \{1,10,100,1000\}$, we generate additional predictor variables $(Z^1, \dots, Z^m)$ from the structural assignments $Z^j := \alpha_j X^3  + N^Z_j$, $j = 1, \dots, m$, where $N_1^Z, \dots, N_m^Z$ are i.i.d.\ standard Gaussian noise variables, and all $\alpha_j$ are drawn independently from a $\text{Uniform}(-1,1)$ distribution. We then perform variable screening by selecting the first 5 predictors included along the Lasso selection path \citep{Tibshirani94}, and run ICPH on the reduced data set. The results in Figure~\ref{fig:extrasim} (middle) suggest that even for a large number of predictors, ICPH is generally able to infer $S^*$ (provided that $S^*$ contains only few variables).

\subsection{Violations of the $h$-Invariance Assumption} \label{sec:violations_h_inv}
The theoretical guarantees of our method rely on the existence of an $h$-invariant set (Assumption~\ref{ass:h_inv}). 
We now empirically investigate the performance of ICPH under violations of this assumption. 
For a fixed sample size of $n=300$, we generate data as described in Section~\ref{sec:SCM}, 
but include direct interventions on $Y$. For increasing values of $\Delta \in \{0, 0.1, \dots, 0.5\}$, 
we change the coefficients $(\beta_{11}^Y, \beta_{21}^Y)$ in the structural assignment of $Y$ to 
$(\beta_{11}^Y + \Delta, \beta_{21}^Y + \Delta)$ in environment $e_2$, and to
$(\beta_{11}^Y- \Delta, \beta_{21}^Y- \Delta)$
in environment $e_3$. As expected, the 
power of our method drops with the strength of intervention (Figure~\ref{fig:extrasim} right).

\bibliography{ref}

\end{document}